\documentclass[12pt]{amsart}
\usepackage{amssymb}
\usepackage{amsthm}
\usepackage{amscd}

\usepackage{amsmath}
\usepackage{epic,eepic}
\usepackage{graphicx}
\DeclareMathAlphabet{\mathpzc}{OT1}{pzc}{m}{it}

\usepackage{latexsym}

\usepackage{pifont}

\theoremstyle{plain}
\newtheorem{lemma}{Lemma}[subsection]
\newtheorem{prop}[lemma]{Proposition}
\newtheorem{thm}[lemma]{Theorem}
\newtheorem{cor}[lemma]{Corollary}
\newtheorem{aplemma}{Lemma~A.\hspace{-1.5mm}}
\newtheorem{approp}{Proposition~A.\hspace{-1.5mm}}
\newtheorem{apthm}{Theorem~A.\hspace{-1.5mm}}
\newtheorem{apcor}{Corollary~A.\hspace{-1.5mm}}
\newtheorem{intthm}{Theorem}

\usepackage[all]{xy}

\newtheorem{conj}[lemma]{Conjecture}

\theoremstyle{definition}

\newtheorem{rema}[lemma]{Remark}

\newtheorem{remb}{Remark}

\newtheorem{defi}[lemma]{Definition}
\newtheorem{exa}[lemma]{Example}
\newtheorem{aprem}{Remark~A.\hspace{-1.5mm}}
\newtheorem{apdefi}{Definition~A.\hspace{-1.5mm}}
\newcommand{\bde}{\begin{defi}}
\newcommand{\ede}{\end{defi}\vspace{1mm}}
\newcommand{\ble}{\begin{lemma}}
\newcommand{\ele}{\end{lemma}}
\newcommand{\bpr}{\begin{prop}}
\newcommand{\epr}{\end{prop}}
\newcommand{\bt}{\begin{thm}}
\newcommand{\et}{\end{thm}}
\newcommand{\bco}{\begin{cor}}
\newcommand{\eco}{\end{cor}}
\newcommand{\bre}{\begin{rema}}
\newcommand{\ere}{\end{rema}}
\newcommand{\brea}{\begin{rema}}
\newcommand{\erea}{\end{rema}\vspace{1mm}}
\newcommand{\breb}{\begin{remb}}
\newcommand{\ereb}{\end{remb}\vspace{1mm}}
\newcommand{\bex}{\begin{exa}}
\newcommand{\eex}{\end{exa}}
\newcommand{\bpf}{\begin{proof}}
\newcommand{\epf}{\end{proof}\vspace{1mm}}

\newcommand{\bade}{\begin{apdefi}}
\newcommand{\eade}{\end{apdefi}}
\newcommand{\bale}{\begin{aplemma}}
\newcommand{\eale}{\end{aplemma}}
\newcommand{\bapr}{\begin{approp}}
\newcommand{\eapr}{\end{approp}}
\newcommand{\bat}{\begin{apthm}}
\newcommand{\eat}{\end{apthm}}
\newcommand{\baco}{\begin{apcor}}
\newcommand{\eaco}{\end{apcor}}
\newcommand{\bare}{\begin{aprem}}
\newcommand{\eare}{\end{aprem}}

\newcommand{\be}{\begin{enumerate}}
\newcommand{\ee}{\end{enumerate}}
\newcommand{\bcd}{\[\begin{CD}}
\newcommand{\ecd}{\end{CD}\]}
\newcommand{\bit}{\begin{itemize}}
\newcommand{\eit}{\end{itemize}}
\newcommand{\bq}{\begin{quote}}
\newcommand{\eq}{\end{quote}}
\newcommand{\ba}{\begin{array}}
\newcommand{\ea}{\end{array}}

\newcommand{\mcE}{\mathcal{E}}
\newcommand{\mcF}{\mathcal{F}}
\newcommand{\mcG}{\mathcal{G}}

\newcommand{\mcK}{\mathcal{K}}
\newcommand{\mcL}{\mathcal{L}}
\newcommand{\mcM}{\mathcal{M}}
\newcommand{\mcN}{\mathcal{N}}
\newcommand{\mcO}{\mathcal{O}}

\newcommand{\mcS}{\mathcal{S}}
\newcommand{\mcT}{\mathcal{T}}

\newcommand{\mcV}{\mathcal{V}}

\newcommand{\mbC}{\mathbb{C}}

\newcommand{\mbF}{\mathbb{F}}

\newcommand{\mbH}{\mathbb{H}}

\newcommand{\mbR}{\mathbb{R}}

\newcommand{\mbZ}{\mathbb{Z}}

\newcommand{\mfU}{\mathfrak{U}}

\newcommand{\mfb}{\mathfrak{b}}

\newcommand{\mfg}{\mathfrak{g}}

\newcommand{\mfl}{\mathfrak{l}}

\newcommand{\mfn}{\mathfrak{n}}

\newcommand{\mfs}{\mathfrak{s}}

\newcommand{\migi}{\rightarrow}
\newcommand{\longmigi}{\longrightarrow}

\newcommand{\isom}{\stackrel{\sim}{\migi}}

\newcommand{\migiincl}{\hookrightarrow}

\newcommand{\migisurj}{\twoheadrightarrow}

\newcommand{\mr}{\mathrm}
\newcommand{\hidden}[1]{\,}

\begin{document}

\title[Symplectic geometry of $p$-adic Teichm\"{u}ller uniformization]{Symplectic geometry   \\  of $p$-adic Teichm\"{u}ller uniformization \\ for ordinary nilpotent indigenous bundles}
\author{Yasuhiro Wakabayashi}
\date{}
\markboth{Yasuhiro Wakabayashi}{}
\maketitle
\footnotetext{Y. Wakabayashi: 
Department of Mathematics, Tokyo Institute of Technology, 2-12-1 Ookayama, Meguro-ku, Tokyo 152-8551, JAPAN;}
\footnotetext{e-mail: {\tt wkbysh@math.titech.ac.jp};}
\footnotetext{2010 {\it Mathematical Subject Classification}: Primary 14H10, Secondary 53D30;}
\footnotetext{Key words: $p$-adic Teichm\"{u}ller theory, hyperbolic curve, indigenous bundle, symplectic structure, canonical lifting, uniformization, crystal.}
\begin{abstract}
The aim of the present paper is to provide a new aspect of the $p$-adic Teichm\"{u}ller theory established by S. Mochizuki. We study the symplectic geometry of the $p$-adic formal stacks $\widehat{\mathcal{M}}_{g, \mathbb{Z}_p}$ (= the moduli classifying $p$-adic formal curves of fixed genus $g>1$) and $\widehat{\mathcal{S}}_{g, \mathbb{Z}_p}$ (= the moduli classifying $p$-adic formal curves of genus $g$ equipped with an indigenous bundle). A major achievement in the (classical) $p$-adic Teichm\"{u}ller theory is the construction of the locus $\widehat{\mathcal{N}}_{g, \mathbb{Z}_p}^{\mathrm{ord}}$ in $\widehat{\mathcal{S}}_{g, \mathbb{Z}_p}$ classifying $p$-adic canonical liftings of ordinary nilpotent indigenous bundles. The formal stack $\widehat{\mathcal{N}}_{g, \mathbb{Z}_p}^{\mathrm{ord}}$ embodies a $p$-adic analogue of uniformization of hyperbolic Riemann surfaces, as well as a hyperbolic analogue of Serre-Tate theory of ordinary abelian varieties. In the present paper, the canonical symplectic structure on the cotangent bundle $T^\vee_{\mathbb{Z}_p} \widehat{\mathcal{M}}_{g, \mathbb{Z}_p}$ of $\widehat{\mathcal{M}}_{g, \mathbb{Z}_p}$ is compared to Goldman's symplectic structure defined on $\widehat{\mathcal{S}}_{g, \mathbb{Z}_p}$ after base-change by the projection $\widehat{\mathcal{N}}_{g, \mathbb{Z}_p}^{\mathrm{ord}} \rightarrow \widehat{\mcM}_{g, \mbZ_p}$. We can think of this comparison as a $p$-adic analogue of certain results in the theory of projective structures on Riemann surfaces proved by S. Kawai and other mathematicians. The key ingredient in our discussion is the $F$-crystal structure on the de Rham/crystalline cohomology associated to the adjoint bundle of each ordinary nilpotent indigenous bundle. We show that the slope decomposition of this  $F$-crystal has a geometric interpretation, i.e., arises as the differential of the $p$-adic Teichm\"{u}ller uniformization. This fact  makes it clear how the two symplectic structures are related.
\end{abstract}
\tableofcontents 

\section*{Introduction}

 The aim of the present paper is to provide a new aspect of the {\it $p$-adic Teichm\"{u}ller theory} established  by S. Mochizuki (cf. ~\cite{Mzk1}, ~\cite{Mzk2}).
We study the symplectic geometry of 
the $p$-adic formal stack defined as the moduli  classifying  $p$-adic formal   curves equipped with an 
 indigenous bundle.
As a consequence of the present paper, we will propose
 a $p$-adic analogue  of the comparison assertion 
 concerning  a canonical symplectic structure
 on the moduli space of projective structures.
 
 Here, recall that a {\it projective structure} on a Riemann surface is a holomorphic  atlas whose transition functions between coordinate charts may be expressed as  M\"{o}bius transformations.
 Also, an {\it indigenous bundle}
   is defined as a sort of algebro-geometric counterpart of a projective structure (cf. \S\,\ref{SS0199} for the precise definition of an indigenous bundle).
 Indigenous bundles, or equivalently, projective structures, 
   have 
    provided a rich  story  in complex (i.e., the usual) Teichm\"{u}ller theory for a long time.
 Canonical examples 
  are constructed by means of various uniformizations such as Fuchsian, 
  Bers,
  or Schottky.
In other words, one may think of an indigenous bundle as an algebraic object encoding (analytic, i.e., non-algebraic) uniformization data for Riemann surfaces.
 
 One subject in the theory of projective structures is to compare  symplectic structures on the relevant spaces via uniformization.
 For example, we shall refer to  works by  S. Kawai, P. Ar\'{e}s-Gastesi, I. Biswas,  B. Loustau,  et al. (cf. ~\cite{Ka}; ~\cite{AGBi}; ~\cite{AGBi2}; ~\cite{Los}). 
To explain some of these works, let us consider
the following spaces associated to  a connected orientable closed surface $\Sigma$ of genus $g>1$:
\begin{align}
 \mcT^\mr{an}_{g, \mbC} := \mr{Conf}(\Sigma)/\mr{Diff}^0(\Sigma) \ \left(\text{resp.,} \    \mcS^\mr{an}_{g, \mbC} := \mr{Proj}(\Sigma)/\mr{Diff}^0(\Sigma) \right),
 \end{align}
where $\mr{Conf}(\Sigma)$ (resp., $\mr{Proj}(\Sigma)$) denotes the space of all 
holomorphic 
structures (resp., all projective structures) on $\Sigma$ compatible with the orientation of $\Sigma$, and $\mr{Diff}^0(\Sigma)$ denotes  the group of all diffeomorphisms of $\Sigma$ homotopic to the identity map of $\Sigma$. 
(In particular, $ \mcT_{g, \mbC}^\mr{an}$ is nothing but the  {\it  Teichm\"{u}ller space}  associated to $\Sigma$.)
It is well-known that $\mcT_{g, \mbC}^\mr{an}$  
admits  the  structure of 
 a complex manifold of dimension $3g-3$ which is 
a universal covering  of the moduli space $\mcM_{g,\mbC}^\mr{an}$ classifying  connected compact Riemann surfaces of genus $g$.
Also,
 $\mcS_{g, \mbC}^\mr{an}$
admits  the  structure of  a complex manifold of dimension $6g-6$ and moreover
   the structure of   a relative  affine space over  $\mcT_{g, \mbC}^\mr{an}$ modeled on (the total space of) the holomorphic cotangent bundle
  $T^\vee_\mbC \mcT_{g, \mbC}^\mr{an}$ of
   $\mcT_{g, \mbC}^\mr{an}$.

Now, let us take  a 
holomorphic  section
\begin{align} 
\sigma  :  \mcT_{g, \mbC}^\mr{an} \migi  \mcS_{g, \mbC}^\mr{an}
 \end{align}
 of the natural projection 
 $\mcS_{g, \mbC}^\mr{an} \migisurj \mcT_{g, \mbC}^\mr{an}$.
Because of the affine structure on $\mcS_{g, \mbC}^\mr{an}$,
this section  may be  extended to   a unique  
biholomorphism
  $\theta_\sigma : T^\vee_\mbC \mcT_{g, \mbC}^\mr{ad} \isom \mcS_{g, \mbC}^\mr{an}$
compatible with the respective affine structures.
It  induces an isomorphism
\begin{align}
\Theta_\sigma : 
\theta_\sigma^*(\bigwedge^2 \Omega_{\mcS^\mr{an}_{g, \mbC}}) \ \left(\cong \bigwedge^2 \theta_\sigma^* (\Omega_{\mcS^\mr{an}_{g, \mbC}}) \right) \isom \bigwedge^2 \Omega_{T^\vee_\mbC \mcT_{g, \mbC}^\mr{an}}.
\end{align}
Notice  that
$T^\vee_\mbC \mcT^\mr{an}_{g, \mbC}$ admits a canonical holomorphic symplectic structure $\omega^\mr{Liou}_{g, \mbC}$ obtained as the differential of the tautological $1$-form (i.e., the so-called Liouville form).
Moreover, 
$\mcS_{g, \mbC}^\mr{an}$ admits a holomorphic symplectic structure $\omega^\mr{PGL}_{g, \mbC}$ induced,  via pull-back by the monodromy map,  from   Goldman's symplectic structure on the $\mr{PGL}_2 (\mbC)$-character variety (cf. ~\cite{Go}).
Thus, we obtain holomorphic symplectic manifolds
\begin{align} \label{W10}
(T^\vee_\mbC \mcT^\mr{an}_{g, \mbC}, \omega^\mr{Liou}_{g, \mbC})
\hspace{5mm} \text{and} \hspace{5mm}
(\mcS_{g, \mbC}^\mr{an}, \omega^\mr{PGL}_{g, \mbC}).
\end{align}
According to a pioneering result proved by S. Kawai,
we can compare these symplectic structures. 
In fact, it follows from ~\cite{Ka}, Theorem, that
if $\sigma$
is any Bers section,
  then
$\theta_\sigma$
preserves the symplectic structure  up to a constant factor; more precisely, the following equality holds:
\begin{align} \label{e06}
 \Theta_{\sigma}(\omega_{g, \mbC}^{\mr{PGL}}) = \pi
 \cdot \omega^{\mr{Liou}}_{g, \mbC}.
 \end{align}
Also,  B. Loustau proved  (cf. ~\cite{Los}, Theorem 6.10) this equality, 
  which may be  described as the equality
$\Theta_{\sigma}(\omega_{g, \mbC}^{\mr{PGL}}) =
 \sqrt{-1} 
 \cdot \omega^{\mr{Liou}}_{g, \mbC}$
with the conventions chosen by him.
Moreover, 
by ~\cite{AGBi}, Theorem 1.1 and  Remark 3.2,
 the same equality holds for the case where  $\sigma$ is taken as a section arising from either the Schottky uniformizaton  or the Earle uniformization.


This article attempts to consider a $p$-adic analogue of these comparison results.
Let $p$ be an odd prime and $\widehat{\mcM}_{g, \mbZ_p}$ (cf. (\ref{W9}))
  denote 
 the $p$-adic formal stack defined as the moduli  classifying (proper, smooth, and geometrically connected) $p$-adic formal  curves of fixed genus $g>1$.
 One obtains a $p$-adic formal stack
  $\widehat{\mcS}_{g, \mbZ_p}$ (cf. (\ref{W5}))
   defined as the moduli classifying  
  pairs of such
  a curve and an indigenous bundle on it.
By the same manner as  the complex case discussed above,
$\widehat{\mcS}_{g, \mbZ_p}$ and 
 the cotangent bundle $T^\vee_{\mbZ_p} \widehat{\mcM}_{g, \mbZ_p}$ of $\widehat{\mcM}_{g, \mbZ_p}$  admit canonical structures
 \begin{align}
\widehat{\omega}^\mr{Liou}_{g, \mbZ_p} \in \Gamma (T^\vee_{\mbZ_p}\widehat{\mcM}_{g, \mbZ_p}, \bigwedge^2 \Omega_{T^\vee_{\mbZ_p}\widehat{\mcM}_{g, \mbZ_p}/\mbZ_p}),
 \hspace{5mm}
\widehat{\omega}^\mr{PGL}_{g, \mbZ_p} \in \Gamma (\widehat{\mcS}_{g, \mbZ_p}, \bigwedge^2 \Omega_{\widehat{\mcS}_{g, \mbZ_p}/\mbZ_p})
\end{align}
respectively (cf. (\ref{W8}) and (\ref{W7})). 
 $\widehat{\omega}^\mr{PGL}_{g, \mbZ_p}$ is, by definition, obtained by composing
 the Killing form on $\mfs \mfl_2$ and the cup product in the de Rham cohomology of the adjoint bundles on indigenous bundles.
In this way, we obtain two symplectic $p$-adic formal  stacks 
\begin{align} \label{W17}
(T^\vee_{\mbZ_p}\widehat{\mcM}_{g, \mbZ_p}, \widehat{\omega}^\mr{Liou}_{g, \mbZ_p})
\hspace{5mm} \text{and} \hspace{5mm}
(\widehat{\mcS}_{g, \mbZ_p}, \widehat{\omega}^\mr{PGL}_{g, \mbZ_p}),
\end{align}
 which have the natural projections onto $\widehat{\mcM}_{g, \mbZ_p}$.


Here, we shall recall the main achievement of the (classical) $p$-adic Techm\"{u}ller theory studied in ~\cite{Mzk1}.
Denote by $\mcN^\mr{ord}_{g, \mbF_p}$ (cf. (\ref{W1})) the locus in the stack $\mcS_{g, \mbF_p} := \widehat{\mcS}_{g, \mbZ_p} \otimes \mbF_p$ over $\mbF_p :=\mbZ/p\mbZ$ classifying {\it ordinary nilpotent indigenous bundles} (cf. \S\,\ref{SS02} for the precise definition of an ordinary nilpotent indigenous bundle). 
 S. Mochizuki proved (cf. ~\cite{Mzk1}, Chap.\,II, \S\,3, Corollary 3.8) that  $\mcN^\mr{ord}_{g, \mbF_p}$ is a nonempty Deligne-Mumford stack which is \'{e}tale and dominant over $\mcM_{g, \mbF_p} := \widehat{\mcM}_{g, \mbZ_p} \otimes \mbF_p$.
This implies that  there exists  
 the unique (up to isomorphism) $p$-adic formal stack
 \begin{align}
\widehat{ \mcN}^\mr{ord}_{g, \mbZ_p}
 \end{align}
(cf. (\ref{W2}))  over $\widehat{\mcM}_{g, \mbZ_p}$ lifting $\mcN^\mr{ord}_{g, \mbF_p}$.
Moreover, he also constructed  (cf. Theorem \ref{T0135} for the precise statement)
 a canonical $p$-adic lifting  $\Phi_\mcN : \widehat{ \mcN}^\mr{ord}_{g, \mbZ_p} \migi \widehat{ \mcN}^\mr{ord}_{g, \mbZ_p}$ of the Frobenius endomorphism of $\mcN^\mr{ord}_{g, \mbF_p}$ together with  an indigenous bundle $(\mcE_\mcN, \nabla_{\mcE_\mcN})$  on the universal family of curves $C_{\mcN}$ over $\widehat{ \mcN}^\mr{ord}_{g, \mbZ_p}$ which is Frobenius invariant in the sense that $\mbF^*(\Phi^*_\mcN (\mcE_\mcN, \nabla_{\mcE_\mcN})) \cong (\mcE_\mcN, \nabla_{\mcE_\mcN})$ (where $\mbF^*(-)$ denotes  renormalized Frobenius pull-back in the sense of  \S\,\ref{SS179}).
This result is used, via restriction to 
various 
 points in $\mcN_{g, \mbF_p}^\mr{ord}$, to
 obtain $p$-adic canonical liftings of curves (endowed  with an ordinary nilpotent indigenous bundle).
We shall refer to $\widehat{ \mcN}^\mr{ord}_{g, \mbZ_p}$ together with both  $\Phi_\mcN$ and $(\mcE_\mcN, \nabla_{\mcE_\mcN})$ as the {\bf classical ordinary $p$-adic Teichm\"{u}ller}
 {\bf uniformization}.


These $p$-adic objects create a situation similar to the complex case.
This means that the indigenous bundle $(\mcE_\mcN, \nabla_{\mcE_\mcN})$ determines its classifying 
 morphism $\sigma : \widehat{ \mcN}^\mr{ord}_{g, \mbZ_p} \migiincl \widehat{\mcS}_{g, \mbZ_p}$ (i.e., a section of the projection $\widehat{\mcS}_{g, \mbZ_p} \migi \widehat{\mcM}_{g, \mbZ_p}$ defined on  the \'{e}tale dominant formal stack $\widehat{ \mcN}^\mr{ord}_{g, \mbZ_p}$ over $\widehat{\mcM}_{g, \mbZ_p}$).
Moreover, since $\widehat{\mcS}_{g, \mbZ_p}$ forms a relative  affine space over $\widehat{\mcM}_{g, \mbZ_p}$ modeled on $T^\vee_{\mbZ_p} \widehat{\mcM}_{g, \mbZ_p}$,
the morphism $\sigma$ gives a trivialization
\begin{align}
\theta : 
T^\vee_{\mbZ_p} \widehat{\mcM}_{g, \mbZ_p} |_\mcN
 \ \left(:= T^\vee_{\mbZ_p} \widehat{\mcM}_{g, \mbZ_p} \times_{\widehat{\mcM}_{g, \mbZ_p}} \widehat{\mcN}^\mr{ord}_{g, \mbZ_p} \right) 
 \isom  \widehat{\mcS}_{g, \mbZ_p} |_\mcN
 \ \left(:= \widehat{\mcS}_{g, \mbZ_p} \times_{\widehat{\mcM}_{g, \mbZ_p}}  \widehat{\mcN}^\mr{ord}_{g, \mbZ_p} \right)
\end{align}
(cf. (\ref{e58}))
of $\widehat{\mcS}_{g, \mbZ_p}$ after  base-change to  $\widehat{\mcN}^\mr{ord}_{g, \mbZ_p}$.
This trivialization  induces an isomorphism
\begin{align}
\Theta : 
\theta^* (\bigwedge^2 \Omega_{\widehat{\mcS}_{g, \mbZ_p}|_\mcN/\mbZ_p}) 
\ \left(\cong \wedge^2 \theta^*(\Omega_{\widehat{\mcS}_{g, \mbZ_p}|_\mcN/\mbZ_p}) \right)
\isom 
\bigwedge^2 \Omega_{T^\vee_{\mbZ_p}\widehat{\mcM}_{g, \mbZ_p} |_\mcN /\mbZ_p}.
\end{align}
Then, the main result of the present paper is described as  the  following  theorem,  asserting  the comparison, via $\Theta$,    between 
the pull-backs  of symplectic structures
\begin{align}
\widehat{\omega}^\mr{Liou}_{g, \mbZ_p} |_\mcN 
  := \widehat{\omega}^\mr{Liou}_{g, \mbZ_p} |_{T^\vee_{\mbZ_p}\widehat{\mcM}_{g, \mbZ_p}|_\mcN}, \hspace{10mm}
 \widehat{\omega}^\mr{PGL}_{g, \mbZ_p} |_\mcN
  := \widehat{\omega}^\mr{PGL}_{g, \mbZ_p} |_{\widehat{\mcS}_{g, \mbZ_p}|_\mcN}
\end{align}
defined on  $T^\vee_{\mbZ_p}\widehat{\mcM}_{g, \mbZ_p}|_\mcN$, 
$\widehat{\mcS}_{g, \mbZ_p}|_\mcN$ respectively.
(We shall refer to ~\cite{Wak3} for the version of this theorem in the case of the moduli classifying {\it dormant} indigenous bundles $\mcM^{^\mr{Zzz...}}_{g, \mbF_p}$.)

\vspace{3mm}
\begin{intthm} [= Theorem \ref{T013}]
\leavevmode\\
 \ \ \ 
If $p>3$, then
the morphism $\theta$ preserves the symplectic structure, i.e., the  following equality holds:
\begin{align}
\Theta(\widehat{\omega}_{g, \mbZ_p}^\mr{PGL} |_{\mcN}) = \widehat{\omega}_{g, \mbZ_p}^\mr{Liou} |_\mcN.
\end{align}
 In particular, the image of $\sigma :  \widehat{\mcN}^\mr{ord}_{g, \mbZ_p} \migiincl  \widehat{\mcS}_{g, \mbZ_p}$ is Lagrangian with respect to the symplectic structure $\widehat{\omega}_{g, \mbZ_p}^\mr{PGL}$.
 \end{intthm}
\vspace{5mm}

In this paragraph, we shall make a brief comment on the proof of the above theorem.
The key ingredient  in our discussion 
  is  the $F$-crystal structure of the  cohomology  associated to  the adjoint bundle of each ordinary   nilpotent indigenous bundle.
Indeed, let us take  
a collection $(X_1, \mcE_1, \nabla_{\mcE, 1})$ classified by $\mcN^\mr{ord}_{g, \mbZ_p}$.
 Denote by $(X, \mcE, \nabla_\mcE)$ the canonical $p$-adic lifting  of $(X_1, \mcE_1, \nabla_{\mcE, 1})$  arising from the c.\,o.\,$p$-Teich.\,uniformization 
and  by $s_\infty$ its classifying  point of $\widehat{\mcN}^\mr{ord}_{g, \mbZ_p}$.
Then, the
 Frobenius invariance of $(\mcE, \nabla_\mcE)$ gives rise to
 an $F$-crystal structure  (cf. (\ref{G951}))
   on the  first  de Rham cohomology $\mbH^1 (\mcK^\bullet [\nabla^\mr{ad}_\mcE])$ (which is isomorphic to the first  crystalline cohomology 
    of the corresponding crystal) of the adjoint flat  bundle associated to $(\mcE, \nabla_\mcE)$.
On the other hand,
since 
 $\mbH^1 (\mcK^\bullet [\nabla^\mr{ad}_\mcE])$ is canonically isomorphic to
 the tangent space of $\widehat{\mcS}_{g, \mbZ_p}$
 at $s_\infty$ (cf. (\ref{EE001})),
 the differential of the embedding  $\sigma : \widehat{\mcN}^\mr{ord}_{g, \mbZ_p} \migi \widehat{\mcS}_{g, \mbZ_p}$ determines
 a direct sum decomposition of $\mbH^1 (\mcK^\bullet [\nabla^\mr{ad}_\mcE])$ (cf. (\ref{G680})).
One important observation is (cf. Corollary \ref{C0338}) that
this decomposition of the $F$-crystal  coincides with (i.e,  gives the geometric interpretation of) the slope decomposition.
It follows (cf. Corollary \ref{G970})  that  both  $\widehat{\omega}^\mr{Liou}_{g, \mbZ_p} |_\mcN$ and $\widehat{\omega}^\mr{PGL}_{g, \mbZ_p} |_\mcN$ turn out to specify  eigenvectors of the $F$-crystal structure defined on 
the second exterior power of the dual   $\mbH^1 (\mcK^\bullet [\nabla^\mr{ad}_\mcE])^\vee$.
This fact makes it clear how the two symplectic structures 
 are related via reduction modulo $p$.
Thus, the proof of the main theorem will be reduced to  an explicit computation of $\mbH^1 (\mcK^\bullet [\nabla^\mr{ad}_\mcE])$ in terms of the \v{C}ech double complex (cf. the discussion in \S\,\ref{SS699}).

In the Appendix of the present paper, we discuss crystals of  torsors (equipped with a structure group) and prove the correspondence between flat torsors and them (cf. Theorem \ref{T043}).
This correspondence  may be thought of as a generalization of the  classical result (cf. ~\cite{BO}, \S\,6.6,  Theorem) for  crystals of vector bundles (i.e., of  $\mr{GL}_n$-torsors).
Moreover, we observe  (cf. Proposition \ref{W22}) the relationship between  the respective  deformations  of a prescribed flat torsor over 
distinct underlying spaces.
Its  application 
 to the case
 of indigenous bundles
(cf. Proposition \ref{W23}) will be used in the proof of the main theorem.
Throughout  the present paper, we shall often  refer to   the Appendix for  some definitions and facts involved.

\vspace{5mm}
\hspace{-4mm}{\bf Acknowledgement} \leavevmode\\
 \ \ \ 
The author  would like to express my sincere gratitude to
Professors Yuichiro Taguchi, Shinichi Mochizuki,  and  Shingo Kawai
(and  various  moduli spaces, e.g.,  $\mcM_{g, \mbF_p}$ and  $\mcN^\mr{ord}_{g, \mbF_p}$)
 for their helpful advice and heartfelt encouragement.
The author was partially  supported by the Grant-in-Aid for Scientific Research (KAKENHI No.\,18K13385).

\vspace{10mm}
\section{Symplectic structures on the moduli of indigenous bundles} \vspace{3mm}

In this section, we shall review various notions and facts concerning our discussion.
In particular, the central characters of the present paper, i.e., the $p$-adic formal stacks and their symplectic structures displayed  in (\ref{W17}) are defined precisely.
Throughout the present paper, we fix an odd prime $p$ and an integer $g$ with $g >1$.

\vspace{5mm}
\subsection{Symplectic structures} \label{SS001}
\leavevmode\\ \vspace{-4mm}

We begin by  reviewing the notion of a symplectic structure.
Let $R$ be a commutative ring with unit and 
$X$  a smooth Deligne-Mumford stack over $R$ of relative dimension $n >0$.
Denote by $\Omega_{X/R}$ the sheaf of $1$-forms on $X$ relative to $R$ and by $\mcT_{X/R}$ its dual.
Hence, both $\Omega_{X/R}$ and $\mcT_{X/R}$ are   vector bundles (i.e., locally free coherent sheaves) on $X$ of rank $n$.
A {\bf symplectic structure}  on $X$  is, by definition,  a nondegenerate closed $2$-form $\omega \in \Gamma (X, \bigwedge^2\Omega_{X/R})$.
Here, we shall   say that a $2$-form $\omega$ is {\bf nondegenerate}  if the $\mcO_X$-linear  morphism $\Omega_{X/R} \migi \mcT_{X/R}$
 induced naturally by $\omega$ is an isomorphism.
Let us fix a symplectic structure $\omega$ on $X$.
Then, we shall say that a smooth substack $Y$ of $X$ is {\bf Lagrangian} (with respect to $\omega$)
if $\omega |_Y =0$ and  $\mr{dim} (Y) = \frac{n}{2}$.

Let us   recall the canonical symplectic structure defined on  the cotangent space. 
Let $X$ be as above.
Denote  the total space of $\Omega_{X/R}$ by
\begin{align}
T^\vee_R X,
\end{align}
and call it
 the {\bf cotangent bundle} of $X$ (over $R$).
One obtains a  symplectic structure
\begin{align} \label{e10}
   \omega^{\mr{Liou}}_{X}  \in \Gamma (T^\vee_R X, {\bigwedge}^2\Omega_{T^\vee_R X/R})
   \end{align}
 on $T^\vee_R X$
 defined as the differential of the tautological $1$-form (called the Liouville form)  on $T^\vee_R X$;
it may be 
   characterized  uniquely  by the condition that
if $q_1, \cdots, q_n$ are any local coordinates in $X$  (relative to $R$) and $p_1, \cdots,  p_n$  are   the dual coordinates 
in $T^\vee_R X$, then 
 $\omega^{\mr{Liou}}_{X}$  has the local expression
$\sum_{i=1}^n dp_i \wedge dq_i$.

Next, denote by
\begin{align}  \label{W25}
0_X : X \migi T^\vee_R X 
\end{align}
the zero section, whose image is immediately  verified to be Lagrangian.
The 
pull-back 
  $0_X^*(\Omega_{T^\vee_R X/X})$ of $\Omega_{T^\vee_R X/X}$ is canonically isomorphic to $\Omega_{X/R}$.
In what follows, let us describe the $\mcO_X$-bilinear map on $0^*_X(\mcT_{T^\vee_R X/R}) \ \left(= 0_X^*(\Omega_{T^\vee_R X/R})^\vee\right)$ corresponding to the restriction $\omega_X^\text{Liou}|_{0_X}$ of  $\omega_X^\text{Liou}$.
Consider the short exact sequence
\begin{align} \label{e732}
   0 \longmigi  0_X^*(\mcT_{T^\vee_R X/X}) \ \left(\cong \Omega_{X/R} \right)\longmigi 0^*_X(\mcT_{T^\vee_R X/R}) \longmigi \mcT_{X/R} \longmigi 0
   \end{align}
   obtained by differentiating 
   the projection $ T^\vee_R X \migi X$
    and successively restricting it to $0_X$.
The differential of  $0_X : X \migi T^\vee_R X$ specifies 
a split injection $\mcT_{X/R} \migiincl 0^*_X(\mcT_{T^\vee_R X/R})$ of this  short exact sequence.
In other words, $0_X$ gives
 a direct sum decomposition
\begin{align} \label{e741} 
0^*_X(\mcT_{T^\vee_R X/R}) 
\isom   \mcT_{X/R} \oplus  \Omega_{X/R}.
\end{align}
Then, the  $\mcO_X$-bilinear map on $0^*_X(\mcT_{T^\vee_R X/R})$
corresponding to $\omega_X^\text{can}|_{0_X}$ is 
given
by the natural  pairing
$\langle-, -\rangle : \mcT_{X/R} \times \Omega_{X/R} \migi \mcO_X$.
More precisely, this bilinear map may be expressed, via
(\ref{e741}),
  as the map given by assigning
\begin{align} \label{e75}
     ((a, b), (a',b')) \mapsto \langle a, b'\rangle - \langle a', b\rangle  \end{align}
 for local sections $a$, $a' \in  \mcT_{X/R}$ and $b$, $b' \in \Omega_{X/R}$.

Let us consider  the case of formal stacks.
Let  $X$ a smooth $p$-adic  formal stack over $\mbZ_p$;  it may be given as 
$X = \varinjlim_n X_n$, where each $X_n$ ($n \geq 1$) is a smooth stack over $\mbZ/p^n \mbZ$ such that  $X_n = X_m \otimes_{\mbZ / p^m \mbZ} (\mbZ/p^n \mbZ)$ (if $n <m$).
We shall write $\mcT_{X/\mbZ_p} := \varprojlim_n \mcT_{X_n/(\mbZ/p^n \mbZ)}$ and 
$\Omega_{X/\mbZ_p} := \varprojlim_n \Omega_{X_n/(\mbZ/p^n \mbZ)}$, that are 
rank $n$ vector bundles on $X$.
By a {\bf symplectic structure} on $X$, we mean a  collection   $\omega := ( \omega_n )_{n \geq 1}$, where each $\omega_n$ denotes  a symplectic structure on $X_n$ such that 
$\omega_{m}|_{X_n} = \omega_n$ (if $n < m$).
Since the natural morphism $\Gamma (X, \bigwedge^2 \Omega_{X/\mbZ_p}) \migi \varprojlim_n \Gamma (X_n, \bigwedge^2 \Omega_{X_n/(\mbZ/p^n \mbZ)})$ is an isomorphism (cf. ~\cite{FGA}, Chap.\,8, \S\,8.2, Corollary 8.2.4),
each  symplectic structure $\omega$  on $X$ may be considered as an element of $\Gamma (X, \bigwedge^2 \Omega_{X/\mbZ_p})$.

Denote by
\begin{align}
T^\vee_{\mbZ_p} X
\end{align}
the (smooth) $p$-adic formal stack defined as $T^\vee_{\mbZ_p} X := \varinjlim_n T^\vee_{\mbZ/p^n \mbZ} X_n$.
Then, the collection
\begin{align} \label{G800}
\omega^\mr{Liou}_X := (\omega^\mr{Liou}_{X_n})_{n \geq 1} \in  \Gamma (T^\vee_{\mbZ_p} X, \bigwedge^2 \Omega_{T^\vee_{\mbZ_p} X/\mbZ_p})
\end{align}
forms a symplectic structure on $T^\vee_{\mbZ_p} X$.
The  fiber  of the projection $T^\vee_{\mbZ_p}X \migi X$ over each point in $X (\mbZ_p)$ is Lagrangian.

\vspace{5mm}
\subsection{Moduli of algebraic curves} \label{SS01W}
\leavevmode\\ \vspace{-4mm}

We shall introduce some notation concerning algebraic curves and their moduli.
By a {\bf curve (of genus $g$)} over a fixed scheme $S$, 
we mean
a geometrically connected, proper, and smooth scheme $f : X \migi S$ over $S$ of relative dimension $1$ such that $f_*(\Omega_{X/S})$ is locally free of constant rank $g$.
We shall
denote by
\begin{align}
\mcM_{g, R}
\end{align}
the moduli stack classifying  curves of genus $g$ over $R$, 
 which is a geometrically connected smooth Deligne-Mumford stack over $R$ of relative dimension $3g-3$.
 Also,  denote by 
\begin{align} 
f_{g, R} : C_{g,R} \migi\mcM_{g,R} 
\end{align}
  the universal family of curves over $\mcM_{g,R}$.
In what follows,  we fix a specific choice of an $\mcO_{\mcM_{g,\mbZ}}$-linear  isomorphism
\begin{align} 
    \int_{C_{g,\mbZ}} : \mbR^1f_{g,\mbZ*}(\Omega_{C_{g,\mbZ}/\mcM_{g,\mbZ}})\isom \mcO_{\mcM_{g,\mbZ}} 
   \end{align}
  (i.e., the {\it trace map}) obtained by Serre duality.
For any family of curves $f : X \migi S$ of genus $g$, 
  we shall write 
 \begin{align} \label{G709}
  \int_{X} : \mbR^1f_{*}(\Omega_{X/S})\isom \mcO_S  
  \end{align}
 for
 the pull-back of
   the isomorphism $ \int_{C_{g,\mbZ}}$  via
   the classifying morphism $S \migi \mcM_{g, \mbZ}$ of this curve.
Here, recall (cf. ~\cite{ILL}, Corollary 5.6) that if $d$ denotes the universal derivation  $\mcO_X \migi \Omega_{X/S}$ (namely, the trivial connection on $\mcO_X$ over $S$), then
   $\mbR^1 f_*(\Omega_{X/S})$ is canonically  isomorphic  to $\mbR^2 f_* (\mcK^\bullet [d])$  (cf. \S\,\ref{SS983} for the definition of $\mcK^\bullet [-]$) via the Hodge to de Rham spectral sequence of $\mcK^\bullet [d]$.
  Accordingly, we sometimes consider  $\int_{X}$ as an $\mcO_S$-linear isomorphism $\mbR^2 f_*(\mcK^\bullet [d]) \isom \mcO_S$.
Also, write
\begin{align} \label{G03}
\oint^\natural_X :  f_*(\Omega_{X/S}^{\otimes 2}) \isom  \mbR^1 f_*(\mcT_{X/S})^\vee 
\end{align}
for the isomorphism 
 arising  from the bilinear map 
\begin{align} \label{G04}
\oint_X : f_* (\Omega_{X/S}^{\otimes 2})  \otimes \mbR^1f_* (\mcT_{X/S})
\xrightarrow{\cup} \mbR^1 f_*(\Omega_{X/S}) \xrightarrow{\int_X} \mcO_S.
\end{align}
By well-known generalities on the deformation theory of curves,
there exists  a canonical isomorphism of $\mcO_S$-modules
\begin{align} \label{G808}
\mcT_{\mcM_{g, R}/R} |_S
  \isom  \mbR^1 f_*(\mcT_{X/S})
\end{align}
(i.e., the {\it Kodaira-Spencer map}), where we use the notation ``$|_S$" to denote pull-back by the classifying morphism  $S \migi \mcM_{g, R}$.
This isomorphism  gives the following composite isomorphism: 
\begin{align} \label{G809}
f_*(\Omega_{X/S}^{\otimes 2}) \xrightarrow{\oint^\natural_X} \mbR^1 f_*(\mcT_{X/S})^\vee \xrightarrow{(\ref{G808})^\vee} \left( (\mcT_{\mcM_{g, R}/R} |_{S})^\vee \cong \right) \ \Omega_{\mcM_{g, R}/R} |_S.
\end{align}

Next, by a {\bf $p$-adic formal curve (of genus $g$)} over a $p$-adic formal scheme $S$, we mean 
a flat $p$-adic formal scheme $X$ over $S$ whose reduction modulo $p^n$ (for each $n \geq 1$) is  
a curve (of genus $g$) over $S \otimes (\mbZ/p^n \mbZ)$.
Denote by 
\begin{align} \label{W9}
\widehat{\mcM}_{g,\mbZ_p}
\end{align}
 the smooth $p$-adic formal stack defined as $\widehat{\mcM}_{g, \mbZ_p} := \varinjlim_n \mcM_{g, \mbZ/p^n \mbZ}$.
Then, $\widehat{\mcM}_{g,\mbZ_p}$ may be identified with the moduli classifying $p$-adic formal curves of genus $g$.
By the discussion in the previous subsection, we obtain a symplectic structure
\begin{align} \label{W8}
\widehat{\omega}^\mr{Liou}_{g, \mbZ_p} \ \left(:= \omega^\mr{Liou}_{\widehat{\mcM}_{g, \mbZ_p}} \right)
\end{align}
on $T^\vee_{\mbZ_p} \widehat{\mcM}_{g, \mbZ_p}$.

\vspace{5mm}
\subsection{Indigenous bundles} \label{SS0199}
\leavevmode\\ \vspace{-4mm}

We shall recall the notion of an indigenous bundle.
Some definitions and notation concerning connections on torsors are introduced in the Appendix of the present paper.
Suppose that $2$ is invertible in $R$.
Denote by $B$ the Borel subgroup of $\mr{PGL}_2$ (:= the projective linear group of rank $2$ over $R$) defined to be the image (via the quotient $\mr{GL}_2 \migisurj \mr{PGL}_2$) of upper triangular matrices.
Let 
 $S$  be  a scheme over $R$  and  $f :X \migi S$
a curve 
  of genus $g$ over $S$.
Recall from, e.g., ~\cite{Mzk1}, Chap.\,I, \S\,2, Definition 2.2, or ~\cite{Wak3},  Definition 2.1.1, that an {\bf indigenous bundle} on $X/S$ is 
 a flat $\mr{PGL}_2$-torsor $(\mcE, \nabla_\mcE)$ over $X/S$ (cf. \S\,\ref{SS6791}), i.e., 
 a pair  
 consisting of 
a $\mr{PGL}_2$-torsor $\pi : \mcE \migi X$ over $X$
and an $S$-connection $\mcT_{X/S} \migi \widetilde{\mcT}_{\mcE/S}$
satisfying the following condition:
there exists a $B$-reduction $\mcE_B$ of $\mcE$ (i.e., a $B$-torsor $\pi_B : \mcE_B\migi X$ over $X$ together with an isomorphism $\mcE_B \times^B \mr{PGL}_2 \isom \mcE$ of $\mr{PGL}_2$-torsors), which induces  an inclusion $\widetilde{\mcT}_{\mcE_B/S} \migiincl \widetilde{\mcT}_{\mcE/S}$,  
such that the composite
\begin{align} \label{G803}
\mr{KS}_{(\mcE, \nabla_\mcE)} : \mcT_{X/S} \xrightarrow{\nabla_\mcE} \widetilde{\mcT}_{\mcE/S} \migisurj
\widetilde{\mcT}_{\mcE/S} /\widetilde{\mcT}_{\mcE_B/S}  
\end{align}
is an isomorphism.
If $(\mcE, \nabla_\mcE)$ is an indigenous bundle, then a $B$-reduction $\mcE_B$ of $\mcE$
satisfying the above condition
 is uniquely determined (up to isomorphism);
  we shall refer to it as the {\bf Hodge reduction} of 
   $(\mcE, \nabla_\mcE)$.
An {\bf isomorphism} $(\mcE, \nabla_{\mcE})\isom (\mcE', \nabla_{\mcE'})$ between indigenous bundles on $X/S$ is defined as an isomorphism $\mcE \isom \mcE'$ of $\mr{PGL}_2$-torsors compatible with the respective connections $\nabla_\mcE$ and $\nabla_{\mcE'}$.


 Let $(\mcE, \nabla_\mcE)$ be an indigenous bundle on $X/S$.
$\nabla_\mcE$ induces an $S$-connection
\begin{align}
\nabla^\mr{ad}_\mcE : \mr{Ad}(\mcE) \migi \Omega_{X/S}\otimes \mr{Ad}(\mcE)
\end{align} 
 (cf. (\ref{W30})) on the adjoint bundle  $\mr{Ad}(\mcE) \ \left(:= \mcE \times^{\mr{PGL}_2}  \mfs \mfl_2\right)$.
According to ~\cite{Wak3}, the discussion in  \S\,2.2 (or  ~\cite{Mzk1}, Chap.\,I, \S\,1, the discussion following Definition 1.8), 
there exist canonical injection and surjection
\begin{align}  \label{FF1}
\zeta^\sharp : \Omega_{X/S} \migiincl \mr{Ad}(\mcE), \hspace{10mm} \zeta^\flat : \mr{Ad}(\mcE) \migisurj \mcT_{X/S}
\end{align}
(i.e., $\overline{\nabla}^\sharp$ and $\overline{\nabla}^\flat$ defined in ~\cite{Wak3}, \S\,2.2)
satisfying the equalities
 $\mr{Im}(\zeta^\sharp) = \mcE_B \times^B \mfn$ and $\mr{Ker} (\zeta^\flat) = \mcE_B \times^B \mfb$,  where $\mfb, \mfn \ \left(\subseteq \mfs \mfl_2 \right)$ denote the Lie algebras of $B$, $[B, B]$ ($\subseteq \mr{PGL}_2$) respectively. 
In particular, if we set
\begin{align}
\mr{Ad}(\mcE)^0 := \mr{Ad}(\mcE), \hspace{5mm}
\mr{Ad}(\mcE)^1 := \mr{Ker}(\zeta^\flat), \hspace{5mm}
\mr{Ad}(\mcE)^2 := \mr{Im}(\zeta^\sharp), \hspace{5mm}
\mr{Ad}(\mcE)^3 :=0,
\end{align}
then $\{ \mr{Ad}(\mcE)^j \}_{j=0}^3$ forms a $3$-step decreasing filtration on $\mr{Ad}(\mcE)$ by subbundles
whose subquotients are line bundles with $\mr{Ad}(\mcE)^j /\mr{Ad}(\mcE)^{j+1} \cong \Omega_{X/S}^{\otimes (j-1)}$ ($j=0, 1,2$).

\vspace{5mm}
\subsection{Moduli of indigenous bundles} \label{SS901}
\leavevmode\\ \vspace{-4mm}

Denote by
\begin{align} \label{W4}
\mcS_{g, R}
\end{align}
the moduli stack classifying collections of data  $(X, \mcE, \nabla_\mcE)$ consisting  of a curve $X$  of genus $g$ over $R$ and an indigenous bundle $(\mcE, \nabla_\mcE)$ on it.
By forgetting the data of an indigenous bundle, we obtain a projection $\mcS_{g, R} \migi \mcM_{g, R}$.
According to  ~\cite{Mzk1}, Chap.\,I, \S\,2, Corollary 2.9  (or   ~\cite{Wak3}, Proposition 2.7),
$\mcS_{g,R}$ admits canonically  the  structure of  a relative affine space over $\mcM_{g,R}$  modeled on $f_{g,R*}(\Omega^{\otimes 2}_{C_{g,R}/\mcM_{g,R}})$ (i.e., modeled on $T^\vee_R \mcM_{g, R}$ under  isomorphism (\ref{G809})) that is compatible with base-change over $R$.  
In particular, $\mcS_{g,R}$ is a geometrically connected  smooth Deligne-Mumford stack over $R$ of relative dimension $6g-6$.



Next, we shall write 
\begin{align} \label{W5}
\widehat{\mcS}_{g, \mbZ_p}
\end{align}
 for the $p$-adic formal stack over $\mbZ_p$ defined as
$\widehat{\mcS}_{g, \mbZ_p}: = \varinjlim_n \mcS_{g, \mbZ/p^n \mbZ}$.
By an {\bf indigenous bundle} on a $p$-adic formal curve $X := \varinjlim_n X_n$ (where $X_n := X \otimes (\mbZ/p^n \mbZ)$), 
 we shall mean 
 a collection $((\mcE_n, \nabla_{\mcE, n}))_{n \geq 1}$, where each $(\mcE_n, \nabla_{\mcE, n})$ denotes an indigenous bundle on $X_n$ such that $(\mcE_m, \nabla_{\mcE, m})|_{X_n} \cong (\mcE_n, \nabla_{\mcE, n})$ (if $n <m$).
Then, $\widehat{\mcS}_{g, \mbZ_p}$ may be identified with the moduli classifying  $p$-adic formal curves over $\mbZ_p$ together with an indigenous bundle on it.
Moreover,
the affine space structures on $\mcS_{g, \mbZ/p^n \mbZ}$'s carry the structure of a relative affine space 
 over $\widehat{\mcM}_{g, \mbZ_p}$ modeled on $T^\vee_{\mbZ_p}\widehat{\mcM}_{g, \mbZ_p}$.

Let $f: X \migi S$ be as in \S\,\ref{SS0199} and $(\mcE, \nabla_\mcE)$ an indigenous bundle on $X/S$.
The collection $(X, \mcE, \nabla_\mcE)$ determines its classifying morphism $S \migi \mcS_{g, R}$.
Let us consider the complex  of sheaves $\mcK^\bullet [\nabla^\mr{ad}_\mcE]$ on $X$.
It follows from ~\cite{Mzk1}, Chap.\,I, \S\,2, Theorem 2.8 (or ~\cite{Wak3}, \S\,2.2), 
that 
\begin{align}
\mbR^0f_* (\mcK^\bullet [\nabla^\mr{ad}_\mcE]) = \mbR^2 f_*(\mcK^\bullet [\nabla^\mr{ad}_\mcE]) =0
\end{align}
 and the sequence
\begin{align} \label{G01}
0 \longmigi f_* (\Omega_{X/S}^{\otimes 2}) \stackrel{\xi^\sharp}{\longmigi} \mbR^1 f_* (\mcK^\bullet [\nabla^\mr{ad}_\mcE]) \stackrel{\xi^\flat}{\longmigi} \mbR^1 f_* (\mcT_{X/S}) \longmigi 0
\end{align}
is exact, where 
$\xi^\sharp$  and $\xi^\flat$ denote the morphisms arising from 
$\mr{id}_{\Omega_{X/S}}\otimes \zeta^\sharp : \Omega^{\otimes 2}_{X/S} \migi \Omega_{X/S}\otimes \mr{Ad}(\mcE)$ and $\zeta^\flat : \mr{Ad}(\mcE) \migisurj \mcT_{X/S}$ respectively.
In particular,  $\mbR^1 f_*(\mcK^\bullet [\nabla^\mr{ad}_\mcE])$ is  a vector bundle on $S$ of rank $6g-6$.
Moreover, according to ~\cite{Wak3}, Proposition 2.8.1, 
there exists a canonical isomorphism
\begin{align} \label{EE001}
\mcT_{\mcS_{g, R}/R} |_S \isom \mbR^1 f_{*} (\mcK^\bullet [\nabla_\mcE^\mr{ad}])
\end{align}
of $\mcO_S$-modules fitting into the following isomorphism of short exact sequences:
\begin{align} \label{E0090}
\begin{CD}
0 @>>> \mcT_{\mcS_{g, R}/\mcM_{g, R}}|_S @>>>  \mcT_{\mcS_{g, R}/R}|_S  @>>> \mcT_{\mcM_{g,R}/R} |_S @>>> 0
\\
@. @V \wr VV @V \wr V (\ref{EE001}) V @V \wr V (\ref{G808})V @.
\\
0 @>>> f_* (\Omega_{X/S}^{\otimes 2}) @> \xi^\sharp >> \mbR^1 f_* (\mcK^\bullet [\nabla^\mr{ad}_\mcE]) @> \xi^\flat >> \mbR^1 f_* (\mcT_{X/S}) @>>> 0,
\end{CD}
\end{align}
where the left-hand vertical arrow  arises from the  affine structure on $\mcS_{g,R}$ mentioned above
and the upper horizontal sequence is  obtained by differentiating the projection $\mcS_{g, R} \migi \mcM_{g, R}$.


\vspace{5mm}
\subsection{Symplectic structure on the moduli of indigenous bundles} \label{SS921}
\leavevmode\\ \vspace{-4mm}

Next, we shall construct a canonical symplectic structure on $\mcS_{g, R}$.
Let $f : X \migi S$, $(\mcE, \nabla_\mcE)$ be as above.
Recall that the Killing form  on the Lie algebra $\mfs \mfl_2$ (defined over $R$)  is a nondegenerate symmetric bilinear map  $\kappa : \mfs \mfl_2 \times \mfs \mfl_2 \migi R$ defined by $\kappa (a, b) = \frac{1}{4} \cdot \mr{tr}(\mr{ad}(a)\cdot \mr{ad}(b))$ ($= \mr{tr}(ab)$) for any $a$, $b \in \mfs \mfl_2$.
By 
 the  change of structure group via $\kappa$,
the $\mr{PGL}_2$-torsor $\mcE$ induces  a symmetric  $\mcO_X$-bilinear  morphism
\begin{align} \label{e52}
 \kappa_{(\mcE, \nabla_\mcE)} : \mr{Ad}(\mcE)\otimes \mr{Ad}(\mcE) \migi \mcO_{X},
 \end{align}
which is nondegenerate, i.e., 
 the associated morphism 
 \begin{align} \label{W87}
\kappa_{(\mcE, \nabla_\mcE)}^\natural : \mr{Ad}(\mcE) \migi \mr{Ad}(\mcE)^\vee
\end{align}
 is an isomorphism.
Let us write 
  $ \nabla^{\mr{ad} \otimes 2}_\mcE$ for the connection on the tensor product $\mr{Ad}(\mcE)\otimes  \mr{Ad}(\mcE)$ induced  naturally by $\nabla^\mr{ad}_\mcE$.
The  morphism $\kappa_{(\mcE, \nabla_\mcE)}$ is compatible with the respective connections  $\nabla^{\mr{ad}\otimes 2}_\mcE$ and $d$.
By composing   $\kappa_{(\mcE, \nabla_\mcE)}$ and the cup product in the de Rham cohomology, we obtain  a skew-symmetric $\mcO_S$-bilinear morphism  on $\mbR^1f_{*}(\mcK^\bullet [{\nabla_\mcE^\mr{ad}}])$:
\begin{align} \label{e54}
\oint_{X, (\mcE, \nabla_\mcE)} : \mbR^1f_{*}(\mcK^\bullet [{\nabla_\mcE^\mr{ad}}]) \otimes  
\mbR^1f_{*}(\mcK^\bullet [{\nabla_\mcE^\mr{ad}}])
& \xrightarrow{ \ \ \, \, \cup \ \ \, \,} \ \mbR^2f_{*} (\mcK^\bullet [{\nabla_\mcE^{\mr{ad} \otimes 2}}]) \\
 & \hspace{-0.8mm}\xrightarrow{\kappa_{(\mcE, \nabla_\mcE)}}  \ \mbR^2f_{*}(\mcK^\bullet [d]) \notag  \\ 
 & \hspace{-0.7mm} \xrightarrow{ \ \ \, \int_{X} \ \ \,}  \ \mcO_S. \notag \end{align}
Denote by
\begin{align}
\oint^\natural_{X, (\mcE, \nabla_\mcE)} : \mbR^1 f_*(\mcK^\bullet [\nabla_\mcE^\mr{ad}]) \migi \mbR^1 f_*(\mcK^\bullet [\nabla_\mcE^\mr{ad}])^\vee
\end{align}
the morphism induced by $\oint_{X, (\mcE, \nabla_\mcE)}$, i.e.,
the morphism determined by the condition that
$\oint_{X, (\mcE, \nabla_\mcE)} (a\otimes b) = \left( \oint^\natural_{X, (\mcE, \nabla_\mcE)} a\right)(b)$
for any local sections $a, b \in \mbR^1 f_*(\mcK^\bullet [\nabla^\mr{ad}_\mcE])$.

\vspace{3mm}
\bpr \label{W32}\leavevmode\\
 \ \ \ 
  The morphism $\oint^\natural_{X, (\mcE, \nabla_\mcE)}$ fits into the following morphism of short exact sequences:
 \begin{align} \label{G851}
 \begin{CD}
0  @>>> f_*(\Omega_{X/S}^{\otimes 2})  @> \xi^\sharp >>  \mbR^1 f_*(\mcK^\bullet [\nabla_\mcE^\mr{ad}]) @> \xi^\flat >> \mbR^1 f_* (\mcT_{X/S}) @>>> 0
 \\
 @. @V \wr V \oint^\natural_X V @V \wr V \oint^\natural_{X, (\mcE, \nabla_\mcE)}V @V \wr V \oint^{\natural \vee}_X V @.
 \\
 0 @>>> \mbR^1 f_* (\mcT_{X/S})^\vee  @> (\xi^\flat)^\vee >>  \mbR^1 f_*(\mcK^\bullet [\nabla_\mcE^\mr{ad}])^\vee @> (\xi^\sharp)^\vee >> f_*(\Omega_{X/S}^{\otimes 2})^\vee @>>> 0.
 \end{CD}
 \end{align} 
In particular,  
$\mr{Im} (\xi^\sharp) \ \left(\subseteq \mbR^1 f_*(\mcK^\bullet [\nabla^\mr{ad}_\mcE])\right)$
 is isotropic with respect to $\oint_{X, (\mcE, \nabla_\mcE)}$.
 \epr
\begin{proof}
Note  that  $\mbR^1f_{*}(\mcK^\bullet [\nabla_\mcE^\mr{ad}])$
  may be, locally on $S$, 
  described 
   as the total cohomology of the \v{C}ech double complex $\check{C}^\bullet (\mfU, \mcK^\bullet [\nabla^\mr{ad}_\mcE])$ (for an affine  open covering $\mfU := \{ U_\alpha \}_\alpha$ of $X$) associated to $\mcK^\bullet [\nabla^\mr{ad}_\mcE]$.
Since 
 $\kappa_{(\mcE, \nabla_\mcE)} |_{\mr{Ad}(\mcE)^2 \otimes \mr{Ad}(\mcE)^1}=0$
(because of the  definition of  $\{ \mr{Ad}(\mcE)^j \}_{j=0}^3$),
 this explicit description of $\mbR^1f_{*}(\mcK^\bullet [\nabla_\mcE^\mr{ad}])$
 implies that
  $\mr{Im} (\xi^\sharp) \ \left(= \mr{Ker} (\xi^\flat) \right)$ is isotropic with respect to $\oint_{X, (\mcE, \nabla_\mcE)}$.
In particular, 
we obtain 
two morphisms
\begin{align} \label{G900}
\mr{Im}(\xi^\sharp) \migi (\mbR^1 f_*(\mcK^\bullet [\nabla_\mcE^\mr{ad}])/\mr{Im}(\xi^\sharp))^\vee, \hspace{5mm} \mbR^1 f_*(\mcK^\bullet [\nabla_\mcE^\mr{ad}])/\mr{Im}(\xi^\sharp) \migi \mr{Im}(\xi^\sharp)^\vee
\end{align}
arising as a restriction and a  quotient of $\oint^\natural_{X, (\mcE, \nabla_\mcE)}$ respectively.
Moreover, 
the following diagram is verified to be  commutative:
\begin{align} \label{G850}
\xymatrix{ 
\mr{Ad}(\mcE)^2 \times (\mr{Ad}(\mcE)/\mr{Ad}(\mcE)^1) 
\ar[rr]_{\sim}^{(\overline{\zeta}^\sharp)^{-1} \times \overline{\zeta}^\flat}  
 \ar[rd]_{\overline{\kappa}_{(\mcE, \nabla_\mcE)}} && \Omega_{X/S} \times \mcT_{X/S}  \ar[ld]^{\langle -, - \rangle}\\
& \mcO_{X}, &
}
\end{align}
 where 
 \begin{itemize}
 \item[$\bullet$]
 $\langle -, - \rangle$ denotes the natural paring $\Omega_{X/S}\times \mcT_{X/S} \migi \mcO_{X}$;
 \item[$\bullet$]
 $\overline{\zeta}^\sharp : \Omega_{X/S} \isom \mr{Ad}(\mcE)^2$ and $\overline{\zeta}^\flat : \mr{Ad}(\mcE)/\mr{Ad}(\mcE)^1 \isom \mcT_{X/S}$ denote the isomorphisms induced naturally by $\zeta^\sharp$ and $\zeta^\flat$ respectively;
 \item[$\bullet$]
 $\overline{\kappa}_{(\mcE, \nabla_\mcE)}$ denotes the morphism $\mr{Ad}(\mcE)^2 \times (\mr{Ad}(\mcE)/\mr{Ad}(\mcE)^1)  \migi \mcO_X$ induced by  $\kappa_{(\mcE, \nabla_\mcE)}$. 
 \end{itemize}
 This implies that, 
 under the identifications
 \begin{align}
 f_*(\Omega_{X/S}^{\otimes 2}) \isom \mr{Im}(\xi^\sharp), 
 \hspace{5mm}
 \mbR^1 f_*(\mcK^\bullet [\nabla_\mcE^\mr{ad}])/\mr{Im}(\xi^\sharp) \isom \mbR^1 f_*(\mcT_{X/S})
 \end{align}
 determined  by $\xi^\sharp$, $\xi^\flat$ respectively,
 the isomorphisms in  (\ref{G900})  coincide with 
  $\oint^\natural_X$
    and $\oint^{\natural \vee}_X$ respectively.
This completes the proof of the assertion.
\end{proof}
\vspace{3mm}

By considering the $\mcO_S$-bilinear maps 
$\oint_{X, (\mcE, \nabla_\mcE)}$
(together with the  isomorphism (\ref{EE001}))
for various  schemes $S$ over $\mcS_{g, R}$,
we obtain
  a $2$-form 
\begin{align} \label{e5f5}
 \omega^{\mr{PGL}}_{{g,R}}   \in \Gamma (\mcS_{g, R}, \bigwedge^2 \Omega_{\mcS_{g, R}/R}).
\end{align}
This $2$-form is, by construction,  compatible with base-change over $R$.
Moreover, it follows from ~\cite{Wak3}, Proposition 4.2.2,  that $\omega^{\mr{PGL}}_{{g,R}}$ specifies a symplectic structure on 
$\mcS_{g,R}$.
Thus, the collection
\begin{align} \label{W7}
\widehat{\omega}_{g, \mbZ_p}^{\mr{PGL}} := (\omega^\mr{PGL}_{g, \mbZ/p^n \mbZ})_{n\geq 1}
\in \Gamma (\widehat{\mcS}_{g, \mbZ_p}, \bigwedge^2 \Omega_{\widehat{\mcS}_{g, \mbZ_p}/\mbZ_p})
\end{align}
specifies a symplectic structure on $\widehat{\mcS}_{g, \mbZ_p}$.
One verifies immediately that the  fiber  of the projection $\widehat{\mcS}_{g, \mbZ_p}  \migi \widehat{\mcM}_{g, \mbZ_p}$ over each point in $\widehat{\mcM}_{g, \mbZ_p} (\mbZ_p)$ is Lagrangian.

\vspace{5mm}
\subsection{Ordinary nilpotent indigenous bundles} \label{SS02}
\leavevmode\\ \vspace{-4mm}

Now, we shall  consider certain indigenous bundles in characteristic $p$  playing central roles in  the (classical) $p$-adic Teichm\"{u}ller theory.
Let $S$ be an $\mbF_p$-scheme and $f : X \migi S$ a  curve of genus $g$ over $S$.
Write
$\Phi_S : S \migi S$ for the absolute Frobenius morphism of $S$, $f^{(1)} : X^{(1)} \migi S$ for the Frobenius twist of $X$ relative to $S$, and $\Phi_{X/S} : X \migi X^{(1)}$ for the relative Frobenius morphism (cf. \S\,\ref{SS7791}).
Also, let us fix
 an indigenous bundle $(\mcE, \nabla_\mcE)$
 on $X/S$.
The connection $\nabla_\mcE$ determines its
 $p$-curvature 
$\psi_{(\mcE, \nabla_\mcE)} : \Phi^*_{X/S}(\mcT_{X^{(1)}/S}) \migi \mr{Ad}(\mcE)$
(cf. \S\,\ref{SS7791}).

Recall that $(\mcE, \nabla_\mcE)$ is called {\bf nilpotent} (cf. ~\cite{Mzk1}, Chap.\,II, \S\,2, Definition 2.4) if 
the composite 
\begin{align}
\Phi^*_{X/S}(\mcT_{X^{(1)}/S}) \xrightarrow{\psi_{(\mcE, \nabla_\mcE)}} \mr{Ad} (\mcE) 
\xrightarrow{\kappa_{(\mcE, \nabla_\mcE)}^\natural}
 \mr{Ad}(\mcE)^\vee \xrightarrow{\psi_{(\mcE, \nabla_\mcE)}^\vee} \Phi^*_{X/S}(\mcT_{X^{(1)}/S})^\vee
\end{align}
coincides with  the zero map.
In particular, if $(\mcE, \nabla_\mcE)$ is nilpotent, then it is $p$-nilpotent in the sense discussed in \S\,\ref{SS7791}, and hence, corresponds to a crystal of $\mr{PGL}_2$-torsors over  the crystalline site  $\mr{Crys} (X/S)$ (cf. Remark \ref{W40} (ii) and Theorem \ref{T043}).

Next, let us consider the composite
\begin{align} \label{K001}
\Phi^{-1}_{X/S}(\mcT_{X^{(1)}/S})
\xrightarrow{\psi^\nabla_{(\mcE, \nabla_\mcE)}}\mr{Ad}(\mcE)  \stackrel{\zeta^\flat}{\migisurj} \mcT_{X/S}.
\end{align}
(cf. (\ref{W118}) for the definition of $\psi^\nabla_{(\mcE, \nabla_\mcE)}$).
By applying the functor  $\mbR^1 f_*(-)$ to  this composite, we  obtain an $\mcO_S$-linear morphism
\begin{align} \label{G02}
 \Phi_S^*(\mbR^1f_* (\mcT_{X/S}))  \left(\isom \mbR^1 f^{(1)}_* (\mcT_{X^{(1)}/S})\right) \migi 
  \mbR^1 f_*(\mcT_{X/S}).
\end{align}
Then, recall that $(\mcE, \nabla_\mcE)$ is called {\bf ordinary}  (cf. ~\cite{Mzk1}, Chap.\,II,  \S\,3,  Definition 3.1) if the morphism
 (\ref{G02})
  is an isomorphism.

Denote by
\begin{align} \label{W1}
\mcN^\mr{ord}_{g, \mbF_p}
\end{align}
the substack of $\mcS_{g, \mbF_p}$ classifying 
ordinary nilpotent indigenous bundles.
It follows from ~\cite{Mzk1}, Chap.\,II, \S\,3, Corollary 3.8, that
$\mcN^\mr{ord}_{g, \mbF_p}$ is a nonempty smooth Deligne-Mumford stack over $\mbF_p$ and
the projection 
$\mcN^\mr{ord}_{g, \mbF_p} \migi \mcM_{g,\mbF_p}$ is  \'{e}tale and quasi-finite  (and hence,  since $\mcM_{g, \mbF_p}$ is irreducible,  it is  dominant when restricted to  each component of $\mcN^\mr{ord}_{g, \mbF_p}$).
Therefore, there exists a unique (up to isomorphism) smooth $p$-adic formal stack 
\begin{align} \label{W2}
\widehat{\mcN}^\mr{ord}_{g, \mbZ_p}
\end{align}
 whose reduction modulo $p$ is $\mcN^\mr{ord}_{g, \mbF_p}$.
Let us recall  here the following assertion, which is 
 one of the main results in ~\cite{Mzk1}.

\vspace{3mm}
\bt[cf. ~\cite{Mzk1}, Chap.\,III, \S\,2, Theorem 2.8] \label{T0135}  \leavevmode\\
 \ \ \ 
 Denote by  $f_{\mcN} : C_{\mcN} \migi  \widehat{\mcN}^\mr{ord}_{g, \mbZ_p}$  the universal family of  curves over $ \widehat{\mcN}^\mr{ord}_{g, \mbZ_p}$.
 Then, there exists a canonical endomorphism
 \begin{align} \label{W3}
\Phi_{\mcN} :  \widehat{\mcN}^\mr{ord}_{g, \mbZ_p} \migi \widehat{\mcN}^\mr{ord}_{g, \mbZ_p}
\end{align}
of $\widehat{\mcN}^\mr{ord}_{g, \mbZ_p}$ together with a canonical indigenous bundle $(\mcE_{\mcN}, \nabla_{\mcE_{\mcN}})$ on $C_{\mcN}/  \widehat{\mcN}^\mr{ord}_{g, \mbZ_p}$
satisfying the following properties:
\begin{itemize}
\item[$\bullet$]
$\Phi_{\mcN}$ is a Frobenius lifting over $\mbZ_p$ (i.e., the reduction modulo $p$ of $\Phi_{\mcN}$ coincides with  the absolute Frobenius morphism of $ \widehat{\mcN}^\mr{ord}_{g, \mbZ_p}$);
\item[$\bullet$]
The reduction modulo $p$ of $(\mcE_{\mcN}, \nabla_{\mcE_{\mcN}})$  is isomorphic to  the indigenous bundle on $(C_{\mcN}\otimes \mbF_p) / \mcN^\mr{ord}_{g, \mbF_p} $  classified by the natural  immersion 
$\mcN^\mr{ord}_{g, \mbF_p}\migiincl \mcS_{g, \mbF_p}$;
\item[$\bullet$]
There exists an isomorphism $\mbF^*\Phi^*_{\mcN}(\mcE_{\mcN}, \nabla_{\mcE_{\mcN}})\isom (\mcE_{\mcN}, \nabla_{\mcE_{\mcN}})$ of indigenous bundles, where $\mbF^* (-)$ denotes   renormalized Frobenius pull-back (cf.   ~\cite{Mzk1}, Chap.\,III, \S\,2, Definition 2.4, or \S\,\ref{SS179} in the present paper) and $\Phi^*_\mcN (-)$ denotes base-change by  $\Phi_\mcN$.
\end{itemize}
Moreover, the collection of data  $(\Phi_{\mcN}, \mcE_{\mcN}, \nabla_{\mcE_{\mcN}})$ is uniquely characterized (up to isomorphism) by the above properties.
   \et


\vspace{5mm}
\subsection{Statement of the main theorem} \label{SS029}
\leavevmode\\ \vspace{-4mm}

In this subsection, we  shall 
describe the main theorem in the present paper.
In what follows,
we shall write $\mcN := \widehat{\mcN}_{g, \mbZ_p}^\mr{ord}$ for simplicity.
The indigenous bundle $(\mcE_{\mcN}, \nabla_{\mcE_{\mcN}})$ obtained in  Theorem \ref{T0135} determines its classifying  morphism
\begin{align} \label{G230}
\sigma : \mcN \migi \widehat{\mcS}_{g, \mbZ_p}
\end{align}
 over $\widehat{\mcM}_{g, \mbZ_p}$, which turns out to be  an immersion;
 it gives, after base-change by $\mcN \migi \widehat{\mcM}_{g, \mbZ_p}$,  a trivialization of 
 the affine space structure on  $\widehat{\mcS}_{g, \mbZ_p}$
 (modeled on $T^\vee_{\mbZ_p} \widehat{\mcM}_{g, \mbZ_p}$).
More precisely, there exists  a unique isomorphism
 \begin{align} \label{e58}
  \theta : 
  T^\vee_{\mbZ_p}  \widehat{\mcM}_{g, \mbZ_p} |_{\mcN} 
  \ \left(:= T^\vee_{\mbZ_p}  \widehat{\mcM}_{g, \mbZ_p} \times_{\widehat{\mcM}_{g, \mbZ_p}} \mcN \right) \isom  \widehat{\mcS}_{g, \mbZ_p} |_{\mcN} \ \left(:= \widehat{\mcS}_{g, \mbZ_p} \times_{\widehat{\mcM}_{g, \mbZ_p}} \mcN \right)
  \end{align}
 which extends  $\sigma$ and  is compatible with the affine space structures  pulled-back from  $T^\vee_{\mbZ_p}  \widehat{\mcM}_{g, \mbZ_p}$ and $\widehat{\mcS}_{g, \mbZ_p}$ respectively.
It induces 
 an isomorphism 
 \begin{align} \label{G240}
 \theta^*(\Omega_{\widehat{\mcS}_{g, \mbZ_p} |_\mcN/\mbZ_p})\isom \Omega_{T^\vee_{\mbZ_p} \widehat{\mcM}_{g, \mbZ_p} |_\mcN/\mbZ_p},
 \end{align} 
 and hence, 
  an isomorphism
 \begin{align} \label{G241}
\Theta :  \theta^*(\bigwedge^2 \Omega_{\widehat{\mcS}_{g, \mbZ_p} |_\mcN/\mbZ_p})
 \ \left(\cong \bigwedge^2 \theta^*( \Omega_{\widehat{\mcS}_{g, \mbZ_p} |_\mcN/\mbZ_p}) \right)
\isom 
\bigwedge^2 \Omega_{T^\vee_{\mbZ_p} \widehat{\mcM}_{g, \mbZ_p} |_\mcN/\mbZ_p}.
 \end{align}
 Since $\mcN$ is \'{e}tale over $\widehat{\mcM}_{g, \mbZ_p}$,
 the projections $T^\vee_{\mbZ_p}\widehat{\mcM}_{g, \mbZ_p} |_\mcN \migi T^\vee_{\mbZ_p}\widehat{\mcM}_{g, \mbZ_p}$
 and 
  $\widehat{\mcS}_{g, \mbZ_p} |_\mcN \migi \widehat{\mcS}_{g, \mbZ_p}$ are \'{e}tale.
 Therefore, the $2$-form 
 \begin{align}
 \widehat{\omega}^{\mr{Liou}}_{g, \mbZ_p} |_{\mcN} \ \left(\text{resp.,} \ \widehat{\omega}^\mr{PGL}_{g, \mbZ_p} |_\mcN \right)
 \end{align}
 on $T^\vee_{\mbZ_p}  \widehat{\mcM}_{g, \mbZ_p} |_{\mcN}$ (resp., $\widehat{\mcS}_{g, \mbZ_p}|_{\mcN}$) defined as the pull-back of $\widehat{\omega}^\mr{Liou}_{g, \mbZ_p}$   (resp., $\widehat{\omega}^\mr{PGL}_{g, \mbZ_p}$)
 specifies  a symplectic structure. 
(Notice that $\widehat{\omega}^\mr{Liou}_{g, \mbZ_p}|_{\mcN} = \omega^\mr{Liou}_{\mcN}$.)
The main result of the present paper is the following  Theorem \ref{T013}, which  describes the relationship between
$\widehat{\omega}^{\mr{Liou}}_{g, \mbZ_p} |_{\mcN}$ and  $\widehat{\omega}^\mr{PGL}_{g, \mbZ_p} |_\mcN$.
(The proof will be given
in \S\,\ref{SS699}.)

\vspace{3mm}
\bt[= Theorem A] \label{T013}  \leavevmode\\
 \ \ \ 
 If $p>3$, then the  morphism $\theta$ preserves the symplectic structures, i.e., the following equality holds:
 \begin{align} \label{W210}
 \Theta (\widehat{\omega}^{\mr{PGL}}_{g, \mbZ_p} |_{\mcN}) = \widehat{\omega}^\mr{Liou}_{g, \mbZ_p} |_\mcN.
 \end{align}
 In particular, the image of $\sigma :  \mcN \migi \widehat{\mcS}_{g, \mbZ_p}$ is Lagrangian  with respect to the symplectic structure $\widehat{\omega}_{g, \mbZ_p}^\mr{PGL}$.
  \et

\vspace{10mm}
\section{$F$-crystals associated to ordinary nilpotent indigenous bundles} \vspace{3mm}

Before proving  Theorem \ref{T013}, we shall study, in this section,  a certain  $F$-crystal structure  (cf. (\ref{G953})) on the cohomology associated to  the adjoint bundle of an ordinary nilpotent indigenous bundle.
One important observation is (cf. Corollary \ref{C0338}) that
the direct sum decomposition  of this cohomology determined by 
(the differential of) $\sigma : \widehat{\mcN}^\mr{ord}_{g, \mbZ_p} \migi \widehat{\mcS}_{g, \mbZ_p}$ 
 coincides with (i.e,  gives the geometric interpretation of) the slope decomposition.
It follows (cf. Corollary \ref{G970} and  (\ref{G401})) that  both  $\widehat{\omega}^\mr{Liou}_{g, \mbZ_p} |_\mcN$ and $\widehat{\omega}^\mr{PGL}_{g, \mbZ_p} |_\mcN$ turn out to specify  eigenvectors of the $F$-crystal structure defined on 
the second exterior power of the dual   $\mbH^1 (\mcK^\bullet [\nabla^\mr{ad}_\mcE])^\vee$.
This fact, being an essential point of our proof of Theorem \ref{T013}, makes it clear how the two symplectic structures 
 are related via reduction modulo $p$.


\vspace{5mm}
\subsection{Renormalized Frobenius pull-back} \label{SS179}
\leavevmode\\ \vspace{-4mm}

First, let us recall (cf. ~\cite{Mzk1}, Chap.\,III, \S\,2, the discussion preceding Definition 2.4) the definition of renormalized Frobenius pull-back, that appeared in the statement of 
Theorem \ref{T0135}.
In what follows, we shall denote, for each positive integer $m$,  the reductions of objects over $\mbZ_p$ to $\mbZ /p^m \mbZ$ by means of a subscripted $m$.
Let $S$ be a $p$-adic formal scheme.
Also, let  $X$ and $Y$ be   curves of genus $g$ over $S$
  such that $Y$ is a $p$-adic  lifting of $X^{(1)}_1 \ \left(:= (X_1)^{(1)}= (X^{(1)})_1\right)$.

We shall fix a flat $\mr{PGL}_2$-torsor $(\mcE, \nabla_\mcE)$ 
 over $Y/S$
   whose reduction modulo $p$ (i.e., $(\mcE_1, \nabla_{\mcE, 1})$) forms a nilpotent  indigenous bundle on $X_1^{(1)} \ (=Y_1)/S_1$.
Let $n$ be an integer  with $n>1$, and 
assume  {\it tentatively} that there exists a  rank $2$  flat vector bundle $(\mcV_n, \nabla_{\mcV, n})$ on $Y_n/S_n$ (i.e., a pair of a rank $2$ vector bundle $\mcV_n$ on $Y_n$ and an $S_n$-connection $\nabla_{\mcV, n}$ on $\mcV_n$)  whose projectivization is isomorphic to $(\mcE_n, \nabla_{\mcE, n})$. 
Denote by $\mcL_1$ the line subbundle of $\mcV_1 \ (:= (\mcV_n)_1)$ corresponding to  the Hodge reduction of the indigenous bundle  $(\mcE_1, \nabla_{\mcE, 1})$.
Since $\nabla_{\mcV, 1}$ has nilpotent $p$-curvature, 
$(\mcV_n, \nabla_{\mcV, n})$ corresponds 
to a crystal  $\mcV^\lozenge_n$ of  vector bundles  on the crystalline site $\mr{Crys}(X^{(1)}_1/S_n)$.
Moreover, it induces a crystal $\Phi_{X_1/S_1}^*(\mcV^\lozenge_n)$  on $\mr{Crys}(X_1/S_n)$ defined as the pull-back of $\mcV^\lozenge_n$  via the relative Frobenius  $\Phi_{X_1/S_1}$.
One may obtain a crystal 
$\widetilde{\mbF}^*(\mcV^\lozenge_n)$
defined as the subsheaf of $\Phi_{X_1/S_1}^*(\mcV^\lozenge_n)$
consisting of sections whose reduction modulo $p$ are contained in the subsheaf $\Phi_{X_1/S_1}^*(\mcL_1) \ \left( \subseteq \Phi_{X_1/S_1}^*(\mcV_1) \right)$. 
If $\nabla'$ denotes the $S_n$-connection on the $\mcO_{X_n}$-module  $\widetilde{\mbF}^*(\mcV^\lozenge_n)_{X_n}$  (i.e., the evaluation of $\widetilde{\mbF}^*(\mcV^\lozenge_n)$ at $X_n$) corresponding to this crystal,
then its  reduction $(\widetilde{\mbF}^*(\mcV^\lozenge_n)_{X_{n-1}}, \nabla'_{n-1})$ modulo $p^{n-1}$ turns out to form a flat {\it vector bundle} on $X_{n-1}/S_{n-1}$.
We shall write 
\begin{align} \label{G801}
\mbF^*(\mcV_n^\lozenge)
\end{align}
 for the crystal of vector bundles on $\mr{Crys}(X_1/S_{n-1})$ corresponding to this flat vector bundle.
The (isomorphism class of the) flat $\mr{PGL}_2$-torsor over $X_{n-1} /S_{n-1}$
 defined as the projectivization of  (the flat vector bundle corresponding to) $\mbF^*(\mcV^\lozenge_n)$
 is independent of
  the choice of $(\mcV_n, \nabla_{\mcV, n})$ (i.e., depends only on $(\mcE_n, \nabla_{\mcE, n})$); thus it makes sense to 
  use the notation
  $\mbF^*(\mcE_n, \nabla_{\mcE, n})$
to denote this flat $\mr{PGL}_2$-torsor.
Moreover,  the  independence of the choice $(\mcV_n, \nabla_{\mcV, n})$
implies that {\it we can construct $\mbF^*(\mcE_n, \nabla_{\mcE, n})$ without
  the existence assumption  of $(\mcV_n, \nabla_{\mcV, n})$ imposed above.} 
By applying the above argument to all $n$, we obtain a flat $\mr{PGL}_2$-torsor
\begin{align} \label{G80300}
\mbF^*(\mcE, \nabla_\mcE)
\end{align}
over $X/S$, which we call the {\bf renormalized Frobenius pull-back} of $(\mcE, \nabla_\mcE)$.
Denote by
\begin{align}
\mbF^*(\mcE, \nabla_\mcE)^\lozenge
\end{align}
the crystal 
  over  $\mr{Crys}(X_1/S)$ (i.e., the compatible system consisting of  crystals over $\mr{Crys}(X_1/S_n)$ for various $n$) corresponding to $\mbF^*(\mcE, \nabla_\mcE)$.
Note that (the isomorphism class of) $\mbF^*(\mcE, \nabla_\mcE)^\lozenge$ does not depend on
 the choice of the $p$-adic lifting $X$ of $X_1$.

\vspace{5mm}
\subsection{Frobenius structure on the cohomology of  the  adjoint bundle} \label{SS379}
\leavevmode\\ \vspace{-4mm}

In this subsection, we shall construct a Frobenius structure on the cohomology associated with the  adjoint bundle of each ordinary nilpotent indigenous bundle.
Let $S$, $\Phi_{S}$, and $X$  be as in the previous subsection, and let us keep some  notational convention 
  as needed (e.g., $(-)_n$ and  $(-)^{(1)}$, etc.).
Assume that
$S$ is endowed  with a morphism $S \migi \mcN$  via which   $\Phi_{S}$ is compatible with $\Phi_{\mcN}$.
Denote by $(\mcE, \nabla_\mcE)$ the ordinary nilpotent indigenous bundle on $X/S$ classified by this morphism.
Since the reduction modulo $p$ of $\nabla_\mcE$ (as well as $\nabla_\mcE^\mr{ad}$) has nilpotent $p$-curvature,
the flat vector bundle $(\mr{Ad}(\mcE_n), \nabla_{\mcE, n}^\mr{ad})$
determines a crystal 
 $\mr{Ad}(\mcE_n)^\lozenge$ on $\mr{Crys}(X_1/S_n)$.
 In particular, we obtain the relative crystalline cohomology sheaf  $\mbR^1 f_{\mr{crys}*} (\mr{Ad}(\mcE_n)^\lozenge)$ on $S_n$ associated to  $\mr{Ad}(\mcE_n)^\lozenge$, and hence, obtain 
 the  $\mcO_S$-module 
 \begin{align}
 \mbR^1 f_{\mr{crys}*} (\mr{Ad}(\mcE)^\lozenge) := \varprojlim_n  \mbR^1 f_{\mr{crys}*} (\mr{Ad}(\mcE_n)^\lozenge).
 \end{align}

Suppose {\it tentatively} that 
there exists  a crystal  $\mcV^\lozenge_n$ of rank $2$ vector bundles on $\mr{Crys}(X_1 /S_n)$ whose projectivization corresponds to $(\mcE_n, \nabla_{\mcE, n})$.
Write  $\Phi_S^*(\mcE_n, \nabla_{\mcE, n})$ (resp.,  $\Phi_{S}^*(\mcV_n^\lozenge)$) for  the base-change of $(\mcE_n, \nabla_{\mcE, n})$ (resp., $\mcV^\lozenge_n$) 
by $\Phi_{S, n} : S_n \migi S_n$, which forms an  indigenous bundle on $X^{(1)}_{n}/S_n$  (resp., a crystal on $\mr{Crys}(X_1^{(1)} /S_n)$). 
Here, notice that  $\mr{Ad}(\mcE_n)^\lozenge$ may be identified with the crystal which assigns, to each  $(U \migiincl T, \delta)$ in $\mr{Crys}(X_1/S_n)$, the sheaf $\mcE nd^0 (\mcV^\lozenge_{n, T})$ of  $\mcO_T$-linear endomorphisms of $\mcV_{n, T}^\lozenge$ with vanishing trace.

Now, let us take a divided power thickening $(U \migiincl T, \delta)$ in $\mr{Crys}(X_1/S_n)$ and 
an $\mcO_T$-linear endomorphism $h$ of $\mcV_{n, T}^\lozenge$ with vanishing trace (i.e., a global section of $\mcE nd^0 (\mcV_{n, T}^\lozenge)$).
After possibly replacing $U$ with its open covering, we suppose
that $T$ admits an endomorphsim $\Phi_T : T \migi T$ compatible with $\Phi_{S, n}$ whose reduction modulo $p$ coincides with
the absolute Frobenius $\Phi_U$.
The endomorphism $p \cdot \Phi_T^*(h)$ of $\Phi^*_T(\mcV^\lozenge_{n, T})$
restricts to an endomorphism  of $\widetilde{\mbF}^*(\Phi_S^*(\mcV^\lozenge_{n}))_T \ \left(\subseteq \Phi^*_T(\mcV^\lozenge_{n, T}) = F_{X_1/S_1}^*(\Phi_S^*(\mcV^\lozenge_n))_T\right)$; we shall denote its reduction modulo $p^{n-1}$ by
$\mbF^*(h)$, which lies in $\mcE nd^0 (\mbF^*(\Phi^*_S(\mcV_n^\lozenge))_{\overline{T}})$ (where $\overline{T} := T_{n-1}$).
Since the assignment $\mbF^*(-)$ is compatible  with base-change over the parameter spaces of   underlying  families of curves,
it follows from
 Theorem \ref{T0135} (and the assumption that $\Phi_S$ is compatible with $\Phi_\mcN$) that 
$\mbF^*(\Phi^*_S(\mcE, \nabla_\mcE))$ is isomorphic to $(\mcE, \nabla_\mcE)$.
Hence, $\mbF^*(\Phi^*_S(\mcV_n^\lozenge))_{\overline{T}}$ is  isomorphic to
$\mcV_{n, \overline{T}}^\lozenge$ 
up to tensoring with  a line bundle, and 
$\mbF^*(h)$ determines a well-defined section of $\mcE nd^0 (\mcV_{n, \overline{T}}^\lozenge) \ \left(= \mr{Ad}(\mcE_{n-1})^\lozenge_{\overline{T}}\right)$. 
Given   a basis $(e_1,  e_2)$ of $\mcV^\lozenge_{n, T}$ such that
 $e_1$ mod $p$ generates $\mcL_1$, we can describe  locally
 $\mbF^*(h)$ by means of this basis.
Indeed, 
 if the matrix representation of $h$ with respect to the basis $(e_1, e_2)$ is of the form
 $\begin{pmatrix} a & b \\ c & -a \end{pmatrix}$,
 then $\mbF^*(h)$ may be expressed as 
  $\begin{pmatrix} p\cdot \Phi^*_T(a) & p^2\cdot \Phi^*_T(b) \\ \Phi^*_T(c) & -p \cdot\Phi^*_T(a) \end{pmatrix}$ (mod $p^{n-1}$) with respect to the basis $(\Phi^*_T(e_1), p\cdot  \Phi^*_T(e_2))$ mod $p^{n-1}$.
 
 Denote by
  \begin{align} \label{W2001}
 (\Phi_{X_1})_{\mr{crys}} : (X_1/S_{n-1})_\mr{crys} \migi (X_1/S_{n})_\mr{crys}
 \end{align}
  the morphism of topoi  induced by the absolute Frobenius endomorphism $\Phi_{X_1}$ of $X_1$ covering  the PD morphism 
  $\Phi_{S, n} |_{S_{ n-1}} : S_{n-1}\migi S_n$.
 Then, the  assignment $h \mapsto \mbF^*(h)$ (for each  $(U \migiincl T, \delta)$ and  $h$ as above) determines 
 a morphism of crystals 
 \begin{align} \label{G740}
 \mr{Ad} (\mbF^*) : \mr{Ad}(\mcE_n)^\lozenge \migi (\Phi_{X_1})_{\mr{crys}*}(\mr{Ad}(\mcE_{n-1})^\lozenge).
 \end{align}
 It induces  
 a $\Phi_{S, n}$-linear morphism  $\mbR^1 f_{n, \mr{crys}*} (\mr{Ad}(\mcE_n)^\lozenge) \migi \mbR^1 f_{n-1, \mr{crys}*} (\mr{Ad}(\mcE_{n-1})^\lozenge)$, or equivalently, an
 $\mcO_{S_n}$-linear morphism
 \begin{align} \label{G201}
 F^\lozenge_{(\mcE_n, \nabla_{\mcE, n})} : \Phi_{S, n}^*(\mbR^1 f_{n, \mr{crys}*} (\mr{Ad}(\mcE_n)^\lozenge)) \migi  \mbR^1 f_{n-1, \mr{crys}*} (\mr{Ad}(\mcE_{n-1})^\lozenge).
 \end{align}
One verifies immediately that this morphism is independent of the choice of $\mcV_n^\lozenge$ (i.e., depends only on $(\mcE_n, \nabla_{\mcE, n})$), which implies that {\it we can remove the existence assumption of $\mcV^\lozenge_n$ imposed above.}
By applying the above argument to all  $n$, we obtain from
$F^\lozenge_{(\mcE_n, \nabla_{\mcE, n})}$'s
  an $\mcO_S$-linear morphism
\begin{align} \label{G953}
 F^\lozenge_{(\mcE, \nabla_{\mcE})} : \Phi_{S}^*(\mbR^1 f_{\mr{crys}*} (\mr{Ad}(\mcE)^\lozenge)) \migi  \mbR^1 f_{\mr{crys}*} (\mr{Ad}(\mcE)^\lozenge).
\end{align}

Here, recall from ~\cite{BO},  Theorem 7.1,  that there exists a canonical isomorphism 
\begin{align}  \label{G952}
\mbR^1 f_{n, \mr{crys}*} (\mr{Ad}(\mcE_n)^\lozenge) \isom \mbR^1 f_{n*}(\mcK^\bullet [\nabla^\mr{ad}_{\mcE_n}])
\end{align}
of $\mcO_{S_n}$-modules.
By taking the inverse  limit $\varprojlim_n (-)$,
we obtain   an isomorphism
\begin{align} \label{G980}
\mbR^1 f_{\mr{crys}*} (\mr{Ad}(\mcE)^\lozenge) \isom  \mbR^1 f_{*}(\mcK^\bullet [\nabla^\mr{ad}_{\mcE}]) \ \left(:= \varprojlim_n  \mbR^1 f_{n*}(\mcK^\bullet [\nabla^\mr{ad}_{\mcE_n}]) \right).
\end{align}
The isomorphism $F^\lozenge_{(\mcE, \nabla_\mcE)}$ becomes,
via (\ref{G980}), an $\mcO_S$-linear morphism
\begin{align} \label{G951}
F_{(\mcE, \nabla_\mcE)} : \Phi^*_S (\mbR^1 f_*(\mcK^\bullet [\nabla^\mr{ad}_\mcE])) \migi  \mbR^1 f_*(\mcK^\bullet [\nabla^\mr{ad}_\mcE]).
\end{align}
In particular,  we apply this construction of $F_{(\mcE, \nabla_\mcE)}$ to the universal case (i.e., the case where the collection $(S, \Phi_S, \mcE, \nabla_\mcE)$ is taken as $(\mcN, \Phi_\mcN, \mcE_\mcN, \nabla_{\mcE_\mcN})$),
  and obtain an $\mcO_{\mcN}$-linear morphism
\begin{align} \label{G807}
F
 : \Phi^*_\mcN (\mbR^1 f_{\mcN*}(\mcK^\bullet [\nabla^\mr{ad}_{\mcE_{\mcN}}])) \migi \mbR^1 f_{\mcN*}(\mcK^\bullet [\nabla^\mr{ad}_{\mcE_{\mcN}}]).
\end{align}

\vspace{5mm}
\subsection{Relationship between $\oint_{X, (\mcE_\mcN, \nabla_\mcN)}$ and $F$} \label{SS0719}
\leavevmode\\ \vspace{-4mm}

By the following assertion, we can see that
the morphism $F$ defined above is compatible, up to multiplication by ``$p^3$",  with the bilinear morphism 
$\oint_{X, (\mcE_\mcN, \nabla_\mcN)}$.


\vspace{3mm}
\bpr \label{P990} \leavevmode\\
 \ \ \ 
 The following square diagram is commutative:
 \begin{align} \label{G731}
\begin{CD}
\Phi^*_{\mcN}(\mbR^1f_{\mcN*} (\mcK^\bullet [\nabla^{\mr{ad}}_{\mcE_{\mcN}}])) \otimes \Phi^*_{\mcN}(\mbR^1 f_{\mcN*} (\mcK^\bullet [\nabla^{\mr{ad}}_{\mcE}]))@> \Phi_\mcN^*(\oint_{X, (\mcE_\mcN, \nabla_\mcN)}) >>  \Phi^*_{\mcN} (\mcO_{\mcN}) \ \left(\cong \mcO_{\mcN} \right)
\\
@V F \otimes F VV @VV [p^3] V
\\
\mbR^1f_{\mcN*} (\mcK^\bullet [\nabla^{\mr{ad}}_{\mcE_{\mcN}}]) \otimes \mbR^1f_{\mcN*} (\mcK^\bullet [\nabla^{\mr{ad}}_{\mcE_{\mcN}}])  @>> \oint_{X, (\mcE_\mcN, \nabla_\mcN)} >
\mcO_{\mcN},
\end{CD}
\end{align}
where the right-hand vertical arrow $[p^3]$ denotes  multiplication by $p^3$.
 \epr
\begin{proof}
Let us keep the notation in \S\,\ref{SS379}.
The Killing form $\kappa$ on $\mfs \mfl_2$ induces, for each $n$, 
a morphism of crystals $\kappa_{(\mcE_n, \nabla_{\mcE, n})}^\lozenge : \mr{Ad}(\mcE_n)^\lozenge \times \mr{Ad}(\mcE_n)^\lozenge \migi \mcO_{X_1/S_n}$.
Moreover,  $\kappa_{(\mcE_n, \nabla_{\mcE, n})}^\lozenge$
 induces
 an $\mcO_{S_n}$-bilinear morphism
\begin{align}
\oint_{X_n, (\mcE_n, \nabla_{\mcE, n})}^\lozenge : \mbR^1 f_{n, \mr{crys}*} (\mr{Ad}(\mcE_n)^\lozenge) \otimes \mbR^1 f_{n, \mr{crys}*} (\mr{Ad}(\mcE_n)^\lozenge) \migi \mcO_{S_n},
\end{align}
which  is  compatible with $\oint_{X_n, (\mcE_n, \nabla_{\mcE, n})}$ via the  isomorphism (\ref{G952}).
Also, by taking account of    the local description of $\mbF^*(-)$  discussed in \S\,\ref{SS379}, 
we see that  the following square   diagram is  commutative:
\begin{align} \label{G704}
\begin{CD}
\mr{Ad}(\mcE_n)^\lozenge \times \mr{Ad}(\mcE_n)^\lozenge @>\kappa_{(\mcE_n, \nabla_{\mcE, n})}^\lozenge >> \mcO_{X_1 /S_n}
\\
@V \mr{Ad} (\mbF^*) \times \mr{Ad}(\mbF^*) VV @VV p^2 \cdot \Phi^*_{X_1/S_1}V
\\
(\Phi_{X_1})_{\mr{crys}*}(\mr{Ad}(\mcE_{n-1})^\lozenge) \times (\Phi_{X_1})_{\mr{crys}*}(\mr{Ad}(\mcE_{n-1})^\lozenge) @>>(\Phi_{X_1})_{\mr{crys}*}(\kappa_{(\mcE_{n-1}, \nabla_{\mcE, n-1})}^\lozenge) > (\Phi_{X_1})_{\mr{crys}*}(\mcO_{X_1 /S_{n-1}})
\end{CD}
\end{align}
(cf. (\ref{W2001}) for the definition of $(\Phi_{X_1})_{\mr{crys}*}$),  where the right-hand vertical arrow $p^2  \cdot \Phi^*_{X_1/S_1}$ denotes 
$p^2$ times
the morphism  $\Phi^*_{X_1/S_1} : \mcO_{X_1 /S_n} \migisurj (\Phi_{X_1})_{\mr{crys}*}(\mcO_{X_1 /S_{n-1}})$   induced by $\Phi_{X_1/S_1}$.
Here, recall from ~\cite{BO}, Theorem 6.12, that we have 
\begin{align} \label{W46}
    \mbR^2 f_{1, \mr{crys}*}(\mcO_{X_1/S_m})\isom
    \left(\mbR^2 f_{m*} (\mcK^\bullet [d]) \isom \right) \mcO_{S_m}
\end{align}
($m=1,2, \cdots$).
 The morphism 
 \begin{align}
 \Phi^*_{S, n}(\mbR^2f_{1, \mr{crys}*}(\mcO_{X_1/S_n})) \migi \mbR^2 f_{1, \mr{crys}*}(\mcO_{X_1/S_{n-1}})
 \end{align}
  induced  by $F^*_{X_1/S_1} : \mcO_{X_1 /S_n} \migisurj (\Phi_{X_1})_{\mr{crys}*}(\mcO_{X_1 /S_{n-1}})$ 
  coincide, via (\ref{W46}),  with the composite of the natural quotient   $\Phi^*_{S, n}(\mcO_{S_n}) \ \left(\cong \mcO_{S_n} \right) \migisurj \mcO_{S_{n-1}}$ and  multiplication by $p$ (cf. ~\cite{B}, Chap.\,VII, \S\,3, Proposition 3.2.4).
Hence, the  diagram (\ref{G704})  gives rise to a commutative diagram of the form
\begin{align} \label{G730}
\begin{CD}
\Phi^*_{S, n}(\mbR^1f_{1, \mr{crys}*} (\mr{Ad}(\mcE_n)^\lozenge) \otimes \Phi^*_{S, n}(\mbR^1 f_{1, \mr{crys}*} (\mr{Ad}(\mcE_n)^\lozenge)@> \Phi_{S, n}^*(\oint^\lozenge_{X_n, (\mcE_n, \nabla_{\mcE, n})}) >>  \Phi^*_{S, n} (\mcO_{S_n}) \ \left(\cong \mcO_{S_n} \right)
\\
@V F^\lozenge_{(\mcE_n, \nabla_{\mcE, n})} \otimes F^\lozenge_{(\mcE_n, \nabla_{\mcE, n})}   VV @VV [p^3] V
\\
\mbR^1f_{1, \mr{crys}*} (\mr{Ad}(\mcE_{n-1})^\lozenge) \otimes \mbR^1f_{1, \mr{crys}*} (\mr{Ad}(\mcE_{n-1})^\lozenge)  @>> \oint^\lozenge_{X_{n-1}, (\mcE_{n-1}, \nabla_{\mcE_{n-1}})} >
\mcO_{S_{n-1}}.
\end{CD}
\end{align}
The  diagram (\ref{G731}) may be obtained, via (\ref{G980}), as 
 the inverse limit (over $n$) of  the diagrams (\ref{G730}) in  the universal case (i.e., the case where the collection $(S, X, \mcE, \nabla_\mcE)$ is taken to be $(\mcN, C_\mcN, \mcE_\mcN, \nabla_{\mcE_\mcN})$).
 This implies the required commutativity, and completes the proof of the assertion.
\end{proof}
\vspace{3mm}


In what follows, we shall give a restatement of  the above proposition.
By passing the isomorphism $\mcT_{\widehat{\mcS}_{g, \mbZ_p}}|_{\mcN} \isom \mbR^1 f_{\mcN *}(\mcK^\bullet [\nabla^\mr{ad}_{\mcE_\mcN}])$ (cf. (\ref{EE001}) in  the universal case over $S = \mcN$), we obtain, from $F$, 
 an $\mcO_\mcN$-linear  morphism
\begin{align}
F^{\mr{PGL}} : \bigwedge^2 \Omega_{\widehat{\mcS}_{g, \mbZ_p}/\mbZ_p} |_{\mcN} \migi 
\left(\bigwedge^2 \Phi^*_{\mcN}(\Omega_{\widehat{\mcS}_{g, \mbZ_p}/\mbZ_p} |_{\mcN})  \cong \right) 
\Phi_\mcN^*(\bigwedge^2 \Omega_{\widehat{\mcS}_{g, \mbZ_p}/\mbZ_p} |_{\mcN}).   
\end{align}
Let us consider 
$\Gamma (\mcN, \bigwedge^2 \Omega_{\widehat{\mcS}_{g, \mbZ_p}/\mbZ_p} |_{\mcN})$
as a submodule of 
   $\Gamma (\mcN, \Phi^*_{\mcN}(\bigwedge^2 \Omega_{\widehat{\mcS}_{g, \mbZ_p}/\mbZ_p} |_{\mcN}))$ via  pull-back by $\Phi_\mcN$.
Then, Proposition \ref{P990} implies  the following assertion, which is essential to complete our proof  
of  Theorem \ref{T013}. (Theorem \ref{P0438} described in the next subsection  may be thought of as 
 another essential  ingredient of the proof.)

\vspace{3mm}
\bco \label{G970} \leavevmode\\
 \ \ \ 
 The following equality holds:
 \begin{align}
 F^{\mr{PGL}}(\widehat{\omega}^\mr{PGL}_{g, \mbZ_p} |_{\mcN} ) = p^3 \cdot \widehat{\omega}^\mr{PGL}_{g, \mbZ_p} |_{\mcN}.
\end{align}

 \eco

\vspace{5mm}
\subsection{Slope decomposition of $(\mbR^1 f_{\mcN*}(\mcK^\bullet [\nabla^\mr{ad}_{\mcE_\mcN}]), F)$} \label{SS079}
\leavevmode\\ \vspace{-4mm}

 Consider  the composite isomorphism
\begin{align} \label{G710}
\mcT_{\mcN/\mbZ_p} \stackrel{\sim}{\longmigi} \mcT_{\widehat{\mcM}_{g,\mbZ_p}/\mbZ_p}|_{\mcN} \stackrel{(\ref{G808})}{\longmigi}  \mbR^1f_{\mcN*}(\mcT_{C_{\mcN}/\mcN}),
\end{align}
where the first arrow arises from the \'{e}taleness of $\mcN/\widehat{\mcM}_{g, \mbZ_p}$.
It induces 
the composite isomorphism
\begin{align}  \label{eE19}
\Omega_{\mcN/\mbZ_p} \ \left(=\mcT_{\mcN/\mbZ_p}^\vee \right)
 \xrightarrow{((\ref{G710})^\vee)^{-1}} \mbR^1f_{\mcN*}(\mcT_{C_{\mcN}/\mcN})^\vee 
  \xrightarrow{\oint^{\natural -1}_{C_\mcN}} f_{\mcN*}(\Omega_{C_{\mcN}/\mcN}^{\otimes 2}).
\end{align}
Denote by
\begin{align} \label{G680}
\Upsilon : \mbR^1 f_{\mcN*}(\mcT_{C_{\mcN}/\mcN})\oplus f_{\mcN*}(\Omega^{\otimes 2}_{C_{\mcN}/\mcN}) \isom  \mbR^1 f_{\mcN*} (\mcK^\bullet [\nabla^\mr{ad}_{\mcE_{\mcN}}])
\end{align}
the unique isomorphism making the following  diagram commute:
\begin{align}
 \xymatrix@!C=140pt{
 0_{\mcN}^*(\mcT_{T^\vee_{\mbZ_p}\mcN/\mbZ_p})  \ar[d]^{\wr}_{d \theta |_{0_{\mcN}}}\ar[r]^{(\ref{e741})} &
  \mcT_{\mcN/\mbZ_p} \oplus \Omega_{\mcN/\mbZ_p} \ar[r]^{\hspace{-15mm}(\ref{G710}) \oplus (\ref{eE19})} &
     \mbR^1 f_{\mcN*}(\mcT_{C_{\mcN}/\mcN})\oplus f_{\mcN*}(\Omega^{\otimes 2}_{C_{\mcN}/\mcN})  \ar[d]_{\wr}^{\Upsilon}
\\
\mcT_{\widehat{\mcS}_{g, \mbZ_p}/\mbZ_p}|_{\mcN}  \ar[rr]_{(\ref{EE001})}^{\sim} & &  \mbR^1 f_{\mcN*} (\mcK^\bullet [\nabla^\mr{ad}_{\mcE_{\mcN}}]),
}
\end{align}
where 
the left-hand vertical arrow $d \theta |_{0_{\mcN}}$ denotes the differential of $\theta$ at $0_{\mcN}$ (i.e., the dual of  (\ref{G240}) restricted to $0_{\mcN}$).
That is to say, the direct sum decomposition $\Upsilon$  arises  from the 
classical ordinary $p$-adic Teichm\"{u}ller
uniformization
(cf. (\ref{G230}) or Introduction).
 
 We shall denote by  
\begin{align} \label{G400}
F^\sharp \ \left(:= p^2 \cdot \Phi_{\mcN}^*\right) : \Phi^*_{\mcN} (f_{\mcN*}(\Omega^{\otimes 2}_{C_{\mcN}/\mcN})) \migi f_{\mcN*}(\Omega^{\otimes 2}_{C_{\mcN}/\mcN})
\end{align}
the morphism defined as $p^2$ times the morphism 
$\Phi^*_{\mcN} : \Phi^*_{\mcN}(\Omega_{\mcN/\mbZ_p}) \migi \Omega_{\mcN/\mbZ_p}$ induced naturally by $\Phi_{\mcN}$
under the identification  $\Omega_{\mcN/\mbZ_p} \isom f_{\mcN*}(\Omega^{\otimes 2}_{C_{\mcN}/\mcN})$ (cf. (\ref{eE19})).
Here, notice that since the reduction modulo $p$ of $\Phi_{\mcN}$ coincides with  the Frobenius endomorphism,  $\Phi^*_{\widehat{\mcN}}$ is divisible by $p$.
%
According to ~\cite{Mzk1},  Chap.\,III, \S\,2, Proposition 2.3, the morphism $\frac{1}{p} \cdot \Phi^*_{\mcN}$ (i.e., $\Phi^*_{\mcN}$ divided by $p$) is an isomorphism.
Thus, we obtain a morphism 
\begin{align}
F^\flat : \Phi^*_{\mcN} (\mbR^1 f_{\mcN*}(\mcT_{C_\mcN/\mcN})) \isom\mbR^1 f_{\mcN*}(\mcT_{C_\mcN/\mcN})
\end{align}
defined to be the inverse to the dual of $\frac{1}{p} \cdot \Phi^*_{\widehat{\mcN}}$ under the identification 
$\mcT_{\mcN/\mbZ_p} \isom \mbR^1 f_{\mcN*}(\mcT_{C_\mcN/\mcN})$
 (cf. (\ref{G710})).
The Frobenius structure  $F$ on $\mbR^1 f_{\mcN*}(\mcK^\bullet [\nabla^\mr{ad}_{\mcE_{\mcN}}])$
  will turn out to be compatible with the Frbenius structures  
$F^\sharp$, $F^\flat$, 
 as   described below.
(The proof  will be given in \S\,\ref{SS499}.)

\vspace{3mm}
\bt \label{P0438} \leavevmode\\
 \vspace{-4mm}
\begin{itemize}
\item[(i)]
Let us consider the short exact sequence
\begin{align}
0 \longmigi f_{\mcN*}(\Omega^{\otimes 2}_{C_\mcN/\mcN})
\stackrel{\xi^\sharp_\mcN}{\longmigi}
\mbR^1 f_{\mcN*}(\mcK^\bullet [\nabla^\mr{ad}_{\mcE_{\mcN}}])
\stackrel{\xi^\flat_\mcN}{\longmigi}
\mbR^1  f_{\mcN*}(\mcT_{C_\mcN/\mcN})
\longmigi 0
\end{align}
defined as the inverse limit (over $n \geq 1$) of the sequence (\ref{G01}) of the case where $(\mcE, \nabla_\mcE)$ is taken to be $(\mcE_{\mcN, n}, \nabla_{\mcE_{\mcN, n}})$.
Then, the morphisms 
$F^\sharp$, $F^\flat$, and $F$
 are compatible with
the morphisms in this short exact sequence.
More precisely, the following diagram is commutative:
\begin{align} \label{W126}
\begin{CD}
\Phi^*_{\mcN}(f_{\mcN*}(\Omega^{\otimes 2}_{C_\mcN/\mcN}))
@>\Phi^*_{\mcN}(\xi_\mcN^\sharp) >> 
\Phi^*_{\mcN}(\mbR^1 f_{\mcN*}(\mcK^\bullet [\nabla^\mr{ad}_{\mcE_{\mcN}}]))
@> \Phi^*_{\mcN}(\xi^\flat_\mcN)>> \Phi^*_{\mcN}(\mbR^1  f_{\mcN*}(\mcT_{C_\mcN/\mcN}))
\\
@VV F^\sharp V @VV F V @VV F^\flat V 
\\
f_{\mcN*}(\Omega^{\otimes 2}_{C_\mcN/\mcN}) @>> \xi^\sharp_\mcN > \mbR^1 f_{\mcN*}(\mcK^\bullet [\nabla^\mr{ad}_{\mcE_{\mcN}}]) @>> \xi^\flat_\mcN > \mbR^1  f_{\mcN*}(\mcT_{C_\mcN/\mcN}).
\end{CD}
\end{align}

\item[(ii)]
The direct sum decomposition  $\Upsilon$ (cf. (\ref{G680})) is compatible with 
$F^\flat \oplus F^\sharp$ and $F$.
More precisely, the following square diagram is commutative:
\begin{align} \label{W125}
\begin{CD}
\Phi^*_{\mcN}(\mbR^1  f_{\mcN*}(\mcT_{C_\mcN/\mcN})) \oplus \Phi^*_{\mcN}(f_{\mcN*}(\Omega^{\otimes 2}_{C_\mcN/\mcN}))@> \Phi^*_{\mcN}(\Upsilon)> \sim > \Phi^*_{\mcN}(\mbR^1 f_{\mcN*}(\mcK^\bullet [\nabla^\mr{ad}_{\mcE_{\mcN}}]))
\\
@V F^\flat \oplus F^\sharp VV @VV F V
\\
 \mbR^1  f_{\mcN*}(\mcT_{C_\mcN/\mcN}) \oplus f_{\mcN*}(\Omega^{\otimes 2}_{C_\mcN/\mcN})  @> \sim > \Upsilon >  \mbR^1 f_{\mcN*}(\mcK^\bullet [\nabla^\mr{ad}_{\mcE_{\mcN}}]).
\end{CD}
\end{align}
\end{itemize}
 \et

\vspace{5mm}
\subsection{$F$-crystals associated to  canonical liftings of indigenous bundles} \label{SS0999}
\leavevmode\\ \vspace{-4mm}

In this subsection, we shall describe a consequence of  the above theorem obtained  by restricting  $\mbR^1 f_{\mcN*}(\mcK^\bullet [\nabla^\mr{ad}_{\mcE_\mcN}])$ to the fiber over each  Frobenius invariant  point in $\mcN$.
Let $k$ be an algebraically closed field of  characteristic $p$.
Write  $W$ for the ring of Witt vectors over  $k$ and 
write $\Phi_W$ for the absolute Frobenius automorphism of
$\mr{Spf}(W)$.
Let us take an arbitrary $k$-rational point $s_1 \in \mcN^\mr{ord}_{g, \mbF_p}(k)$.
Then, there exists a canonical  lifting $s_\infty  : \mr{Spf} (W) \migi \mcN$ of $s_1$  characterized uniquely by the equality   $\Phi_{\mcN} \circ s_\infty = s_\infty \circ \Phi_W$ (cf. ~\cite{Mzk1}, Chap.\,III, the discussion preceding Definition 1.9).
The point $s_\infty$ classifies a curve $X$ over $W$ and an  indigenous bundle $(\mcE, \nabla_\mcE)$ on it.
Denote by 
\begin{align}
\mbH^1 (\mcK^\bullet [\nabla^\mr{ad}_\mcE]) \ \left(= \varprojlim_n  \mbH^1 (\mcK^\bullet [\nabla^\mr{ad}_{\mcE, n}])\right)
\end{align}
 the first  hypercohomology associated to  the complex $\mcK^\bullet [\nabla^\mr{ad}_\mcE]$ (which may be  obtained by restricting $\mbR^1 f_{\mcN*}(\mcK^\bullet [\nabla^\mr{ad}_{\mcE_\mcN}])$ to $s_\infty$).
The isomorphism  $\Upsilon$ restricts to a direct sum decomposition
\begin{align} \label{G657}
\Upsilon_{(\mcE, \nabla_\mcE)} \ \left(:= s_\infty^*(\Upsilon) \right) : H^1 (X, \mcT_{X/W}) \oplus \Gamma (X, \Omega^{\otimes 2}_{X/W}) \isom \mbH^1 (\mcK^\bullet [\nabla^\mr{ad}_{\mcE}]).
\end{align}

Next, consider  the following sequence of isomorphisms
\begin{align} \label{G658}
s_\infty^*(\Phi^*_{\mcN}(\mbR^1 f_{\mcN*}(\mcK^\bullet [\nabla^\mr{ad}_{\mcE_{\mcN}}]))) &
\isom 
(\Phi_{\mcN} \circ s_\infty)^* (\mbR^1 f_{\mcN*}(\mcK^\bullet [\nabla^\mr{ad}_{\mcE_{\mcN}}])) \\
& \isom (s_\infty \circ \Phi_W)^* (\mbR^1 f_{\mcN*}(\mcK^\bullet [\nabla^\mr{ad}_{\mcE_{\mcN}}]))  \notag \\
& \isom  \Phi_W^*(s_\infty^*(\mbR^1 f_{\mcN*}(\mcK^\bullet [\nabla^\mr{ad}_{\mcE_{\mcN}}])) ). \notag 
\end{align}
By means of this composite isomorphism, the morphism $F$  restricts to a $W$-linear morphism
\begin{align} \label{G659}
F_{(\mcE, \nabla_\mcE)} \ \left(:= s_\infty^* (F)\right) : \Phi_W^*(\mbH^1 (\mcK^\bullet [\nabla^\mr{ad}_\mcE])) \migi \mbH^1 (\mcK^\bullet [\nabla^\mr{ad}_\mcE]).
\end{align}
Also, the morphism $F^\flat$ (resp., $F^\sharp$) restricts to a $W$-morphism
\begin{align}
F^\flat_{(\mcE, \nabla_\mcE)} : \Phi^*_W (H^1 (X, \mcT_{X/W})) \migi H^1 (X, \mcT_{X/W}) \\
  \left(\text{resp.,} \  F^\sharp_{(\mcE, \nabla_\mcE)} : \Phi^*_W (\Gamma (X, \Omega^{\otimes 2}_{X/W})) \migi \Gamma (X, \Omega_{X/W}^{\otimes 2})\right), \notag
\end{align}
for which the pair 
\begin{align}
(H^1 (X, \mcT_{X/W}), F^\flat_{(\mcE, \nabla_\mcE)}) \  \left(\text{resp.,}  \ 
(\Gamma (X, \Omega_{X/W}^{\otimes 2}), F^\sharp_{(\mcE, \nabla_\mcE)})\right)
\end{align}
 forms an {\it isoclinic $F$-crystal  over $k$ of rank $3g-3$ and its unique Newton slope is    $0$ (resp., $3$)}.
Then, Theorem \ref{P0438} asserted above implies the following assertion, which gives a geometric interpretation of the slope decomposition on $\mbH^1 (\mcK^\bullet [\nabla_\mcE^\mr{ad}])$ with respect to $F_{(\mcE, \nabla)}$.

\vspace{3mm}
\bco \label{C0338} \leavevmode\\
 \ \ \ 
 The pair 
 \begin{align}
 (\mbH^1 (\mcK^\bullet [\nabla^\mr{ad}_\mcE]), F_{(\mcE, \nabla_\mcE)})
 \end{align}
 forms  an $F$-crystal over $k$ of  rank $6g-6$ and all its   Newton slopes  are $0$ and $3$.
 Moreover, 
 the following equalities hold:
 \begin{align}
 \mbH^1 (\mcK^\bullet [\nabla^\mr{ad}_\mcE])^{F =p^0} = \Upsilon_{(\mcE, \nabla_\mcE)} (H^1 (X, \mcT_{X/W})), \hspace{5mm}
 \mbH^1 (\mcK^\bullet [\nabla^\mr{ad}_\mcE])^{F =p^3} = \Upsilon_{(\mcE, \nabla_\mcE)} (\Gamma (X, \Omega^{\otimes 2}_{X/W})),
 \end{align}
where $\mbH^1 (\mcK^\bullet [\nabla^\mr{ad}_\mcE])^{F =p^m}$ (for an  integer $m$) denotes
the isoclinic component of  $\mbH^1 (\mcK^\bullet [\nabla^\mr{ad}_\mcE])$ of  slope $m$ (with respect to  $F_{(\mcE, \nabla)}$).
In particular,  $\Upsilon_{(\mcE, \nabla_\mcE)}$ may be regarded as the slope decomposition of $(\mbH^1 (\mcK^\bullet [\nabla^\mr{ad}_\mcE]), F_{(\mcE, \nabla_\mcE)})$.
 
\eco

\vspace{3mm}
\begin{rema} \label{R012}\leavevmode\\
 \ \ \ 
In this remark, we shall describe $F^\flat_{(\mcE, \nabla_\mcE)}$ in terms of deformations of $(\mcE, \nabla_\mcE)$.
 Let us keep the above notation.
Denote by  $ \mcT_{\mcN/\mbZ_p, s_\infty}$  the tangent space of $\mcN$ (over $\mbZ_p$) at $s_\infty$, which 
 may be identified with the 
 deformation space of $(\mcE, \nabla_\mcE)$ over $W_\epsilon := W [\epsilon]/(\epsilon^2)$.
 Let us consider the  differential 
 \begin{align} \label{W200}
 d \Phi_\mcN |_{s_\infty} : \mcT_{\mcN/\mbZ_p, s_\infty} \migi \left(\Phi^*_\mcN(\mcT_{\mcN/\mbZ_p, s_\infty}) =\right) \ \Phi^*_W(\mcT_{\mcN/\mbZ_p, s_\infty})
 \end{align}
  of $\Phi_\mcN$ at $s_\infty$.
Under the identification 
$ \mcT_{\mcN/\mbZ_p, s_\infty} \isom H^1 (X, \mcT_{X/W})$ obtained by restricting  (\ref{G710}),
the equality $(F^\flat_{(\mcE, \nabla_\mcE)})^{-1} = \frac{1}{p} \cdot  d \Phi_\mcN |_{s_\infty}$ holds.
In particular, we have
 \begin{align} \label{W201}
\mr{Im}(d \Phi_\mcN |_{s_\infty}) = p\cdot \Phi^*_W(\mcT_{\mcN/\mbZ_p, s_\infty}) = \Phi^*_W(p \mcT_{\mcN/\mbZ_p, s_\infty}) \ \left(\subseteq \Phi^*_W(\mcT_{\mcN/\mbZ_p, s_\infty}) \right),
\end{align}
and  the morphism $\mcT_{\mcN/\mbZ_p, s_\infty} \migi \Phi^*_W(p \mcT_{\mcN/\mbZ_p, s_\infty})$ induced by
$d \Phi_\mcN |_{s_\infty}$ is an isomorphism.
 
 Now, 
 let us 
 take  an element 
 $v \in \mcT_{\mcN/\mbZ_p, s_\infty}$ (resp., $v' \in p\mcT_{\mcN/\mbZ_p, s_\infty}$).
 By the above discussion, one may find a unique $\widetilde{v}' \in \mcT_{\mcN/\mbZ_p, s_\infty}$ with
 $d \Phi_\mcN |_{s_\infty} (\widetilde{v}') = v'$, or equivalently, $\Phi_\mcN \circ \widetilde{v}' = v' \circ \Phi_{W_\epsilon}$, where $\Phi_{W_\epsilon}$ denotes
  the base-change of $\Phi_W$ over $W_\epsilon$.
Denote by $(X_\epsilon^v, \mcE^v_\epsilon, \nabla^v_{\mcE, \epsilon})$,  $(X^{v'}_\epsilon, \mcE^{v'}_{\epsilon}, \nabla^{v'}_{\mcE, \epsilon})$, and $(X^{\widetilde{v}'}_\epsilon, \mcE^{\widetilde{v}'}_{\epsilon}, \nabla^{\widetilde{v}'}_{\mcE, \epsilon})$  the collections classified by $v$, $v'$, and $\widetilde{v}'$ respectively.
In particular, $(X_\epsilon^{\overline{v}'})_1 \cong (X_\epsilon^{v'})_1^{(1)}$.
Then, we have the following sequence of isomorphisms  over $X_\epsilon^{\overline{v}'}$:
\begin{align}
(\mcE^{\widetilde{v}'}_{\epsilon}, \nabla^{\widetilde{v}'}_{\mcE, \epsilon})
&\cong \widetilde{v}'^*(\mcE_\mcN, \nabla_{\mcE_\mcN})  \\
&\cong
\widetilde{v}'^* (\mbF^*(\Phi_\mcN^*(\mcE_\mcN, \nabla_{\mcE_\mcN})))  \notag \\
&\cong
\mbF^*((\Phi_\mcN \circ \widetilde{v}')^*(\mcE_\mcN, \nabla_{\mcE_\mcN})) \notag \\
&\cong
\mbF^*((v' \circ \Phi_{W_\epsilon})^*(\mcE_\mcN, \nabla_{\mcE_\mcN})) \notag \\
& \cong
\mbF^*(\Phi_{W_\epsilon}^* (v'^*(\mcE_\mcN, \nabla_{\mcE_\mcN}))) \notag \\
&\cong
\mbF^*(\Phi_{W_\epsilon}^*(\mcE^{v'}_{\epsilon}, \nabla^{v'}_{\mcE, \epsilon})). \notag
\end{align}
It follows consequently that 
  {\it $\left(p \cdot (F^\flat_{(\mcE, \nabla_\mcE)})^{-1} (v) =\right) \ d\Phi_\mcN |_{s_\infty} (v) = v'$ (or equivalently, $v = \widetilde{v}'$) if and only if
$(\mcE^v_\epsilon, \nabla^v_{\mcE, \epsilon}) \cong  \mbF^*(\Phi_{W_\epsilon}^*(\mcE^{v'}_{\epsilon}, \nabla^{v'}_{\mcE, \epsilon}))^\lozenge_{X^v_\epsilon}$}. 
\end{rema}

\vspace{5mm}
\subsection{Relationship between $p$-curvature and the differential of $\sigma$} \label{SS0919}
\leavevmode\\ \vspace{-4mm}

In this subsection, we shall prove a proposition (cf. Proposition \ref{p452} asserted below), which describes 
the reduction modulo $p$ of the  differential of $\sigma$ by means of $p$-curvature.
That proposition   will be used in the proof of Theorem \ref{T013} (cf. the proof of Lemma \ref{p401}).

Let $k$ be as before and $s_1$ a $k$-rational point of $\mcN^\mr{ord}_{g, \mbF_p}$, and $(X, \mcE, \nabla_\mcE)$ the data classified by $s_1$.
The morphism   $\psi_{(\mcE, \nabla_{\mcE})}^\nabla$ (cf. (\ref{W118})) induces a morphism of complexes $\Phi^{-1}_{X/k}(\mcT_{X^{(1)}/k})[0] \migi  \mcK^\bullet [\nabla^\mr{ad}_{\mcE}]$ (cf. \S\,\ref{SS983} for the definition of $(-)[0]$).
By applying the functor $\mbH^1 (-)$ to this morphism, 
we obtain a morphism
\begin{align} \label{W120}
\mbH^1(\psi_{(\mcE, \nabla_{\mcE})}^\nabla) : H^1 (X^{(1)}, \mcT_{X^{(1)}/k}) \migi 
 \mbH^1 (\mcK^\bullet [\nabla^\mr{ad}_{\mcE}]).
\end{align}
On the other hand, 
let us consider the composite  
\begin{align}
d \sigma |_{s_1} : H^1 (X, \mcT_{X/k}) \migiincl  H^1 (X, \mcT_{X/k}) \oplus \Gamma (X, \Omega_{X/k}^{\otimes 2}) 
\stackrel{\sim}{\longmigi}
\mbH^1 (\mcK^\bullet [\nabla_{\mcE}^\mr{ad}]),
\end{align}
 where the first arrow denotes the inclusion into the first factor and the second arrow denotes the 
 restriction to $s_1$ of the direct sum decomposition (\ref{G680}).
 It   may also be  obtained as the differential of the immersion $\mcN^\mr{ord}_{g, \mbF_p} \migi \mcS_{g, \mbF_p}$ at $s_1$
 under the identifications (\ref{EE001}) and  (the reduction modulo $p$ of) (\ref{G710}). 
 Then, the following lemma holds.

\vspace{3mm}
\bpr \label{p452}\leavevmode\\
 \ \ \  
The following square diagram is commutative:
\begin{align}
\begin{CD}
\Phi^*_k (H^1 (X, \mcT_{X/k}))@> \varsigma >\sim > H^1(X^{(1)}, \mcT_{X^{(1)}/k})
\\
@V - F_{(\mcE, \nabla_{\mcE})}^\flat V \wr V @VV\mbH^1(\psi_{(\mcE, \nabla_{\mcE})}^\nabla) V
\\
H^1(X, \mcT_{X/k})@>> d \sigma |_{s_1}> \mbH^1 (\mcK^\bullet [\nabla^\mr{ad}_{\mcE}]),
\end{CD}
\end{align}
where $\varsigma$ denotes the natural isomorphism induced by $\mr{id}_X \times \Phi_k :X^{(1)} \migi X$ 
and $F^\flat_{(\mcE, \nabla_{\mcE})}$ denotes  the 
restriction of $F^\flat$ to $s_1$.
In particular, $\mbH^1(\psi_{(\mcE_1, \nabla_{\mcE,1})}^\nabla)$ is injective.
 \epr
\begin{proof}
By definition, the composite 
\begin{align}
\xi^\flat \circ \mbH^1(\psi_{(\mcE, \nabla_{\mcE})}^\nabla) \circ \varsigma : \Phi^*_k (H^1 (X, \mcT_{X/k})) \migi H^1 (X, \mcT_{X/k})
\end{align}
coincides with the morphism (\ref{G02}) in the present  case.
Hence, it follows from ~\cite{Mzk1}, Chap.\,III, \S\,2, Proposition 2.3,
that  
\begin{align}
\xi^\flat \circ ( \mbH^1(\psi_{(\mcE_1, \nabla_{\mcE, 1})}^\nabla) \circ \varsigma) =  - F_{(\mcE_1, \nabla_{\mcE, 1})}^\flat \ \left(= - \xi^\flat \circ (d \sigma |_{s_1} \circ F_{(\mcE_1, \nabla_{\mcE, 1})}^\flat) \right).
\end{align}
This implies that, to complete the proof of the assertion,   it suffices to prove the equality 
\begin{align}
\mr{Im} (\mbH^1(\psi_{(\mcE, \nabla_{\mcE})}^\nabla)) = \mr{Im} (d \sigma |_{s_1}).
\end{align}

First, let $v$ be an element of $\mr{Im} (\mbH^1(\psi_{(\mcE, \nabla_{\mcE})}^\nabla))$, and 
denote by $(X_\epsilon^v, \mcE^v_\epsilon, \nabla^v_{\mcE, \epsilon})$ 
 the deformation over $k_\epsilon := k [\epsilon]/(\epsilon^2)$ of $(X, \mcE, \nabla_\mcE)$ corresponding to $(\eta^\mr{ad}_B)^{-1} (v) \in \mbH^1 (\mcK^\bullet [\widetilde{\nabla}^\mr{ad}_{\mcE_B}])$ (cf. (\ref{W356}) for the definition of $\eta^\mr{ad}_B$).
By the definition of $\psi^\nabla_{(\mcE, \nabla_\mcE)}$,  $v$ may be represented as the data  (\ref{W105}) with $b_\alpha = 0$ (for any $\alpha$).
It follows from the construction of the bijection (\ref{W100}) that
$(\mcE^v_\epsilon, \nabla^v_{\mcE, \epsilon})$  is,
 locally on $X$,
isomorphic to the trivial deformation.
Hence, it has nilpotent $p$-curvature, and hence, forms a nilpotent  indigenous bundle.
According to  ~\cite{Mzk1},  Chap.\,II, \S\,2, Theorem 2.13,  $\mcN^\mr{ord}_{g, \mbF_p}$ coincides with the \'{e}tale locus (relative to  $\mcM_{g, \mbF_p}$) in the substack of $\mcS_{g, \mbF_p}$ classifying  nilpotent indigenous bundles.
This implies  that  the element $v$,  when considered as an element of $\mcS_{g, \mbF_p} (k_\epsilon)$, lies in $\mcN^\mr{ord}_{g, \mbF_p} (k_\epsilon)$.
That is to say, $v$ is contained in $\mr{Im}(d \sigma |_{s_1})$.

Conversely, let $u$ be an element of 
$\mr{Im}(d \sigma |_{s_1})$ and denote by $(X_\epsilon^u, \mcE^n_\epsilon, \nabla^u_{\mcE, \epsilon})$ the deformation of $(X, \mcE, \nabla_\mcE)$ over $k_\epsilon$ determined by $(\eta^\mr{ad}_B)^{-1}(u) \in \mbH^1 (\mcK^\bullet [\widetilde{\nabla}^\mr{ad}_{\mcE_{1, B}}])$.
Recall here  the notion of an  FL-(vector) bundle introduced in  ~\cite{Mzk1}, Chap.\,II, \S\,1, Definition 1.3. 
Denote by $(\mcG, \nabla_\mcG)$ (resp.,  $(\mcG_\epsilon^u, \nabla^u_{\mcG, \epsilon})$) the FL-bundle on $X$ (resp.,  $X_\epsilon^u$)
corresponding to $(\mcE, \nabla_\mcE)$ (resp., $(\mcE^n_\epsilon, \nabla^u_{\mcE, \epsilon})$) (cf. ~\cite{Mzk1}, Chap.\,II, \S\,2, Proposition 2.5).
The flat bundle $(\mcG^u_\epsilon, \nabla^u_{\mcG, \epsilon})$ forms  
an extension of $(\mcO_{X_\epsilon^u}, d)$ (i.e., the trivial flat bundle) by $(\Phi^*_{X_{\epsilon}^u/k_\epsilon}(\mcT_{(X^{u}_{\epsilon})^{(1)}/k_\epsilon}), \nabla^{\mr{can}}_{(X^{u}_{\epsilon})^{(1)}})$ (cf. (\ref{W110}) for the definition of $\nabla^\mr{can}_{(-)}$). 
Let us consider the short exact sequence
\begin{align} \label{W90}
0 \longmigi H^1 ((X_\epsilon^u)^{(1)}, \mcT_{(X_\epsilon^u)^{(1)}/k_\epsilon}) \stackrel{\nu}{\longmigi} 
\mbH^1 (\mcK^\bullet [\nabla^{\mr{can}}_{(X^{u}_{\epsilon})^{(1)}}]) 
\stackrel{\mu}{\longmigi} \Gamma (X_\epsilon^u, \mcO_{X_\epsilon}) \ \left(= k_\epsilon \right) \longmigi 0
\end{align}
discussed in ~\cite{Mzk1}, Chap.\,II, \S\,1, Proposition 1.1, in the case where the curve ``$X^\mr{log}/S^\mr{log}$" in {\it loc. cit.} is taken to be  $X_\epsilon^u/k_\epsilon$;
it is obtained by applying the functor $\mbH^1 (-)$ to the following short exact sequence of complexes:
\begin{align} \label{W91}
0 \longmigi \Phi^{-1}_{X_{\epsilon}^u/k_\epsilon}(\mcT_{(X^{u}_{\epsilon})^{(1)}/k_\epsilon})[0]
\longmigi  \mcK^\bullet [\nabla^{\mr{can}}_{(X^{u}_{\epsilon})^{(1)}}]
\longmigi \Phi^{-1}_{X_{\epsilon}^u/k_\epsilon}(\mcO_{(X_\epsilon^u)^{(1)}})[1]
\longmigi 0
\end{align}
(cf. \S\,\ref{SS983} for the definition of $(-)[n]$), where the second arrow arises from the inclusion $\Phi^{-1}_{X_{\epsilon}^u/k_\epsilon}(\mcT_{(X^{u}_{\epsilon})^{(1)}/k_\epsilon}) \migiincl \Phi^*_{X_{\epsilon}^u/k_\epsilon}(\mcT_{(X^{u}_{\epsilon})^{(1)}/k_\epsilon})$ and the third arrow arises from
the Cartier operator 
\begin{align}
\Omega_{X_\epsilon^u/k_\epsilon}\otimes \Phi_{X_\epsilon^u/k_\epsilon}^*(\mcT_{X_\epsilon^u/k_\epsilon}) \migisurj \Phi_{X_\epsilon^u/k_\epsilon }^{-1}(\Omega_{(X_\epsilon^u)^{(1)}/k_\epsilon} \otimes \mcT_{(X_\epsilon^u)^{(1)}/k_\epsilon}) \ \left(= \Phi_{X_\epsilon^u/k_\epsilon }^{-1}(\mcO_{(X_\epsilon^u)^{(1)}}) \right)
\end{align}
 (cf. ~\cite{Og}, Proposition 1.2.4 together with   the discussion following that proposition).
If $\mr{ex}_{(\mcG_\epsilon^u, \nabla^u_{\mcG, \epsilon})}$ $\left(\in \mbH^1 (\mcK^\bullet [\nabla^{\mr{can}}_{(X^{u}_{\epsilon})^{(1)}}]) \right)$ denotes the extension class determined by
$(\mcG_\epsilon^u, \nabla^u_{\mcG, \epsilon})$,
then it follows from the definition of an FL-bundle that $\mu (\mr{ex}_{(\mcG_\epsilon^u, \nabla^u_{\mcG, \epsilon})} ) \in k_\epsilon^{\times}$.
Also, let us denote by $\mr{ex}_{(\mcG_\epsilon, \nabla_{\mcG, \epsilon})}$ $\left(\in \mbH^1 (\mcK^\bullet [\nabla^\mr{can}_{X_\epsilon^{(1)}}]) \right)$ the extension class determined by  the trivial deformation $(\mcG_\epsilon, \nabla_{\mcG, \epsilon})$.

In what follows, we shall denote the base-changes to $k_\epsilon$ of objects over $k$ by means of a subscripted $\epsilon$.
 Let us
 take an affine open covering $\mfU := \{ U_\alpha \}_{\alpha \in I}$ of $X$.
In the \v{C}ech double complex $\mr{Tot}^\bullet (\check{C} (\mfU, \mcK^\bullet [\nabla^\mr{ad}_\mcE]))$,  the element $u$ may be represented by 
   $(\{ a_{\alpha \beta} \}_{\alpha \beta}, \{ b_\alpha \}_\alpha)$ as in (\ref{W105}).
Since $X_\epsilon^u$ and $\mcE_\epsilon^u$ may be obtained by gluing together $U_{\alpha, \epsilon}$'s and $\mcE_{\epsilon} |_{U_\alpha}$'s respectively,
there exist  natural isomorphisms $\iota_{X, \alpha} : U_{\alpha, \epsilon} \isom X^u_\epsilon |_{U_\alpha}$ and $\iota_{\mcE, \alpha} : \mcE_{\epsilon} |_{U_\alpha} \isom \mcE^u_{\epsilon} |_{U_\alpha}$ (for each $\alpha \in I$) respectively, whose reductions modulo $\epsilon$ are the identity morphisms.
Under the isomorphism $\iota_{\mcE, \alpha}$, the restriction $(\mcE^u_{\epsilon} |_{U_\alpha}, \nabla^u_{\epsilon} |_{U_\alpha})$ may be identified with $(\mcE_\epsilon |_{U_\alpha}, \nabla_{\mcE, \epsilon} + \epsilon \cdot b_\alpha)$.

 Let us fix  $\alpha \in I$, and consider  the following morphism of short exact sequences induced by restriction via $U_{\alpha, \epsilon} \xrightarrow{\iota_{X, \alpha}}X^u_{\epsilon}|_{U_\alpha} \migiincl X^u_\epsilon$:
 \begin{align} \label{W94}
 \begin{CD}
 0 @>>>  H^1 ((X_\epsilon^u)^{(1)}, \mcT_{(X_\epsilon^u)^{(1)}/k_\epsilon})  @> \nu >>  
 \mbH^1 (\mcK^\bullet [\nabla^{\mr{can}}_{(X^u_\epsilon)^{(1)}}]) 
 @> \mu >> k_\epsilon @>>> 0
 \\
 @. @VVV @VV (-) |_{U_\alpha}V @VVV @.
 \\
 0 @>>> H^1 (U_{\alpha, \epsilon}, \mcT_{U_{\alpha, \epsilon}^{(1)}/k_\epsilon})
 @>>\nu_\alpha>
 \mbH^1 (\mcK^\bullet [\nabla^\mr{can}_{U_{\alpha, \epsilon}^{(1)}}])
 @>>\mu_\alpha >
 \Gamma (U_{\alpha, \epsilon}^{(1)}, \mcO_{U_{\alpha, \epsilon}^{(1)}})
 @>>> 0,
 \end{CD}
 \end{align}
where the lower horizontal arrow is
(\ref{W90}) (i.e., the upper horizontal sequence) with  $X_\epsilon^u$ and $\nabla^\mr{can}_{(X^u_\epsilon)^{(1)}}$ replaced by $U_{\alpha, \epsilon}$ and   $\nabla^\mr{can}_{U_{\alpha, \epsilon}^{(1)}}$ respectively.
Since $\nu_\alpha =0$ (which implies that $\mu_\alpha$ is an isomorphism),
$k_\epsilon \ \left(\subseteq  \Gamma (U_{\alpha, \epsilon}^{(1)}, \mcO_{U_{\alpha, \epsilon}^{(1)}}) \right)$ may be thought, via $\mu_\alpha$, of as a submodule of 
$ \mbH^1 (\mcK^\bullet [\nabla^\mr{can}_{U_{\alpha, \epsilon}^{(1)}}])$.
The restriction  $\mr{ex}_{(\mcG_\epsilon^u, \nabla^u_{\mcG, \epsilon})} |_{U_\alpha} \ \left(\in  \mbH^1 (\mcK^\bullet [\nabla^\mr{can}_{U_{\alpha, \epsilon}^{(1)}}])\right)$ lies in $k_\epsilon^\times$.
On the other hand, the restriction $\mr{ex}_{(\mcG_\epsilon, \nabla_{\mcG, \epsilon})} |_{U_\alpha} \ \left(\in  \mbH^1 (\mcK^\bullet [\nabla^\mr{can}_{U_{\alpha, \epsilon}^{(1)}}]) \right)$ of
$\mr{ex}_{(\mcG_\epsilon, \nabla_{\mcG, \epsilon})}$
to $U_\alpha$
   lies in $k^\times$.
Hence, 
$\mr{ex}_{(\mcG_\epsilon^u, \nabla^u_{\mcG, \epsilon})} |_{U_\alpha}$
and
$\mr{ex}_{(\mcG_\epsilon, \nabla_{\mcG, \epsilon})} |_{U_\alpha}$ 
differ at most by a  constant factor in $k_\epsilon^\times$, which implies that
the restrictions $(\mcG_\epsilon^u |_{U_\alpha}, \nabla^u_{\mcG, \epsilon} |_{U_\alpha})$
and 
$(\mcG_\epsilon  |_{U_\alpha}, \nabla_{\mcG, \epsilon} |_{U_\alpha})$
are isomorphic.
The resulting  isomorphism
$(\mcE_\epsilon |_{U_\alpha}, \nabla_{\mcE, \epsilon}) \isom (\mcE_\epsilon |_{U_\alpha}, \nabla_{\mcE, \epsilon} + \epsilon \cdot b_\alpha) \ \left(\stackrel{\iota_{\mcE, \alpha}}{\cong}  (\mcE^u_\epsilon |_{U_\alpha}, \nabla^u_\epsilon |_{U_\alpha})\right)$ may be expressed as $\mr{id}_{\mcE_\epsilon |_{U_\alpha}} + \epsilon \cdot c_\alpha$ for some $c_\alpha \in \Gamma (U_\alpha, \mr{Ad}(\mcE))$.
By the definition of $\nabla^\mr{ad}_{\mcE}$, the equality $\nabla^\mr{ad}_\mcE (c_\alpha) = b_\alpha$ holds.
 Hence,
 if  $\nabla^{\mr{ad}, \mr{Im}}_\mcE$  denotes the morphism $\mr{Ad}(\mcE) \migi \mr{Im}(\nabla^\mr{ad}_\mcE)$ obtained from $\nabla^\mr{ad}_\mcE$  by restricting   its codomain, then
  $u$ lies in the image of the morphism $\mbH^1 (\mcK^\bullet [\nabla^{\mr{ad}, \mr{Im}}_\mcE]) \migi \mbH^1 (\mcK^\bullet [\nabla^\mr{ad}_\mcE])$ induced by the natural injection
$\mcK^\bullet [\nabla^{\mr{ad}, \mr{Im}}_\mcE] \migiincl \mcK^\bullet [\nabla^{\mr{ad}}_\mcE]$.
 Observe that the morphism 
  $H^1 (X, \mr{Ker}(\nabla^\mr{ad}_\mcE)) \migi \mbH^1 (\mcK^\bullet [\nabla^{\mr{ad}, \mr{Im}}_\mcE])$ induced by the natural morphism $\mr{Ker}(\nabla^\mr{ad}_\mcE)[0] \migi \mcK^\bullet [\nabla_\mcE^{\mr{ad}, \mr{Im}}]$ is an isomorphism. 
Moreover,   it follows from   ~\cite{Mzk1}, Chap.\,II, \S\,2,  Proposition 2.7,
that  $\psi^\nabla_{(\mcE, \nabla_\mcE)}$ restricts to an isomorphism 
$\Phi^{-1}_{X/k}(\mcT_{X^{(1)}/k}) \isom \mr{Ker}(\nabla^\mr{ad}_{\mcE})$, which induces an isomorphism $H^1 (X^{(1)}, \mcT_{X^{(1)}/k}) \isom  H^1 (X, \mr{Ker}(\nabla^\mr{ad}_\mcE))$.
Hence, $u$ lies in the image of the composite injection
\begin{align}
H^1 (X^{(1)}, \mcT_{X^{(1)}/k}) \isom
 H^1 (X, \mr{Ker}(\nabla^\mr{ad}_\mcE))
\isom
\mbH^1 (\mcK^\bullet [\nabla^{\mr{ad}, \mr{Im}}_\mcE])
\migiincl
\mbH^1 (\mcK^\bullet [\nabla^{\mr{ad}}_\mcE]).
\end{align}
Since this composite  is nothing but $\mbH^1 (\psi^\nabla_{(\mcE, \nabla_{\mcE})})$,
we have  $u \in \mr{Im}(\mbH^1 (\psi^\nabla_{(\mcE, \nabla_{\mcE})}))$.
 This completes the proof of the assertion.
\end{proof}


\vspace{10mm}
\section{Proofs of theorems} \vspace{3mm}

This section is devoted to prove Theorem \ref{T013}  
 and Theorem \ref{P0438} 
  described earlier.

\vspace{5mm}
\subsection{Proof of Theorem \ref{P0438}} \label{SS499}
\leavevmode\\ \vspace{-4mm}

First, we prove Theorem  \ref{P0438}.
Let $k$, $W$, $s_1$, $s_\infty$, $X$, and  $(\mcE, \nabla_\mcE)$ be as in \S\,\ref{SS0999}.
By considering various $s_1$, we see that, in order to 
prove assertions  (i) and (ii),
it suffices to verify  the required commutativities  of the diagrams (\ref{W126}), (\ref{W125})   restricted  to $s_\infty \in \mcN (W)$.
Write  $\pi_B : \mcE_B \migi X$  for the Hodge reduction of $(\mcE, \nabla_\mcE)$, and write  
 $\widetilde{\nabla}^\mr{ad}_{\mcE_B} := \varprojlim_n \widetilde{\nabla}^\mr{ad}_{\mcE_{B, n}}$,
$\nabla^\mr{ad}_{\mcE} := \varprojlim_n \nabla^\mr{ad}_{\mcE_{n}}$.
Denote by $\Phi_{W_\epsilon}$ the base-change of $\Phi_W$ over $W_\epsilon  := W[\epsilon]/(\epsilon^2)$.
Let  $(\mcE^\dagger_\epsilon, \nabla^\dagger_{\mcE, \epsilon})$ be a deformation of  $(\mcE, \nabla_\mcE)$ over $X_\epsilon := X \times_W W_\epsilon$  classified by  $p  \mbH^1 (\mcK^\bullet [\nabla^\mr{ad}_\mcE])$ $\left(\subseteq  \mbH^1 (\mcK^\bullet [\nabla^\mr{ad}_\mcE]) \right)$.
Then, (since $\mbF^*(\Phi_W^*(\mcE, \nabla_\mcE)) \cong (\mcE, \nabla_\mcE)$) the flat $\mr{PGL}_2$-torsor  $\mbF^*(\Phi^*_{W_\epsilon}(\mcE^\dagger_\epsilon, \nabla^\dagger_{\mcE, \epsilon}))$ over $X_\epsilon/W_\epsilon$
forms a deformation  of $(\mcE, \nabla_\mcE)$, i.e., specifies an element of $\mbH^1 (\mcK^\bullet [\nabla^\mr{ad}_\mcE])$.
The assignment $(\mcE^\dagger_\epsilon, \nabla^\dagger_{\mcE, \epsilon}) \mapsto \mbF^*(\Phi_{W_\epsilon}^*(\mcE^\dagger_\epsilon, \nabla^\dagger_{\mcE, \epsilon}))$ defines a $W$-linear map
\begin{align} \label{W191}
\phi_\mcE : \Phi^*_W (p \mbH^1 (\mcK^\bullet [\nabla^\mr{ad}_\mcE])) \migi \mbH^1 (\mcK^\bullet [\nabla^\mr{ad}_\mcE]).
\end{align}
On the other hand, let $(X^\ddagger_\epsilon, \mcE^\ddagger_\epsilon, \nabla^\ddagger_{\mcE, \epsilon})$ be the deformation of  $(X, \mcE, \nabla_\mcE)$ over $W_\epsilon$  classified by  $p  \mbH^1 (\mcK^\bullet [\widetilde{\nabla}^\mr{ad}_{\mcE_B}]) \ \left(\subseteq  \mbH^1 (\mcK^\bullet [\widetilde{\nabla}^\mr{ad}_{\mcE_B}]) \right)$.
According to  Corollary \ref{W150}, there exists a unique deformation 
$(X_\epsilon^{\ddagger '}, \mcE^{\ddagger '}_\epsilon, \nabla^{\ddagger '}_{\mcE, \epsilon})$ classified by 
$\mbH^1 (\mcK^\bullet [\widetilde{\nabla}^\mr{ad}_{\mcE_B}])$ such that
$(\mcE^{\ddagger '}_\epsilon, \nabla^{\ddagger '}_{\mcE, \epsilon}) \cong \mbF^*(\Phi_{W_\epsilon}^*(\mcE^\dagger_\epsilon, \nabla^\dagger_{\mcE, \epsilon}))^\lozenge_{X_\epsilon^{\ddagger '}}$.
The assignment $(X^\ddagger_\epsilon, \mcE^\ddagger_\epsilon, \nabla^\ddagger_{\mcE, \epsilon}) \mapsto (X_\epsilon^{\ddagger '}, \mcE^{\ddagger '}_\epsilon, \nabla^{\ddagger '}_{\mcE, \epsilon})$ defines a $W$-linear map
\begin{align} \label{W190}
\widetilde{\phi}_{\mcE_B} : \Phi^*_W (p \mbH^1 (\mcK^\bullet [\widetilde{\nabla}^\mr{ad}_{\mcE_B}])) \migi \mbH^1 (\mcK^\bullet [\widetilde{\nabla}^\mr{ad}_{\mcE_B}]).
\end{align}
By construction, the following square diagram is commutative:
\begin{align}
\begin{CD}
\Phi^*_W (p \mbH^1 (\mcK^\bullet [\widetilde{\nabla}^\mr{ad}_{\mcE_B}])) 
@> \widetilde{\phi}_{\mcE_B} >>
\mbH^1 (\mcK^\bullet [\widetilde{\nabla}^\mr{ad}_{\mcE_B}])
\\
@V \eta_B^\mr{ad} V \wr V @V \wr V \eta^\mr{ad}_B V
\\
\Phi^*_W (p \mbH^1 (\mcK^\bullet [\nabla^\mr{ad}_\mcE])) 
@>> \phi_\mcE >
\mbH^1 (\mcK^\bullet [\nabla^\mr{ad}_\mcE]).
\end{CD}
\end{align}
Observe that  the composite 
\begin{align}
H^1 (X, \mcT_{X/W}) \migiincl H^1 (X, \mcT_{X/W})\oplus \Gamma (X, \Omega_{X/W}^{\otimes 2}) \stackrel{\Upsilon_{(\mcE, \nabla_\mcE)}}{\stackrel{\sim}{\longmigi}}
\mbH^1 (\mcK^\bullet [\nabla_\mcE^\mr{ad}]) 
\stackrel{(\eta_B^{\mr{ad}})^{-1}}{\stackrel{\sim}{\longmigi}}
\mbH^1 (\mcK^\bullet[\widetilde{\nabla}^\mr{ad}_{\mcE_B}])
\end{align}
coincides (via (\ref{G710}) and (\ref{D003})) with the differential of $\sigma :\mcN \migi \widehat{\mcS}_{g, \mbZ_p}$ at the point $s_\infty$.
If we consider $H^1 (X, \mcT_{X/W})$ as a submodule of $\mbH^1 (\mcK^\bullet[\widetilde{\nabla}^\mr{ad}_{\mcE_B}])$ via this composite,
then,  by the definition of $\widetilde{\phi}_{\mcE_B}$, we have $\widetilde{\phi}_{\mcE_B} (\Phi^{-1}_W(pH^1 (X, \mcT_{X/W}))) \subseteq H^1 (X, \mcT_{X/W})$.
Moreover, it follows from the discussion  in Remark \ref{R012}   that 
the restriction $\Phi^*_W(pH^1 (X, \mcT_{X/W})) \migi H^1 (X, \mcT_{X/W})$ of $\widetilde{\phi}_{\mcE_B}$ coincides  with the inverse of the differential $d \Phi_{\mcN} |_{s_\infty}$ (cf. (\ref{W200})) of $\Phi_\mcN$ at $s_\infty$.
That is to say, the equality 
\begin{align}
F^\flat_{(\mcE, \nabla_\mcE)} =\widetilde{\phi}_{\mcE_B}\circ \Phi^*_W([p]) 
\end{align}
holds, where, for each $W$-module $H$, 
  we shall write   $[p]$ for the morphism $H \migi pH \ \left(\subseteq H \right)$  given by multiplication by $p$.
Also, 
if we consider $\Gamma (X, \Omega^{\otimes 2}_{X/W})$
as a submodule of $\mbH^1 (\mcK^\bullet [\widetilde{\nabla}^\mr{ad}_\mcE])$ via the composite injection 
$(\eta_B^\mr{ad})^{-1}\circ \xi^\sharp : \Gamma (X, \Omega^{\otimes 2}_{X/W}) \migiincl \mbH^1 (\mcK^\bullet [\widetilde{\nabla}^\mr{ad}_\mcE])$ (cf. (\ref{G01})),
then the definition of $\widetilde{\phi}_{\mcE_B}$ implies that 
$\widetilde{\phi}_{\mcE_B} (\Phi^{-1}_W(p \Gamma (X, \Omega^{\otimes 2}_{X/W}))) \subseteq \Gamma  (X, \Omega^{\otimes 2}_{X/W})$.
We shall write
\begin{align}
\phi_0 : \Phi^*_{W} (\Gamma (X, \Omega^{\otimes 2}_{X/W})) \migi \Gamma (X, \Omega^{\otimes 2}_{X/W})
\end{align}
 for the composite of $\Phi_W^*([p]) : \Phi^*_{W} (\Gamma (X, \Omega^{\otimes 2}_{X/W})) \migi  \Phi^*_{W} (p\Gamma (X, \Omega^{\otimes 2}_{X/W}))$ and the restriction of $\widetilde{\phi}_{\mcE_B}$ to $ \Phi^*_{W} (p\Gamma (X, \Omega^{\otimes 2}_{X/W}))$.
By the above discussion,   the following   square diagram turns out to be commutative:
\begin{align} \label{W137}
\begin{CD}
\Phi^*_W(H^1 (X, \mcT_{X/W})) \oplus \Phi^*_W(\Gamma (X, \Omega_{X/W}^{\otimes 2}))  @> \Phi^*_W(\Upsilon_{(\mcE, \nabla_\mcE)})> \sim > \Phi^*_W(\mbH^1 (\mcK^\bullet [\nabla^\mr{ad}_\mcE]))
\\
@V  F^\flat_{(\mcE, \nabla_\mcE)} \oplus \phi_0
 VV @VV \phi_\mcE \circ (\Phi^*_W([p])) V
\\
H^1 (X, \mcT_{X/W}) \oplus \Gamma (X, \Omega_{X/W}^{\otimes 2})  @> \sim > \Upsilon_{(\mcE, \nabla_\mcE)}> \mbH^1 (\mcK^\bullet [\nabla^\mr{ad}_\mcE]).
\end{CD}
\end{align}
But,  
it follows from the definitions of $\phi_\mcE$ and $F_{(\mcE, \nabla_\mcE)}$  that $\phi_\mcE \circ (\Phi^*_W([p]))$ (i.e., the right-hand vertical arrow in the above  diagram) coincides with $F_{(\mcE, \nabla_\mcE)}$.
Thus, 
to obtain the required commutativities  of the diagrams,
  it suffices to verify the equality
$\phi_0 
= F^\sharp_{(\mcE, \nabla_\mcE)}$.

The commutative dagram (\ref{G731})
 induces a commutative square diagram of the form
\begin{align} \label{W135}
\xymatrix@!C=144pt{\Phi^*_{W}(\mbH^1 (\mcK^\bullet [\nabla^{\mr{ad}}_{\mcE}]))  
\ar[r]^{\Phi^*_{W}(\oint^\natural_{X, (\mcE, \nabla_\mcE)})}
\ar[d]_{F_{(\mcE, \nabla_\mcE)}}
&
 \Phi^*_{W}(\mbH^1 (\mcK^\bullet [\nabla^{\mr{ad}}_{\mcE}]))^\vee
 \ar[r]^{[p^3]}
 &
  \Phi^*_{W}(\mbH^1 (\mcK^\bullet [\nabla^{\mr{ad}}_{\mcE}]))^\vee
  \\
 \mbH^1 (\mcK^\bullet [\nabla^{\mr{ad}}_{\mcE}])
  \ar[rr]_{\oint^\natural_{X, (\mcE, \nabla_\mcE)}} &&
   \mbH^1 (\mcK^\bullet [\nabla^{\mr{ad}}_{\mcE}])^\vee \ar[u]_{F_{(\mcE, \nabla_\mcE)}^\vee},
}
\end{align}
where $[p^3]$ denotes multiplication by $p^3$.
The submodule $\Gamma (X, \Omega^{\otimes 2}_{X/W}) \subseteq \mbH^1 (\mcK^\bullet [\nabla^\mr{ad}_\mcE])$ is isotropic with respect to $\oint_{X, (\mcE, \nabla_\mcE)}$ (cf. Proposition \ref{W32}),
and the above  diagram restricts to a commutative square diagram of the form
\begin{align} \label{W134}
\xymatrix@!C=124pt{\Phi^*_{W}(\Gamma (X, \Omega_{X/W}^{\otimes 2}))  
\ar[r]^{\Phi^*_W(\int^\natural_X)}
\ar[d]_{\phi_0}
&
 \Phi^*_{W}(H^1 (X, \mcT_{X/W}))^\vee
 \ar[r]^{[p^3]}
 &
  \Phi^*_{W}(H^1 (X, \mcT_{X/W}))^\vee
  \\
 \Gamma (X, \Omega_{X/W}^{\otimes 2})
  \ar[rr]_{\int^\natural_X} &&
   H^1 (X, \mcT_{X/W})^\vee \ar[u]_{F_{(\mcE, \nabla_\mcE)}^{\flat \vee}}.
}
\end{align}
Since $(F^{\flat \vee}_{(\mcE, \nabla_\mcE)})^{-1} \circ [p^3] = F^\sharp_{(\mcE, \nabla_\mcE)}$ via $\oint_X^\natural$,
the above diagram
  implies 
the equality $\phi_0$ $= F^\sharp_{(\mcE, \nabla_\mcE)}$, as desired.
This completes the proof of the assertion.

\vspace{5mm}
\subsection{Proof of Theorem \ref{T013}} \label{SS699}
\leavevmode\\ \vspace{-4mm}


Next, we shall prove  Theorem \ref{T013}.
Denote by 
\begin{align}
0^*_\mcN(\widehat{\omega}^\mr{Liou}_{g, \mbZ_p} |_{\mcN}) 
\in \Gamma (\mcN, \bigwedge^2 0^*_{\mcN} (\Omega_{T^\vee_{\mbZ_p}\mcN/\mbZ_p})) 
\
\left(\text{resp.,} \ \sigma^*(\widehat{\omega}^\mr{PGL}_{g, \mbZ_p}|_\mcN) \in \Gamma (\mcN, \bigwedge^2 \sigma^*(\Omega_{\widehat{\mcS}_{g, \mbZ_p}/\mbZ_p})) \right)
\end{align}
 the 
restriction  of $\widehat{\omega}^\mr{Liou}_{g, \mbZ_p} |_\mcN$ (resp., $\widehat{\omega}^\mr{PGL}_{g, \mbZ_p}|_\mcN$) via the zero section $0_{\mcN} : \mcN \migi T^\vee_{\mbZ_p} \mcN$ (resp., $\sigma : \mcN \migi \widehat{\mcS}_{g, \mbZ_p}$).
Also, denote by
\begin{align}
\Lambda : \Gamma  (\mcN, \bigwedge^2 \sigma^*(\Omega_{\widehat{\mcS}_{g, \mbZ_p}/\mbZ_p})) \isom  \Gamma (\mcN, \bigwedge^2 0^*_{\mcN} (\Omega_{T^\vee_{\mbZ_p}\mcN/\mbZ_p})) 
\end{align}
the morphism obtained by  restricting  (\ref{G241}) to $0_{\mcN}$.
Then, let us consider the following lemma.

\vspace{3mm}
\ble \label{p401}\leavevmode\\
 \ \ \  
The following equality holds:
\begin{align} \label{Le020}
\Lambda (
\sigma^*(\widehat{\omega}^\mr{PGL}_{g, \mbZ_p}|_\mcN)
 ) = 
0^*_\mcN(\widehat{\omega}^\mr{Liou}_{g, \mbZ_p} |_{\mcN})
\end{align}
 \ele
\begin{proof}
Let us consider 
$F^\flat \oplus F^\sharp$ as a morphism 
$\Phi^*_\mcN(0^*_{\mcN}(\mcT_{T^\vee_{\mbZ_p} \mcN/\mbZ_p})) \migi 0^*_{\mcN}(\mcT_{T^\vee_{\mbZ_p} \mcN/\mbZ_p})$
via  (\ref{e741}).
This morphism induces a morphism
\begin{align}
F^\mr{Liou} :
\bigwedge^2 0^*_{\mcN}(\Omega_{T^\vee_{\mbZ_p} \mcN/\mbZ_p}) 
\migi \left(\bigwedge^2 \Phi^*_\mcN(0^*_{\mcN}(\Omega_{T^\vee_{\mbZ_p} \mcN/\mbZ_p}))  \cong \right) \ \Phi^*_\mcN(\bigwedge^2 0^*_{\mcN}(\Omega_{T^\vee_{\mbZ_p} \mcN/\mbZ_p})). 
\end{align}
Let us consider 
$\Gamma (\mcN, 0^*_{\mcN}(\Omega_{T^\vee_{\mbZ_p} \mcN/\mbZ_p}))$
as a submodule of $\Gamma (\mcN, \Phi^*_\mcN(\bigwedge^2 0^*_{\mcN}(\Omega_{T^\vee_{\mbZ_p} \mcN/\mbZ_p})))$ via pull-back by $\Phi_\mcN$.
By the definitions of 
$F^\flat$,  $F^\sharp$, and the explicit description of 
$0^*_\mcN(\widehat{\omega}^\mr{Liou}_{g, \mbZ_p} |_{\mcN})$,
 the following equality holds:
\begin{align} \label{G401}
F^\mr{Liou} (0^*_\mcN(\widehat{\omega}^\mr{Liou}_{g, \mbZ_p} |_{\mcN})) = p^3 \cdot  0^*_\mcN(\widehat{\omega}^\mr{Liou}_{g, \mbZ_p} |_{\mcN}).
\end{align}

Here, 
let us write $H := \Gamma (\mcN, \bigwedge^2 0^*_\mcN (\Omega_{T^\vee_{\mbZ_p}\mcN/\mbZ_p}))$, and suppose that 
$0^*_\mcN(\widehat{\omega}^\mr{Liou}_{g, \mbZ_p} |_{\mcN}) - \Lambda (\sigma^*(\widehat{\omega}^\mr{PGL}_{g, \mbZ_p} |_{\mcN}))$ $\in p^m H$ for some $m>0$, i.e., 
$0^*_\mcN(\widehat{\omega}^\mr{Liou}_{g, \mbZ_p} |_{\mcN}) - \Lambda (\sigma^*(\widehat{\omega}^\mr{PGL}_{g, \mbZ_p} |_{\mcN})) \in p^m \cdot h$
for some $h \in H$.
Then, 
\begin{align}
p^3 \cdot (0^*_\mcN(\widehat{\omega}^\mr{Liou}_{g, \mbZ_p}  |_{\mcN})
 - \Lambda (\sigma^*(\widehat{\omega}^\mr{PGL}_{g, \mbZ_p}  |_{\mcN}))) 
& = 
p^3 \cdot 0^*_\mcN(\widehat{\omega}^\mr{Liou}_{g, \mbZ_p}  |_{\mcN})
- \Lambda (p^3 \cdot \sigma^*(\widehat{\omega}^\mr{PGL}_{g, \mbZ_p}  |_{\mcN})) \\
& =  F^{\mr{Liou}} (0^*_\mcN(\widehat{\omega}^\mr{Liou}_{g, \mbZ_p}  |_{\mcN}))
- \Lambda (F^{\mr{PGL}} (\sigma^*(\widehat{\omega}^\mr{PGL}_{g, \mbZ_p}  |_{\mcN}))) \notag \\
& =  F^{\mr{Liou}} (0^*_\mcN(\widehat{\omega}^\mr{Liou}_{g, \mbZ_p}  |_{\mcN}))
-F^{\mr{Liou}} ( \Lambda (\sigma^*(\widehat{\omega}^\mr{PGL}_{g, \mbZ_p}  |_{\mcN}))) \notag \\
& =  F^{\mr{Liou}} (0^*_\mcN(\widehat{\omega}^\mr{Liou}_{g, \mbZ_p}  |_{\mcN})
-\Lambda (\sigma^*(\widehat{\omega}^\mr{PGL}_{g, \mbZ_p}  |_{\mcN}))) \notag \\
& =  F^{\mr{Liou}} ( p^m \cdot h) \notag \\
&  =  p^{pm} \cdot  F^{\mr{Liou}} (h). \notag
\end{align}
This implies  that
$0^*_\mcN(\widehat{\omega}^\mr{Liou}_{g, \mbZ_p}  |_{\mcN}) - \Lambda (\sigma^*(\widehat{\omega}^\mr{PGL}_{g, \mbZ_p}  |_{\mcN})) \in p^{pm-3} H \subseteq p^{m+1}H$ (where we recall the assumption  $p>3$).
By induction on $m$,
we see that 
\begin{align}
0^*_\mcN(\widehat{\omega}^\mr{Liou}_{g, \mbZ_p}  |_{\mcN}) - \Lambda (\sigma^*(\widehat{\omega}^\mr{PGL}_{g, \mbZ_p}  |_{\mcN}))  \in \bigcap_{m >0} p^m  H = \{ 0\},
\end{align}
 that is to say, 
 the equality $0^*_\mcN(\widehat{\omega}^\mr{Liou}_{g, \mbZ_p}  |_{\mcN}) = \Lambda (\sigma^*(\widehat{\omega}^\mr{PGL}_{g, \mbZ_p}  |_{\mcN}))$ holds.
Thus, in order to complete the proof of the assertion,  it suffices to varify the equality (\ref{Le020}) modulo $p$.

Denote by
$\langle -, - \rangle^{\mr{PGL}}$
  the bilinear map on 
  $\mbR^1 f_{\mcN_1*}(\mcT_{C_{\mcN, 1}/\mcN_1}) \oplus f_{\mcN _1*}(\Omega_{C_{\mcN,1}/\mcN_1}^{\otimes 2})$ (where $\mcN_1 := \mcN^\mr{ord}_{g, \mbF_p}$)
  corresponding to $\sigma^*(\widehat{\omega}^\mr{PGL}_{g, \mbZ_p}  |_{\mcN})$ mod $p$ via
  $\Upsilon$ (cf. (\ref{EE001}) and (\ref{G680})).
Let us fix  local sections $a, a' \in \mbR^1 f_{\mcN_1*}(\mcT_{C_{\mcN, 1}/\mcN_1})$ and   $b, b' \in f_{\mcN _1*}(\Omega_{C_{\mcN,1}/\mcN_1}^{\otimes 2})$.
The result of Proposition \ref{W32} implies that
\begin{align} \label{G712}
\langle a+ b, a'+b' \rangle^\mr{PGL} &= \langle a, a' \rangle^\mr{PGL} + \langle a, b' \rangle^\mr{PGL} + \langle b, a' \rangle^\mr{PGL} + \langle b, b' \rangle^\mr{PGL} \\
& =  \langle a, a' \rangle^\mr{PGL} + \langle a, b' \rangle - \langle a', b \rangle,
 \notag 
\end{align}
where $ \langle -, - \rangle$ denotes the natural pairing $\mbR^1 f_{\mcN_1*}(\mcT_{C_{\mcN, 1}/\mcN_1}) \times f_{\mcN _1*}(\Omega_{C_{\mcN,1}/\mcN_1}^{\otimes 2}) \migi \mcO_{\mcN_1}$.
Also,
according to Proposition \ref{p452}, 
any local section of $\mbR^1 f_{\mcN_1*}(\mcT_{C_{\mcN, 1}/\mcN_1})$
 (considered as a local section of $\mbR^1f_{\mcN_1*}(\mcK^\bullet [\nabla^\mr{ad}_{\mcE_{\mcN_1}}])$ via $\Upsilon$) may be represented, locally on $\mcN_1$,  by  a collection $(\{ a_{\alpha \beta} \}_{\alpha, \beta}, \{ b_\alpha \}_\alpha)$ as in (\ref{W105})   with $b_\alpha =0$ (for any $\alpha$).
In particular, $\mbR^1 f_{\mcN_1*}(\mcT_{C_{\mcN, 1}/\mcN_1})$ is isotropic with respect to $\langle -, - \rangle^\mr{PGL}$,  and hence, $\langle a, a' \rangle^\mr{PGL} =0$.
Thus, by (\ref{G712}), the equality 
\begin{align}
\langle a+ b, a'+b' \rangle^\mr{PGL} =  \langle a, b' \rangle - \langle a', b \rangle
\end{align}
holds.
It follows from the definition of $0^*_\mcN(\widehat{\omega}^\mr{Liou}_{g, \mbZ_p} |_\mcN)$ modulo $p$ (cf. (\ref{e75})) and the above equality that the equality (\ref{Le020}) modulo $p$ holds, as desired.
This completes the proof of the assertion.
\end{proof}
\vspace{3mm}



Now, let $n$ be a positive integer and write $R = \mbZ /p^n \mbZ$.
Let $S$ be an $R$-scheme admitting an  \'{e}tale  morphism $v : S \migi \mcN_n \ \left(:= \mcN \otimes R \right)$ which dominates any component of $\mcN_n$.
Denote by $\widehat{\omega}^\text{Liou}_{g, \mbZ_p} |_{S}$, $\widehat{\omega}^\text{PGL}_{g, \mbZ_p} |_{S}$,  and $\Theta_S$ the base-changes by the composite   $S \stackrel{v}{\migi} \mcN_n \migiincl \mcN$ of $\widehat{\omega}^\text{Liou}_{g, \mbZ_p} |_{\mcN}$, $\widehat{\omega}^\text{PGL}_{g, \mbZ_p} |_{\mcN}$,  and $\Theta$ respectively. 
Since the natural map 
\begin{align} 
\Gamma (T^\vee_R \mcN_{n}, \bigwedge^2\Omega_{T^\vee_R \mcN_{n}/R}) \migi \Gamma (T^\vee_R S, \bigwedge^2\Omega_{T^\vee_R S/R})\notag
\end{align}
induced by  $v$
 is injective, the proof of Theorem \ref{T013} may be  reduced to proving the equality (\ref{W210}) restricted to $S$, i.e., the equality 
\begin{align} \label{e68}
\Theta_S (\widehat{\omega}^\mr{PGL}_{g, \mbZ_p} |_S)  = \widehat{\omega}^\text{Liou}_{g, \mbZ_p} |_{S}
\end{align}
(for any $n$ and $S$).
Moreover, for the same reason,  we are always free to replace $S$ by any \'{e}tale covering of $S$. 

 Next, let us take
  $A \in \Gamma (S, \Omega_{S/R})$.
  We denote by  $\sigma_A : S \migi T^\vee_R S$ the section 
   corresponding to $A$.
  It follows from an argument similar to the argument in ~\cite{Wak3}, \S\,5.3 (and \S\,5.1),  that  the following  equality holds:
  \begin{align} \label{EE099}
  \sigma^*_S(\widehat{\omega}^\mr{PGL}_{g, \mbZ_p} |_S)  - 0_S^*(\widehat{\omega}^{\mr{Liou}}_{g, \mbZ_p}|_{S}) = \sigma_{A}^{*} (\Theta_S^*(\widehat{\omega}^{\mr{PGL}}_{g, \mbZ_p}|_S) - \widehat{\omega}^{\mr{Liou}}_{g, \mbZ_p}|_{S}),
  \end{align}
where $\sigma_S$ denotes the morphism $S \migi \mcN \times_{\widehat{\mcM}_g} S$ induced by $\sigma$.
After possibly replacing $S$ by 
its \'{e}tale covering,
 we may assume  that {\it $S$ is affine and  the vector bundle $\Omega_{S/R}$ is free}.
Under this assumption, 
$\Theta_S^*(\widehat{\omega}^{\mr{PGL}}_{g, \mbZ_p}|_S) - \widehat{\omega}^{\mr{Liou}}_{g, \mbZ_p}|_{S}=0$ if and only if
$ \sigma^*_S(\widehat{\omega}^\mr{PGL}_{g, \mbZ_p} |_S)  - 0_S^*(\widehat{\omega}^{\mr{Liou}}_{g, \mbZ_p}|_{S}) =0$ for all
$A \in \Gamma (S, \Omega_{S/R})$.
Thus, in order to
verify the equality (\ref{e68}),
   it suffices 
   (by (\ref{EE099}))
    to prove  the equality
   $\sigma^*_S(\widehat{\omega}^\mr{PGL}_{g, \mbZ_p} |_S)  = 0_S^*(\widehat{\omega}^{\mr{Liou}}_{g, \mbZ_p}|_{S})$.
In particular, it suffices to prove the equality
\begin{align}
\sigma^*(\widehat{\omega}^\mr{PGL}_{g, \mbZ_p} |_{\mcN}) = \widehat{0}_\mcN^*(\widehat{\omega}^\mr{Liou}_{g, \mbZ_p}|_\mcN).
\end{align}
But, this equality holds by Lemma  \ref{p401}.
This completes the proof of Theorem \ref{T013}.
\vspace{10mm}
\section{Appendix (Crystals of  torsors and connections)} \label{A2}\vspace{3mm}


In this Appendix, we study  crystals of torsors (equipped with a structure group)
  and prove the bijective correspondence between crystals of torsors and quasi-nilpotent  flat torsors  (cf. Theorem \ref{T043}).
This correspondence  enable us to understand  the relationship (cf. Proposition \ref{W22})  between  the respective deformations  of a prescribed flat torsor over 
distinct underlying spaces.
Notice that its  application 
 to the case
 of indigenous bundles
(cf. Proposition \ref{W23}) was used in the proof of the main theorem.

\vspace{5mm}
\subsection{Connections on torsors} \label{SS6791}
\leavevmode\\ \vspace{-4mm}

Let $R$ be a commutative ring with unit, 
$G$ a geometrically connected  smooth algebraic group over $R$ with Lie algebra $\mfg$,
$S$  a scheme over $R$, and  $f : X \migi S$ a  smooth scheme over $S$ of relative dimension $n>0$.
Suppose that we are given a $G$-torsor 
$\pi : \mcE \migi X$
   over $X$.
 Denote by 
 \begin{align} \label{W230}
 \mr{Ad} (\mcE)  \ \left(:= \mcE \times^{G} \mfg\right)
 \end{align}
  the adjoint vector bundle on $X$ associated to $\mcE$ (i.e., the vector bundle obtained from $\mcE$ by the change of structure group via the adjoint representation $G \migi \mr{GL}(\mfg)$).
Also, denote by
 \begin{align} \label{W31}
 \widetilde{\mcT}_{\mcE/S} \ \left( := (\pi_*(\mcT_{\mcE/S}))^G \right)
 \end{align}
  the subsheaf of $\pi_*(\mcT_{\mcE/S})$ consisting of $G$-invariant sections.
Then, the differential of $\pi$ induces a surjection $\widetilde{\mcT}_{\mcE/S} \migisurj \mcT_{X/S}$, which we denote by $d \pi$.
The kernel of $d \pi$ may be naturally  identified with
 $\mr{Ad}(\mcE) \ \left(\cong (\pi_*(\mcT_{\mcE/X}))^G \right)$.
Thus, we have a short exact sequence
\begin{align} \label{H001}
0 \longmigi \mr{Ad}(\mcE) \longmigi \widetilde{\mcT}_{\mcE/S} \stackrel{d \pi}{\longmigi} \mcT_{X/S} \longmigi 0.
\end{align}
An {\bf $S$-connection} on $\mcE$ is, by definition, a split injection of (\ref{H001}), i.e., an $\mcO_X$-linear morphism $\nabla_\mcE : \mcT_{X/S} \migi \widetilde{\mcT}_{\mcE/S}$ satisfying the equality  
$d \pi \circ \nabla_\mcE = \mr{id}_{\mcT_{X/S}}$.
If $G=\mr{GL}_n$ for some $n>0$,
then the above definition of an $S$-connection  is   equivalent to the classical
definition of an $S$-connection on
  the corresponding vector bundle $\mcV := \mcE \times^{\mr{GL}_n} R^{\oplus n}$ (cf. ~\cite{Wak5}, \S\,4.2), i.e., an $f^{-1} (\mcO_S)$-linear morphism $\mcV \migi \Omega_{X/S} \otimes \mcV$ satisfying the Leibniz rule.
  In this case, we shall not distinguish between these definitions of an $S$-connection.
Denote by 
\begin{align} \label{W30}
\nabla^\mr{ad}_\mcE : \mr{Ad}(\mcE) \migi \Omega_{X/S} \otimes \mr{Ad}(\mcE)
\end{align}
 the $S$-connection  on the vector bundle  $\mr{Ad}(\mcE)$ induced by $\nabla_\mcE$ by the  change of structure group via the adjoint representation $G \migi \mr{GL}(\mfg)$.


Let us fix an $S$-connection $\nabla_\mcE$ on $\mcE$.
The {\bf curvature} of $\nabla_\mcE$ is, by definition,  the $\mcO_X$-linear morphism
$\bigwedge^2 \mcT_{X/S} \migi \mr{Ad}(\mcE) \ \left(\subseteq \widetilde{\mcT}_{\mcE/S} \right)$ determined by assigning $\partial_1 \wedge \partial_2 \mapsto [\nabla_\mcE (\partial_1), \nabla (\partial_2)] - \nabla_\mcE ([\partial_1, \partial_2])$ (for any local sections  $\partial_1, \partial_2 \in \mcT_{X/S}$).
We shall say that $\nabla_\mcE$ is {\bf flat} if its curvature  vanishes identically on $X$.
(If $X/S$ is of relative dimension $1$, which implies that $\mcT_{X/S}$ is a line bundle, then any $S$-connection on $\mcE$   is automatically flat.)
It is verified that
$\nabla_\mcE$ is flat if and only if $\nabla^\mr{ad}_\mcE$ is flat.
By a  {\bf flat $G$-torsor} over $X/S$, we mean 
 a pair $(\mcE, \nabla_\mcE)$ consisting of a $G$-torsor over $X$ and a flat  $S$-connection $\nabla_\mcE$ on $\mcE$.
An {\bf isomorphism} of flat $G$-torsors  from $(\mcE, \nabla_\mcE)$ to $(\mcE', \nabla_{\mcE'})$
is an isomorphism $\mcE \isom \mcE'$ of $G$-torsors compatible with the respective connections $\nabla_\mcE$ and $\nabla_{\mcE'}$.
Thus, flat $G$-torsors over $X/S$ and isomorphisms between them forms a groupoid.

\vspace{5mm}
\subsection{$p$-curvature of flat torsors} \label{SS7791}
\leavevmode\\ \vspace{-4mm}

In this subsection, suppose further  that
 $p \cdot \mcO_S = 0$.
Write $\Phi_S : S \migi S$ for the absolute Frobenius morphism of $S$, $f^{(1)} : X^{(1)} \ \left(:= X \times_{\Phi_S, S}S \right)\migi S$ for the Frobenius twist of $X$ relative to $S$, and 
$\Phi_{X/S}:X \migi X^{(1)}$ for the relative Frobenius morphism.
 Recall that
for each $\mcO_{X^{(1)}}$-module $\mcF$, the pull-back $\Phi^*_{X/S}(\mcF)$ 
 admits a canonical $S$-connection 
\begin{align} \label{W110}
\nabla^\mr{can}_\mcF : \Phi^*_{X/S}(\mcF) \migi \Omega_{X/S} \otimes \Phi^*_{X/S}(\mcF)
\end{align}
 determined uniquely by the condition that
the sections of $\Phi^{-1}_{X^{(1)}/S} (\mcF)$ are horizontal.

 Now, let us fix a flat $G$-torsor $(\mcE, \nabla_\mcE)$ over $X/S$.
  The {\bf $p$-curvature} of 
    $(\mcE, \nabla_\mcE)$ 
     is defined to be the $\mcO_X$-linear morphism 
\begin{align} \label{Dd010}
\psi_{(\mcE, \nabla_\mcE)} : \Phi_{X/S}^*(\mcT_{X^{(1)}/S}) \migi \mr{Ad}(\mcE) \ \left(\subseteq \widetilde{\mcT}_{\mcE/S} \right)
\end{align}
 determined uniquely  by $F^{-1}_{X/S}(\partial) \mapsto  \nabla_\mcE (\partial)^{[p]} - \nabla_\mcE (\partial^{[p]})$ for any local section $\partial \in \mcT_{X/S}$ (cf. ~\cite{Kal}, \S\,5.0).
Here, $(-)^{[p]}$
denotes the operator on  vector fields  given by taking the $p$-th iterates of the corresponding derivations, by which both $\mcT_{X^{(1)}/S}$ and $\widetilde{\mcT}_{\mcE/S}$ form  sheaves  of $p$-Lie algebras.
As  is well-known,
  $\psi_{(\mcE, \nabla_\mcE)}$
 is compatible with the respective connections $\nabla^\mr{can}_{\mcT_{X^{(1)}/S}}$ and $\nabla^\mr{ad}_\mcE$.
In particular, 
if 
\begin{align} \label{W118}
\psi_{(\mcE, \nabla_{\mcE})}^\nabla : \Phi^{-1}_{X/S}(\mcT_{X^{(1)}/S}) \migi  
\mr{Ad}(\mcE)
\end{align}
denotes the restriction of $\psi_{(\mcE, \nabla_{\mcE})}$,
then its image lies in $\mr{Ker} (\nabla^\mr{ad}_\mcE)$.

Finally, we shall say that $(\mcE, \nabla_\mcE)$ is {\bf $p$-nilpotent} if $\psi_{(\mcE, \nabla_\mcE)}$ has nilpotent image.
It is verified that
$(\mcE, \nabla_\mcE)$ is $p$-nilpoent  if and only if $(\mr{Ad}(\mcE), \nabla^\mr{ad}_\mcE)$ is $p$-nilpotent.

\vspace{5mm}
\subsection{Quasi-nilpotence and crystals} \label{Ss1}
\leavevmode\\ \vspace{-4mm}

In  this subsection,  we shall suppose  that 
$p^N \cdot \mcO_S=0$ for some $N>0$.
 Fix a flat $G$-torsor $(\mcE, \nabla_\mcE)$ over $X/S$.
 Let $U_+ := (U, \{x_i \}_{i=1}^n)$ be a collection, where $U$ denotes  an open subscheme of $X$ and $\{ x_i \}_{i=1}^n \ (\subseteq \Gamma (U, \mcO_{X}))$ denotes a local coordinate system defined on $U$ relative to $S$;
 we shall refer to such a collection as a {\bf coordinate chart} for $X/S$.
 Given 
 a coordinate chart 
 $U_+ = (U, \{x_i \}_{i=1}^n)$ for $X/S$,
 we shall consider the following condition: 
\begin{itemize}
\item[$(*)_{U_+}$ :]
 For each $s \in \Gamma (U, \mr{Ad}(\mcE))$, there exist an open covering $\{ U_\alpha \}_\alpha$ of $U$ and a set  of positive integers $\{ e_{i, \alpha} \}_{i, \alpha}$
 such that $\mr{ad}(\nabla_\mcE (\frac{\partial}{\partial x_i}))^{e_{i, \alpha}} (s |_{U_\alpha})=0$ for all $i$ and $\alpha$, where $\mr{ad}(v)$ (for each $v \in \widetilde{\mcT}_{\mcE/S}$) denotes the adjoint operator $[v, -] : \widetilde{\mcT}_{\mcE/S} \migi \widetilde{\mcT}_{\mcE/S}$.
\end{itemize}


\vspace{3mm}
\bde \label{DD011} \leavevmode\\
 \ \ \ 
We shall say that $(\mcE, \nabla_\mcE)$ is {\bf quasi-nilpotent}
 if
there exists a collection $\mfU_+ := \{ U_{\gamma, +} \}_\gamma$,
 where each 
 $U_{\gamma, +}$ denotes a coordinate chart  $(U_{\gamma}, \{ x_{\gamma, i}\}_{i=1}^n)$ for $X/S$
  satisfying the condition $(*)_{U_{\gamma, +}}$ such that $\{ U_\gamma \}_{\gamma}$ forms an open covering of $X$.
 \ede

\vspace{3mm}
\begin{rema}  \label{W40}\leavevmode\\
\vspace{-5mm}
\begin{itemize}
\item[(i)]
Suppose that  $\nabla_\mcE$ is flat.
Then, it is  verified that 
{\it any coordinate chart  $U_{+} := (U, \{ x_i \}_{i=1}^n)$ for $X/S$
 satisfies the condition $(*)_{U_+}$} (cf. ~\cite{BO}, Remark 4.11).
\item[(ii)]
It is verifies that {\it $\nabla_\mcE$ is quasi-nilpotent if and only if its reduction modulo $p$ is quasi-nilpotent.}
Also, the quasi-nilpotence of a connection may be related to the $p$-nilpotence of its reduction modulo $p$.
In fact,
let $X_1/S_1$ and $(\mcE_1, \nabla_{\mcE, 1})$ be the reductions modulo $p$ of $X/S$ and $(\mcE, \nabla_\mcE)$ respectively.
 Let us fix a coordinate chart  $U_+ := (U, \{ x_i \}_{i=1}^n)$
  for $X_1/S_1$.
Since $\left( \frac{\partial}{\partial x_i}\right)^p =0$,
the following sequence of equalities holds:
\begin{align}
\mr{ad}(\psi_{(\mcE_1, \nabla_{\mcE,1})}(F^{-1}_{X_1/S_1} \left(\frac{\partial}{\partial x_i}\right)))
 &  = \mr{ad}(\nabla_{\mcE, 1} (\left(\frac{\partial}{\partial x_i}\right))^{[p]}- \nabla_{\mcE, 1} (\left(\frac{\partial}{\partial x_i}\right)^{[p]}) )\\
  &  = \mr{ad}(\nabla_{\mcE,1} (\left(\frac{\partial}{\partial x_i}\right))^{[p]}) \notag \\
  &  = \mr{ad}(\nabla_{\mcE,1} (\frac{\partial}{\partial x_i}))^p. \notag
\end{align}
Hence, {\it $(\mcE, \nabla_\mcE)$ is quasi-nilpotent (or equivalently, $(\mcE_1, \nabla_{\mcE, 1})$ is quasi-nilpotent) if and only if  ($\mr{ad}(\nabla_{\mcE,1} (\frac{\partial}{\partial x_i}))$ is nilpotent for any $i$, or equivalently) the reduction  of $(\mcE, \nabla_\mcE)$ modulo $p$  is $p$-nilpotent}.
\end{itemize}
 \end{rema}
\vspace{3mm}

Let  $(S, I, \gamma)$ be  a PD scheme over $R$
(with $I$ a quasi-coherent ideal) and 
$X$ an $S$-scheme to which $\gamma$ extends.
Denote by 
\begin{align}
\mr{Crys} (X/S)
\end{align}
 the crystalline site, which is  the site whose objects are 
 divided power thickenings, i.e.,
pairs $(U \migiincl T, \delta)$, where $U$ is a Zariski open subscheme of $X$, $U \migiincl T$ is a closed $S$-immersion defined by an ideal $J$, and $\delta$ is a PD structure on $J$ which is compatible with $\gamma$ in the evident sense.
We shall often abuse notation by writing $(U, T, \delta)$ for $(U \migiincl T, \delta)$, or even by just writing $T$ for the whole thing.
(We shall call $(U \migiincl T, \delta)$ an $S$-PD thickening of $U$.)
The morphisms and the covering families in $\mr{Crys}(X/S)$ are defined in the usual manner.
Recall that  a {\bf crystal of  $G$-torsors} over   $\mr{Crys}(X/S)$ (resp., a {\bf crystal} of vector bundles on $\mr{Crys}(X/S)$)  is a 
cartesian section  of the fibered category of $G$-torsors over $\mr{Crys}(X/S)$ (resp., the fibered category of vector bundles on $\mr{Crys}(X/S)$).
For each crystal $\mcF^\lozenge$ (of either $G$-torsors or vector bundles) on $\mr{Crys}(X/S)$ and each  divided power thickening  $T$ in $\mr{Crys}(X/S)$,
 we shall write $\mcF_{T}^\lozenge$ for 
  the evaluation of $\mcF^\lozenge$ on  this thickening.

\vspace{5mm}
\subsection{Correspondence between crystals and quasi-nilpotent flat torsors} \label{TSS179}
\leavevmode\\ \vspace{-4mm}



Let us suppose that  the following condition $(**)_G$ on $G$ is satisfied:
\begin{itemize}
\item[$(**)_G$ :]
$G$ is a simple algebraic group over $R$ of adjoint type satisfying the inequality $p >h$,
where $h$ denotes the Coxeter number of $G$.
\end{itemize}
(For instance, $G  = \mr{PGL}_n$ with $n <p$.)
In particular, 
the morphism of algebraic groups $G \migi \mr{Aut}^0 (\mfg)$
  obtained as the adjoint representation of $G$ is an isomorphism, where $\mr{Aut}^0 (\mfg) \ \left( \subseteq \mr{GL}(\mfg)\right)$ denotes  the identity component of the group of Lie algebra automorphisms of $\mfg$.
(Indeed, it follows from a result in ~\cite{On} that the  reduction modulo $p$ of this morphism is an isomorphism.)


Let $X$ be a smooth scheme over $S$ and denote by  $D_{X/S}(1)$
the divided power envelope of $X$ in $X \times_S X$ via the diagonal embedding $\Delta :X \migi X \times_S X$.
For $i =1, 2$, we shall write $\mr{pr}_i : D_{X/S} (1) \migi X$ for the composite of the natural morphism $D_{X/S} (1) \migi X \times_S X$ and the projection $X \times_S X \migi X$ onto the $i$-th factor.

\vspace{3mm}
\bde \label{FD1} \leavevmode\\
 \ \ \
Let $\mcE$ be a $G$-torsor over $X$.
 An {\bf HPD  stratification} on $\mcE$ is an isomorphism
$e : \mr{pr}^*_2 (\mcE) \isom \mr{pr}^*_1 (\mcE)$ of $G$-torsors over  $D_{X/S}(1)$  whose restriction to $X$ via $\Delta$
coincides with the identity morphism of $\mcE$ and which 
 satisfies the cocycle condition  in the usual  sense (cf. ~\cite{BO},  \S\,2, the comment following Definition 2.10).
  \ede


\bt \label{T043} \leavevmode\\
 \ \ \ 
Let $J$ be a sub-PD ideal of $I$.
Write $\overline{S}$ the closed subscheme of $S$ defined by $J$  and write $\overline{X} := X \times_S \overline{S}$.
 Then,  the following categories are naturally equivalent:
 \begin{itemize}
 \item[(i)]
 The category of crystals of  $G$-torsors  on $\mr{Crys}(\overline{X}/S)$;
 \item[(ii)]
 The category of $G$-torsors $\mcE$ over  $X$ together  with an HPD stratification $\mr{pr}_2^*(\mcE) \isom \mr{pr}_1^*(\mcE)$ on it;
 \item[(iii)]
 The category of $G$-torsors $\mcE$ over $X$ together with an HPD stratification $\mr{pr}_2^*(\mr{Ad}(\mcE)) \isom \mr{pr}_1^*(\mr{Ad}(\mcE))$ on its adjoint bundle $\mr{Ad}(\mcE)$ (in the sense of ~\cite{BO}, \S\,4, Definition 4.3H) that is compatible with the respective Lie bracket structures pulled-back from $\mr{Ad}(\mcE)$;
 \item[(iv)]
 The category of $G$-torsors  $\mcE$  over $X$ together with a quasi-nilpotent flat $S$-connection on $\mr{Ad}(\mcE)$ compatible with the Lie bracket structure.
 \item[(v)]
 The category of quasi-nilpotent  flat $G$-torsors over $X/S$.
 \end{itemize}
  \et
\begin{proof}
The equivalence of (ii) and (iii), as well as (iv) and (v),  follows from the isomorphism $G \isom \mr{Aut}^0(\mfg)$.
The equivalent of (iii) and (iv) follows from ~\cite{BO}, \S\,4, Theorem 4.12 (and its proof).
Thus, it suffices to
check the equivalence of (i) and (ii).
First, let  $\mcE^\lozenge$ be a crystal of $G$-torsors on $\mr{Crys}(\overline{X}/S)$.
The morphism $\mr{pr}_i : D_{X/S}(1) \migi X$ (for each $i=1,2$) induces 
an isomorphism $e_i : \mr{pr}_i^*(\mcE^\lozenge_X) \isom \mcE^\lozenge_{D_{X/S}(1)}$.
Then, the composite  isomorphism  $e_{\mcE^\lozenge} := e_1^{-1} \circ e_2 : \mr{pr}_2^*(\mcE^\lozenge_X)\isom \mr{pr}_2^*(\mcE^\lozenge_X)$  specifies an HPD stratification on the $G$-torsor $\mcE^\lozenge_X$.

Conversely, given a $G$-torsor $\mcE$ over $X$ together with an  HPD stratification $e : \mr{pr}_2^*(\mcE) \isom \mr{pr}_1^*(\mcE)$, we construct a crystal $\mcE^\lozenge$ of  $G$-torsors on $\mr{Crys} (\overline{X}/S)$.
To this end, 
it suffices to specify $\mcE_T^\lozenge$ for sufficiently small $(U, T, \delta) \in \mr{Ob}(\mr{Crys}(\overline{X}/S))$,
e.g., 
so that
 there exists an $S$-morphism $h : T \migi X$ over $(S, I, \gamma)$ extending  the open immersion  $U \migiincl \overline{X}$.
Then, let us define $\mcE_T^\lozenge$ to be the pull-back $h^*(\mcE)$.
For a morphism $h' : T \migi X$   as $h$, 
we obtain $(h, h') : T \migi D_{X/S}(1)$.
The pull-back of $\epsilon$ via this morphism determines an isomorphism $(h, h')^*(\epsilon): h^*(\mcE) \isom h'^*(\mcE)$.
The fact that $\mcE^\lozenge_T$ does not depend on $h$ (up to canonical isomorphism) comes from 
the isomorphisms $(h, h')^*(\epsilon)$ (for  $(h, h')$'s).
Thus, $\mcE_T^\lozenge$'s for various $(T, T, \delta)$'s  forms a crystal $\mcE_e^\lozenge$ of $G$-torsors.
One verifies immediately that the assignments $\mcE^\lozenge \mapsto e_{\mcE^\lozenge}$, $e \mapsto \mcE_e^\lozenge$ define the equivalence of (i) and (ii).
This completes the proof of the assertion.
 \end{proof}
\vspace{3mm}

For each quasi-nilpotent flat $G$-torsors $(\mcE, \nabla_\mcE)$ over $X/S$,
we shall write 
\begin{align}
(\mcE, \nabla_\mcE)^\lozenge
\end{align} 
for the crystal of $G$-torsors corresponding to $(\mcE, \nabla_\mcE)$ via the equivalence of categories between (i) and (v) in Theorem \ref{T043}.


\vspace{3mm}
\begin{rema}  \label{W41}\leavevmode\\
 \ \ \ 
Let us consider the case where  $G = \mr{PGL}_n$ with $n <p$.
 Then,  the equivalences of categories obtained in Theorem \ref{T043} are  compatible, via projectivization,  with those obtained in   the corresponding classical result  for crystals of rank $n$ vector bundles (i.e., crystals of $\mr{GL}_n$-torsors) described, e.g., ~\cite{BO}, \S\,6, Theorem 6.6.
To be precise, let $(\mcV, \nabla_\mcV)$ be a quasi-nilpotent  flat vector bundle of rank $n$ and denote by
$\mcV^\lozenge$ the corresponding crystal obtained by the result in {\it loc.\,cit.}.
Then, the flat $\mr{PGL}_n$-torsor obtained from $(\mcV, \nabla_\mcV)$ via projectivization (i.e., via the change of structure group by the quotient $\mr{GL}_n \migisurj \mr{PGL}_n$)
corresponds, via the equivalence of categories in Theorem \ref{T043}, to the crystal given by assigning, to each $(U \migiincl T, \gamma) \in \mr{Ob} (\mr{Crys} (X/S))$, the projectivization  of $\mcV^\lozenge_T$.
 \end{rema}
\vspace{3mm}

Let us describe the following two 
corollaries of Theorem \ref{T043}.
We shall   keep the notation in
that theorem.

\vspace{3mm}
\bco \label{W20} \leavevmode\\
 \ \ \
 There exists an equivalence of categories between the category of crystals of  $G$-torsors on $\mr{Crys}(\overline{X}/S)$ and the category of crystals of $G$-torsors on $\mr{Crys}(X/S)$.
  \eco

\vspace{3mm}
\bco \label{W21} \leavevmode\\
 \ \ \
Suppose that we are given another smooth scheme $X'$ over $S$ whose reduction modulo $J$ is isomorphic to $\overline{X}$.
Then,
for each quasi-nilpotent flat $G$-torsor $(\mcE, \nabla_\mcE)$ over $X/S$,
there exists a unique (up to isomorphism) quasi-nilpotent flat $G$-torsor 
\begin{align} \label{Le001}
(\lambda_{X'} (\mcE), \lambda_{X'}(\nabla_{\mcE}))
\end{align}
 over $X'/S$ such that the crystals of $G$-torsors over  $\mr{Crys}(\overline{X}/S)$ corresponding, via the equivalence of categories in Theorem \ref{T043}, to $(\mcE, \nabla_\mcE)$ and  $(\lambda_{X'} (\mcE), \lambda_{X'}(\nabla_{\mcE}))$ are isomorphic.
Moreover, the assignment $(\mcE, \nabla_\mcE) \mapsto (\lambda_{X'} (\mcE), \lambda_{X'}(\nabla_{\mcE}))$ determines
 an equivalence of categories between the category of quasi-nilpotent flat $G$-torsors on $X/S$ and the category of quasi-nilpotent flat $G$-torsors on $X'/S$.
  \eco

\vspace{5mm}
\subsection{Deformation space of flat torsors} \label{SS983}
\leavevmode\\ \vspace{-4mm}

 In this subsection, we describe the change of flat torsors $(\mcE, \nabla_\mcE) \mapsto (\lambda_{X'} (\mcE), \lambda_{X'}(\nabla_{\mcE}))$ obtained in Corollary \ref{W21}
   in terms of 
    de Rham cohomology of complexes.
For each morphism of sheaves $\nabla : \mcK^0 \migi \mcK^1$ on $X$, 
we shall write $\mcK^\bullet [\nabla]$ for $\nabla$   regarded as a complex concentrated at degree $0$ and $1$.
Also, for each sheaf $\mcF$ on $X$ and each $n \in \mbZ$, we shall write $\mcF [n]$ for $\mcF$ considered as a complex concentrated at degree $n$.

  Let us keep the  notation in the previous subsection
 and suppose that $S = \mr{Spec}(R)$ and $X$ is a curve over $S$ (cf. \S\,\ref{SS01W}).
Denote by 
\begin{align}
\widetilde{\nabla}^\mr{ad}_\mcE : \widetilde{\mcT}_{\mcE/R} \migi \Omega_{X/R} \otimes \mr{Ad} (\mcE) 
\end{align}
the  unique  
$R$-linear 
morphism determined  by the condition that
\begin{align}
\langle \partial_1, \widetilde{\nabla}_\mcE^\mr{ad} (\partial_2)  \rangle = [\nabla_\mcE (\partial_1), \partial_2] - \nabla_\mcE ([\partial_1, d \pi (\partial_2)]) 
\end{align}
(cf. (\ref{H001}) for the definition of $d\pi$) for any local sections $\partial_1 \in \mcT_{X/R}$ and $\partial_2 \in \widetilde{\mcT}_{\mcE/R}$,
 where
   $\langle -, - \rangle$ denotes the pairing $\mcT_{X/R} \times (\Omega_{X/R} \otimes \mr{Ad} (\mcE)) \migi \mr{Ad} (\mcE)$ arising from the natural pairing  $\mcT_{X/R} \times \Omega_{X/R} \migi \mcO_{X}$.
The short exact sequence (\ref{H001}) induces naturally  the following  short exact sequence of complexes:
\begin{align}
0 \longmigi \mcK^\bullet [\nabla^\mr{ad}_\mcE] \longmigi \mcK^\bullet [\widetilde{\nabla}^\mr{ad}_\mcE] \longmigi \mcT_{X/R}[0] \longmigi 0. 
\end{align}
This sequence  gives rise to 
a short exact sequence
\begin{align} \label{Ee1000}
0 \longmigi \mbH^1 (\mcK^\bullet [\nabla^\mr{ad}_\mcE]) \stackrel{\iota^\mbH}{\longmigi} \mbH^1 (\mcK^\bullet [\widetilde{\nabla}^\mr{ad}_\mcE]) \stackrel{\pi^\mbH}{\longmigi} H^1 (X, \mcT_{X/R}) \longmigi 0.
\end{align}
By passing to the   injection $\iota^\mbH$, we shall consider $\mbH^1 (\mcK^\bullet [\nabla^\mr{ad}_\mcE])$ as a submodule of $\mbH^1 (\mcK^\bullet [\widetilde{\nabla}^\mr{ad}_\mcE])$.

Here,  we shall  write $R_\epsilon := R [\epsilon]/(\epsilon^2)$, in which the ideal $I R_\epsilon$ is endowed with a divided power structure extended  from $I \subseteq R$.
  In what follows, we shall denote the base-changes  to $R_\epsilon$ of objects over $R$  by means of a subscripted $\epsilon$.
It is well-known that there exists a canonical bijection
\begin{align} \label{W102}
H^1 (X, \mcT_{X/R}) \isom \mr{Def}_{R_\epsilon} (X)
\end{align}
 between the set $H^1 (X, \mcT_{X/R})$ and the set $\mr{Def}_{R_\epsilon} (X)$ consisting of  isomorphism classes of  deformations of $X$ over $R_\epsilon$.
Also, denote by
\begin{align}
\mr{Def}_{R_\epsilon} (X, \mcE, \nabla_\mcE)
\end{align}
the set of isomorphism classes of deformations over $R_\epsilon$  of $(X, \mcE, \nabla_\mcE)$ (as a data consisting of  a curve  and a flat $G$-torsor over it).
By well-known  generalities on the deformation theory of connections, 
there exists a canonical bijections
\begin{align} \label{W100}
 \mbH^1 (\mcK^\bullet [\widetilde{\nabla}^\mr{ad}_\mcE]) \isom \mr{Def}_{R_\epsilon} (X, \mcE, \nabla_\mcE)
\end{align}
 making the following square diagram commute:
\begin{align}
\begin{CD}
\mbH^1 (\mcK^\bullet [\widetilde{\nabla}^\mr{ad}_\mcE]) @> (\ref{W100})> \sim >
\mr{Def}_{R_\epsilon} (X, \mcE, \nabla_\mcE)
\\
@V \pi^\mbH VV @VVV
\\
H^1 (X, \mcT_{X/R}) @> \sim > (\ref{W102}) > \mr{Def}_\epsilon (X),
\end{CD}
\end{align}
where the right-hand vertical arrow denotes the projection induced by forgetting the data of a deformation of $(\mcE, \nabla_\mcE)$.
Moreover, denote by
\begin{align}
\mr{Def}_{R_\epsilon} (\mcE, \nabla_\mcE)
\end{align}
the set of isomorphism classes of deformations over $X_\epsilon$ of $(\mcE, \nabla_\mcE)$ (as a flat $G$-torsor), which may be though of as a subset of $\mr{Def}_{R_\epsilon} (X, \mcE, \nabla_\mcE)$, i.e., the subset consisting of deformations whose underlying curves are the trivial deformation of $X$.
One verifies immediately that the bijection  (\ref{W100}) restricts  to a bijection
\begin{align} \label{W109}
\mbH^1 (\mcK^\bullet [\nabla^\mr{ad}_\mcE]) \isom \mr{Def}_{R_\epsilon} (\mcE, \nabla_\mcE).
\end{align}

\vspace{3mm}
\begin{rema}  \label{W101}\leavevmode\\
 \ \ \ 
 In this remark, we describe the bijection (\ref{W100}) in terms of \v{C}ech cohomology.
Let us take an affine open covering $\mfU := \{ U_\alpha \}_{\alpha \in I}$ (where $I$ is an index set) of $X$.
We shall write $I_2$ for the set of pairs $(\alpha, \beta) \in I \times I$ with $U_{\alpha \beta} := U_\alpha \cap U_\beta \neq \emptyset$.
One may calculate $\mbH^1 (\mcK^\bullet[\widetilde{\nabla}^\mr{ad}_\mcE])$ as the total cohomology   of the \v{C}ech double complex   $\mr{Tot}^\bullet (\check{C}^\bullet (\mfU, \mcK^\bullet [\widetilde{\nabla}^\mr{ad}_\mcE]))$
 associated to $\mcK^\bullet [\widetilde{\nabla}^\mr{ad}_\mcE]$.
Each element $v$ of $\mbH^1 (\mcK^\bullet[\widetilde{\nabla}^\mr{ad}_\mcE])$  may be given by a collection of data
\begin{align}\label{W105}
 v = (\{ a_{\alpha \beta} \}_{\alpha, \beta}, \{ b_\alpha \}_{\alpha})
\end{align}
consisting of a \v{C}ech $1$-cocycle 
$\{ a_{\alpha\beta} \}_{\alpha, \beta} \in \check{C}^1 (\mfU, \widetilde{\mcT}_{\mcE/R})$
 (where $a_{\alpha\beta} \in \Gamma (U_{\alpha \beta}, \widetilde{\mcT}_{\mcE/R})$) and a \v{C}ech $0$-cochain $\{b_\alpha \}_\alpha \in \check{C}^0 (\mfU, \Omega_{X/R} \otimes \mr{Ad}(\mcE))$
(where $b_\alpha \in \Gamma (U_\alpha, \Omega_{X/R} \otimes \mr{Ad}(\mcE))$ $= \mr{Hom}_{\mcO_{U_\alpha}}(\mcT_{U_\alpha/R}, \mr{Ad}(\mcE)|_{U_\alpha})$)
which agree under $\widetilde{\nabla}^\mr{ad}_\mcE$ and the \v{C}ech coboundary map.
(The elements in $\mbH^1 (\mcK^\bullet [\nabla^\mr{ad}_\mcE])$ may be represented by $v$ as above such that $d \pi (a_{\alpha \beta})=0$ for any  pair $(\alpha, \beta)$.)
The $R_\epsilon$-schemes $U_{\alpha, \epsilon}$ (for various $\alpha \in I$) may be glued together by means of the isomorphisms
\begin{align} \label{W240}
\tau^v_{X, \alpha\beta} : =\mr{id}_{U_{\alpha \beta, \epsilon}} +\epsilon \cdot  d \pi (a_{\alpha \beta}) : 
U_{\beta, \epsilon} |_{U_{\alpha \beta}} \isom  U_{\alpha, \epsilon} |_{U_{\alpha \beta}}
\end{align}
 (for $(\alpha, \beta) \in I_2$).
 The resulting $R_\epsilon$-scheme, which we denote by  $X_\epsilon^v$,
 specifies  the deformation corresponding to $\pi^\mbH (v)$  via (\ref{W102}). 
 Moreover,  the flat $G$-torsors  
 $(\mcE_\epsilon |_{U_\alpha}, \nabla_{\mcE, \epsilon} |_{U_\alpha}+ \epsilon \cdot b_\alpha)$
 may be glued together by means of the isomorphisms
 \begin{align}
\tau^v_{\mcE, \alpha \beta} := \mr{id}_{\mcE_\epsilon |_{U_{\alpha\beta}}} + \epsilon \cdot a_{\alpha \beta}
 :
  (\mcE_\epsilon |_{U_\beta}, \nabla_{\mcE, \epsilon} |_{U_\beta}+ \epsilon  \cdot b_\beta) |_{U_{\alpha \beta}} 
\isom 
(\mcE_\epsilon |_{U_\alpha}, \nabla_{\mcE, \epsilon} |_{U_\alpha}+ \epsilon \cdot  b_\alpha) |_{U_{\alpha \beta}}
\end{align}
over $\tau^v_{X, \alpha \beta}$ (for $(\alpha, \beta) \in I_2$).
The data consisting of $X_\epsilon^v$  and the resulting flat $G$-torsor, which we denote by $(\mcE_\epsilon^v, \nabla_{\mcE, \epsilon}^v)$,
 specifies the deformation of $(X, \mcE, \nabla_\mcE)$ corresponding to $v$ via (\ref{W100}).
 The assignment $v \mapsto (X_\epsilon^v, \mcE^v_\epsilon, \nabla^v_{\mcE, \epsilon})$
obtained in this way gives the bijection (\ref{W100}).

\end{rema}
\vspace{3mm}



Since $\widetilde{\nabla}^\mr{ad}_\mcE \circ \nabla_\mcE =0$,
the pair of $\nabla_\mcE$ and the zero map $0 \migi \Omega_{X/R} \otimes \mr{Ad}(\mcE)$ induces a morphism of complexes $\mcT_{X/R}[0] \migi \mcK^\bullet [\widetilde{\nabla}^\mr{ad}_\mcE]$.
By applying this morphism to the functor $\mbH^1 (-)$, we obtain a split injection of (\ref{Ee1000}):
\begin{align} \label{W230}
\nabla_\mcE^\mbH  : H^1 (X, \mcT_{X/R}) \migiincl \mbH^1 (\mcK^\bullet [\widetilde{\nabla}^\mr{ad}_\mcE]).
\end{align}



\vspace{3mm}
\bpr  \label{W22}\leavevmode\\
 \ \ \
Let us take an element $v \in \mbH^1 (\mcK^\bullet [\widetilde{\nabla}^\mr{ad}_\mcE])$
 and denote by 
$(X_\epsilon^{v}, \mcE_\epsilon^v, \nabla_{\mcE, \epsilon}^v)$
 the deformation   over $R_\epsilon$  
 of $(X, \mcE, \nabla_{\mcE})$
 determined by $v$ via (\ref{W100}).
(In particular, $X_\epsilon^{v}$ is the deformation of $X$ determined by  $\pi^\mbH (v) \in  H^1 (X, \mcT_{X/R})$ via (\ref{W102}).)
Also, let us take an element $s \in J H^1 (X, \mcT_{X/R})$ (where we recall that $J$ is a sub-PD ideal of $I$)
  and denote by
  $X_\epsilon^{v, +s}$ the deformation of $X$ determined  by $\pi^\mbH (v) +s$.
  (Notice that the reductions modulo $J$ of $X_\epsilon^v$ and $X_\epsilon^{v, +s}$ are isomorphic.)
  Then, the deformation
  \begin{align} 
  (X_\epsilon^{v, +s}, \lambda_{X_\epsilon^{v, +s}} (\mcE^v_\epsilon), \lambda_{X_\epsilon^{v, +s}} (\nabla_{\mcE, \epsilon}^v))
  \end{align}
  (cf. Corollary \ref{W21} for the definition of $\lambda_{(-)}(-)$) corresponds to $v + \nabla^\mbH_\mcE (s) \in \mbH^1 (\mcK^\bullet [\widetilde{\nabla}^\mr{ad}_\mcE])$ via (\ref{W100}).
\epr
\begin{proof}
Let us fix an affine open covering $\mfU := \{ U_\alpha \}_{\alpha \in I}$ of $X$ and take a representative 
$(\{ a_{\alpha\beta} \}_{\alpha, \beta}, \{ b_{\alpha} \}_\alpha)$ of the class
 $v \in \mbH^1 (\mcK^\bullet [\widetilde{\nabla}^\mr{ad}_\mcE])$ as displayed in (\ref{W105}).
Also,  $s$ may be represented by a \v{C}ech $1$-cocycle $\{ s_{\alpha \beta}\}_{\alpha, \beta} \in \check{C} (\mfU, \mcT_{X/R})$.
In the following, we shall use the  notation  in Remark \ref{W101} and apply the discussion there  to the present  $v$.
In particular, 
$X_\epsilon^v$  denotes the curve over $R_\epsilon$  obtained by gluing together $U_{\alpha, \epsilon}$'s by means of the isomorphisms $\tau^v_{X, \alpha\beta}$'s  (cf. (\ref{W240})).
Since $\mr{id}_{U_{\alpha \beta, \epsilon}}  + \epsilon \cdot (s_{\alpha \beta} + d\pi  (a_{\alpha \beta})) = (\mr{id}_{U_{\alpha \beta, \epsilon}} + \epsilon \cdot  s_{\alpha \beta}) \circ \tau_{X, \alpha \beta}^{v}$,
the deformation $X^{v, +s}_\epsilon$ 
 may be obtained by gluing together $U_{\alpha, \epsilon}$'s 
  by means of
the isomorphisms  $(\mr{id}_{U_{\alpha \beta, \epsilon}} + \epsilon \cdot  s_{\alpha \beta}) \circ \tau_{X, \alpha \beta}^{v}$ (for various $(\alpha, \beta) \in I_2$).
Now,
let  $e : \mr{pr}^*_2(\mcE_\epsilon^v) \isom \mr{pr}^*_1 (\mcE_\epsilon^v)$ be 
 the HPD stratification on $\mcE_\epsilon^v$ corresponding to $\nabla^v_{\mcE, \epsilon}$ via the equivalence of categories between (ii) and (v) in Theorem \ref{T043}.
 If $\iota^v_{\alpha \beta} : D_{U_{\alpha\beta, \epsilon}/R_\epsilon} (1) \migi D_{X^v_{\epsilon}/R_\epsilon} (1)$ (for each $(\alpha, \beta)\in I_2$)
 denotes the morphism arising from the natural open immersion $U_{\alpha \beta, \epsilon} \migiincl X^v_{\epsilon}$,
 then the following equality holds:
 \begin{align}
  (\iota_{\alpha \beta}^v   \circ (\mr{id}_{U_{\alpha\beta, \epsilon}}, \mr{id}_{U_{\alpha\beta, \epsilon}} + \epsilon \cdot s_{\alpha \beta}) )^*(e) & = \mr{id}_{\mcE_\epsilon |_{U_{\alpha\beta}}} + (\nabla_{\mcE, \epsilon}|_{U_{\alpha}} + \epsilon \cdot b_\alpha)(\epsilon \cdot s_{\alpha \beta}) \\
  & =\mr{id}_{\mcE_\epsilon |_{U_{\alpha\beta}}} + \epsilon \cdot  \nabla_{\mcE} (s_{\alpha \beta}), \notag
 \end{align}
 where $(\mr{id}_{U_{\alpha\beta, \epsilon}}, \mr{id}_{U_{\alpha\beta, \epsilon}} + \epsilon \cdot s_{\alpha \beta})$ denotes the unique morphism $U_{\alpha\beta, \epsilon} \migi D_{U_{\alpha\beta, \epsilon}/R_\epsilon} (1)$ whose composites with the first and the second  projections  $D_{U_{\alpha\beta, \epsilon}/R_\epsilon} (1) \migi U_{\alpha\beta, \epsilon}$ coincide with 
 $\mr{id}_{U_{\alpha\beta, \epsilon}}$ and $\mr{id}_{U_{\alpha\beta, \epsilon}} + \epsilon \cdot s_{\alpha \beta}$ respectively.
The following square diagram is commutative:
\begin{align}
\begin{CD}
\mcE_{\epsilon} |_{U_{\alpha \beta}}@> (\mr{id} + \epsilon \cdot \nabla_\mcE (s_{\alpha \beta})) \circ  \tau_{\mcE, \alpha \beta}^v > \sim > \mcE_{\epsilon} |_{U_{\alpha \beta}}
\\
@VVV @VVV
\\
U_{\beta, \epsilon} |_{U_{\alpha \beta}} @> \sim > (\mr{id} + \epsilon \cdot  s_{\alpha \beta}) \circ  \tau_{X, \alpha \beta}^{v} > U_{\alpha, \epsilon} |_{U_{\alpha \beta}}, 
\end{CD}
\end{align}
where the both sides of the vertical arrows denote the natural projections.
It follows from the definition of 
$\lambda_{(-)}(-)$
that 
$(\lambda_{X_\epsilon^{v, +s}} (\mcE^v_\epsilon), \lambda_{X_\epsilon^{v, +s}} (\nabla_{\mcE, \epsilon}^v))$ may be obtained by gluing together 
$(\mcE_\epsilon |_{U_\alpha}, \nabla_{\mcE, \epsilon} |_{U_\alpha} + \epsilon \cdot b_\alpha)$'s by means of the isomorphisms 
\begin{align}
(\mr{id}_{\mcE_\epsilon |_{U_{\alpha \beta}}} + \epsilon \cdot \nabla_\mcE (s_{\alpha \beta})) \circ  \tau_{\mcE, \alpha \beta}^v \ \left(=\mr{id}_{\mcE_\epsilon |_{U_{\alpha \beta}}}+\epsilon \cdot (a_{\alpha \beta} + \nabla_\mcE (s_{\alpha\beta}))\right)
\end{align}
for various $(\alpha, \beta) \in I_2$.
 Hence,  the element of $\mbH^1 (\mcK^\bullet [\widetilde{\nabla}^\mr{ad}_\mcE])$ corresponding to the deformation $(X_\epsilon^{v, +s}, \lambda_{X_\epsilon^{v, +s}} (\mcE^v_\epsilon), \lambda_{X_\epsilon^{v, +s}} (\nabla_{\mcE, \epsilon}^v))$ may be represented by 
 the collection of data
 \begin{align}
 (\{ a_{\alpha \beta} + \nabla_\mcE (s_{\alpha \beta}) \}_{\alpha, \beta}, \{ b_\alpha \}_\alpha),
 \end{align}
which specifies $v + \nabla^\mbH_\mcE (s)$.
This completes the proof of the assertion.
\end{proof}
\vspace{3mm}


\vspace{5mm}
\subsection{Deformation space of  indigenous bundles} \label{SS1083}
\leavevmode\\ \vspace{-4mm}

Suppose further that 
$G = \mr{PGL}_2$  and $(\mcE, \nabla_\mcE)$ forms an  indigenous bundle on the curve $X/R$.
Denote by $\pi_B : \mcE_B \migi X$ the Hodge reduction of $(\mcE, \nabla_\mcE)$ and by
\begin{align} \label{W50}
\widetilde{\nabla}^\mr{ad}_{\mcE_B} : \widetilde{\mcT}_{\mcE_B/R} \migi \Omega_{X/R} \otimes \mr{Ad}(\mcE)
\end{align}
 the morphism obtained  from $\widetilde{\nabla}^\mr{ad}_\mcE$ by restricting its domain. 
Then, we obtain the following morphism of short exact sequences:
\begin{align} \label{W51}
\begin{CD}
0 @>>> \widetilde{\mcT}_{\mcE_B/R} @>\mr{incl.}>> \widetilde{\mcT}_{\mcE/R}  @> d\breve{\pi}   >> \mcT_{X/R} @>>> 0
\\
@. @VVV @VVV @VVV @.
\\
0 @>>> \mr{Ad}(\mcE)\otimes \Omega_{X/R}@>> \mr{id} > \mr{Ad}(\mcE) \otimes \Omega_{X/R} @>>> 0 @>>> 0,
\end{CD}
\end{align}
where 
$d\breve{\pi}$
 denotes the composite of the quotient $\widetilde{\mcT}_{\mcE/R} \migisurj \widetilde{\mcT}_{\mcE/R}/\widetilde{\mcT}_{\mcE_B/R}$ and the inverse of $\mr{KS}_{(\mcE, \nabla_\mcE)} : \mcT_{X/R} \isom \widetilde{\mcT}_{\mcE/R}/\widetilde{\mcT}_{\mcE_B/R}$ (cf. (\ref{G803})).
Since $H^0 (X, \mcT_{X/R}) =0$,  this morphism of sequences  induces the following short exact sequence:
\begin{align} \label{W52}
0 \longmigi \mbH^1 (\mcK^\bullet [\widetilde{\nabla}^\mr{ad}_{\mcE_B/R}]) \stackrel{\breve{\iota}^\mbH}{\longmigi} 
\mbH^1 (\mcK^\bullet [\widetilde{\nabla}^\mr{ad}_{\mcE/R}])
\stackrel{\breve{\pi}^\mbH}{\longmigi} 
H^1 (X, \mcT_{X/R})
\longmigi 
0.
\end{align}
If we consider $ \mbH^1 (\mcK^\bullet [\widetilde{\nabla}^\mr{ad}_{\mcE_B/R}])$
as a submodule of $ \mbH^1 (\mcK^\bullet [\widetilde{\nabla}^\mr{ad}_{\mcE/R}])$ via $\breve{\iota}^\mbH$,
then
the elements of $ \mbH^1 (\mcK^\bullet [\widetilde{\nabla}^\mr{ad}_{\mcE_B/R}])$ 
classifies the deformations 
equipped with a deformation of the $B$-reduction $\mcE_B$.
This implies that
if 
\begin{align} \label{W53}
\mr{Def}^{\mr{ind}}_{R_\epsilon} (X, \mcE, \nabla_\mcE)
\end{align}
the set of isomorphism classes of deformations over $R_\epsilon$ of $(X, \mcE, \nabla_\mcE)$ (as a data consisting of a curve and an indigenous bundle on it),
then the bijection (\ref{W100}) restricts to a bijection
\begin{align} \label{D003}
\mbH^1 (\mcK^\bullet [\widetilde{\nabla}^\mr{ad}_{\mcE_B}]) \isom \mr{Def}^{\mr{ind}}_{R_\epsilon} (X, \mcE, \nabla_\mcE).
\end{align}

 Denote by
 \begin{align} \label{W54}
 \eta : \widetilde{\mcT}_{\mcE/R} \migi  \widetilde{\mcT}_{\mcE/R}
 \end{align}
 the $\mcO_X$-linear  endomorphism of $\widetilde{\mcT}_{\mcE/R}$ given by 
 $s \mapsto s - (\nabla_\mcE \circ d \pi)(s)$ for any local section $s \in \widetilde{\mcT}_{\mcE/R}$.
One verifies that $\eta$ is an isomorphism and its inverse may be given by $s - (\nabla_\mcE \circ  d \breve{\pi})(s)$.
Since the equality $\widetilde{\nabla}_\mcE^\mr{ad} \circ \eta = \widetilde{\nabla}_\mcE^\mr{ad}$ holds,
 the pair of morphisms $(\eta, \mr{id}_{\Omega_{X/R}\otimes \mr{Ad}(\mcE)})$
   specifies
an automorphism of the complex $\mcK^\bullet [\widetilde{\nabla}^\mr{ad}_\mcE]$.
Moreover, one verify that this automorphism restricts to  an isomorphism $\mcK^\bullet [\widetilde{\nabla}^\mr{ad}_{\mcE_B}] \isom \mcK^\bullet [\nabla^\mr{ad}_{\mcE}]$.
By applying the functor $\mbH^1 (-)$ to  these isomorphisms of complexes, we  obtain 
the following square diagram:
\begin{align} \label{W356}
\begin{CD}
 \mbH^1 (\mcK^\bullet [\widetilde{\nabla}^\mr{ad}_{\mcE_B}]) @> \eta^\mr{ad}_B>\sim>\mbH^1 (\mcK^\bullet [\nabla^\mr{ad}_\mcE])
\\
@V \breve{\iota}^\mbH VV @VV \iota^\mbH V
\\
\mbH^1 (\mcK^\bullet [\widetilde{\nabla}^\mr{ad}_\mcE]) @> \sim >\eta^\mr{ad} >  \mbH^1 (\mcK^\bullet [\widetilde{\nabla}^\mr{ad}_\mcE]).
\end{CD}
\end{align}
 In particular, by restricting $\eta_B^\mr{ad}$, we obtain an isomorphism 
 \begin{align} \label{W430}
 \eta^\mr{ad}_{B, J} : J\mbH^1 (\mcK^\bullet [\widetilde{\nabla}^\mr{ad}_{\mcE_B}]) \isom J\mbH^1 (\mcK^\bullet [\nabla^\mr{ad}_\mcE]).
 \end{align}
 
 Now, let us introduce some notation to describe the statement of Proposition \ref{W23} below.
 Fix an element $u$  of
  $\mbH^1 (\mcK^\bullet [\widetilde{\nabla}^\mr{ad}_{\mcE_B}])\otimes \overline{R}$, where  $\overline{R} := R/J$.
 Denote by
  \begin{align} \label{W401}
  \mr{Def}_{R_\epsilon}^\mr{ind} (X, \mcE, \nabla_\mcE)_u
  \end{align}
 the subset of  $ \mr{Def}^\mr{ind}_{R_\epsilon}(X, \mcE, \nabla_\mcE)$
 consisting of deformations whose reductions modulo $J$ correspond to  $u$ via  the bijection (\ref{D003}).
  Notice that the $\mbH^1 (\mcK^\bullet [\widetilde{\nabla}^\mr{ad}_{\mcE_B}])$-action  on   $ \mr{Def}^\mr{ind}_{R_\epsilon}(X, \mcE, \nabla_\mcE)$ arising  from 
   (\ref{D003})  gives the structure of  $J \mbH^1 (\mcK^\bullet [\widetilde{\nabla}^\mr{ad}_{\mcE_B}])$-torsor on 
$\mr{Def}_{R_\epsilon}^\mr{ind} (X, \mcE, \nabla_\mcE)_u$.
 Also, let us fix an element  $s$ of $H^1 (X, \mcT_{X/R})$ whose image in 
$H^1 (X, \mcT_{X/R}) \otimes \overline{R}$
coincides with $(\pi^\mbH \otimes \mr{id}_{\overline{R}})(u)$. 
 Denote by
 \begin{align} \label{W402}
 \mr{Def}_{R_\epsilon} (X, \mcE, \nabla_\mcE)_{u, s}
 \end{align}
  the subset of  $\mr{Def}_{R_\epsilon} (X, \mcE, \nabla_\mcE)$ consisting of deformations  which are contained   in 
  $(\pi^\mbH)^{-1}(s)$ via (\ref{W100}) and whose reductions modulo $J$ correspond to $u$.
   By means of   (\ref{W100}), it has canonically  the structure of $J \mbH^1 (\mcK^\bullet [\nabla^\mr{ad}_{\mcE}])$-torsor.
Then, the following assertion holds.
  

\vspace{3mm}
\bpr \label{W23}\leavevmode\\
 \ \ \
 Let us keep the notation in Proposition \ref{W22}.
 Suppose that $(\mcE, \nabla_\mcE)$ is quasi-nilpotent.
Then, the assignment 
\begin{align} \label{W500}
(X^v_\epsilon, \mcE^v_\epsilon, \nabla^v_{\mcE, \epsilon}) \mapsto
 (X_\epsilon^{v, +s \odot v}, \lambda_{X_\epsilon^{v, +s\odot v}} (\mcE^v_\epsilon), \lambda_{X_\epsilon^{v, + s \odot v}} (\nabla^v_{\mcE, \epsilon})),
\end{align}
where $s \odot v :=s-\pi^\mbH (v)$,  defines  a bijection  
\begin{align} \label{W59}
 \mr{Def}_{R_\epsilon}^\mr{ind} (X, \mcE, \nabla_\mcE)_u
  \isom
   \mr{Def}_{R_\epsilon} (X, \mcE, \nabla_\mcE)_{u, s},
\end{align}
that is compatible with the respective torsor  structures via $\eta^\mr{ad}_{B, J}$ (cf. (\ref{W430})).
  \epr
\begin{proof}
By  definition,  the assignment (\ref{W500}) is compatible with the respective torsor structures on 
 $\mr{Def}_{R_\epsilon}^\mr{ind} (X, \mcE, \nabla_\mcE)_u$ and $\mr{Def}_{R_\epsilon} (X, \mcE, \nabla_\mcE)_{u, s}$.
Thus, the assertion follows from the fact that $\eta^\mr{ad}_{B, J}$ is bijective.
\end{proof}
\vspace{3mm}

By Proposition \ref{W23}, we obtain the following  corollary, which was used in the proof of Theorem \ref{P0438}.

\vspace{3mm}
\bco \label{W150} \leavevmode\\
 \ \ \
Denote by $X_\epsilon^s$ deformation over $R_\epsilon$ corresponding to $s$ and by  $(\overline{X}^u_\epsilon, \overline{\mcE}^u_{\epsilon}, \overline{\nabla}^u_{\mcE, \epsilon})$
the deformation over $\overline{R}_\epsilon \ (:= R_\epsilon \times_R \overline{R})$ of $(X, \mcE, \nabla_\mcE)$ corresponding to $u$.
(Hence, the reduction modulo $J$ of $X_\epsilon^s$ is isomorphic to $\overline{X}^u_\epsilon$, and $(\overline{\mcE}^u_{\epsilon}, \overline{\nabla}^u_{\mcE, \epsilon})$ forms an indigenous bundle on $\overline{X}^u_\epsilon$.)
Now, suppose that $(\mcE, \nabla_\mcE)$ is quasi-nilpotent.
 Also, let $\mcE^\lozenge$ be  a crystal of $\mr{PGL}_2$-torsors
 on $\mr{Crys}(\overline{X}^{u}_\epsilon /R_\epsilon)$
such that the associated  crystal on $\mr{Crys}(\overline{X}^{u}_\epsilon /\overline{R}_\epsilon)$ corresponds to $(\overline{\mcE}^u_{\epsilon}, \overline{\nabla}^u_{\mcE, \epsilon})$ via the equivalence of categories between (i) and (v) in Theorem \ref{T043}.
Then, there exists  a unique (up to isomorphism) deformation over $X^s_\epsilon$ of $(\overline{\mcE}^u_{\epsilon}, \overline{\nabla}^u_{\mcE, \epsilon})$ which forms an indigenous bundle and corresponds to   $\mcE^\lozenge$ (via the equivalence between (i) and (v) as above).
  \eco

\end{document}